\documentclass[opre, nonblindrev]{informs3}
\OneAndAHalfSpacedXI

\usepackage{natbib}
\bibpunct[, ]{(}{)}{,}{a}{}{,}%
\def\bibfont{\small}%
\def\bibsep{\smallskipamount}%

\TheoremsNumberedThrough
\EquationsNumberedThrough

\usepackage{amsmath, comment, tikz, graphicx, algorithm, multirow, multicol, subcaption, bbm, epsfig, enumerate, amsfonts, enumitem, subeqnarray, url, color, array, xcolor, booktabs, endnotes, float, mathpazo,pdflscape}

\usepackage[figuresright]{rotating}
\usepackage[noend]{algpseudocode}
\usetikzlibrary{decorations.pathreplacing}
\usepackage[skip=2pt, font=footnotesize, labelfont=bf, labelsep=period]{caption}

\usepackage{pgfplots}
    \pgfplotsset{
        compat=1.9,
    }
\usepackage[backref = false, bookmarks, breaklinks=true, colorlinks = true, plainpages = false, citecolor = DarkBlue , urlcolor = DarkBlue, filecolor = DarkBlue, linkcolor = DarkBlue]{hyperref}
\definecolor{DarkBlue}{rgb}{0,0.08,0.45}

\usetikzlibrary{positioning, arrows, snakes, backgrounds}
\usetikzlibrary{calc}

\let\footnote=\endnote

\renewcommand{\theARTICLETOP}{}
\renewcommand{\theRRHFirstLine}{}
\renewcommand{\theRRHSecondLine}{}
\renewcommand{\theLRHFirstLine}{}
\renewcommand{\theLRHSecondLine}{}

\algnewcommand\algorithmicinitialize{\textbf{Initialize:}}
\algnewcommand\Initialize{\item[\algorithmicinitialize]}
\algnewcommand\algorithmicrepeati{\textbf{Repeat:}}
\algnewcommand\Repeati{\item[\algorithmicrepeati]}
\algnewcommand\algorithmicutlcon{\textbf{Until Convergence.}}
\algnewcommand\Utlcon{\item[\algorithmicutlcon]}

\allowdisplaybreaks

\begin{document}

\TITLE{Optimized Dimensionality Reduction for Moment-based Distributionally Robust Optimization}

\ARTICLEAUTHORS{ \vspace{-0.2cm}
\AUTHOR{Shiyi Jiang$^{\textup a}$, Jianqiang Cheng$^{\textup b}$, Kai Pan$^{{\textup a}}$, Zuo-Jun Max Shen$^{\textup{c,d}}$} 
\vspace{0.05cm}

\AFF{$^{\textup a}$Faculty of Business, The Hong Kong Polytechnic University, Kowloon, Hong Kong}
\AFF{$^{\textup b}$College of Engineering, University of Arizona, Tucson, AZ 85721, USA}
\AFF{$^{\textup c}$College of Engineering, University of California, Berkeley, California 94720, USA}
\AFF{$^{\textup d}$Faculty of Engineering and Faculty of Business and Economics, The University of Hong Kong, Hong Kong}
\vspace{0.05cm}

\AFF{Contact: \EMAIL{shiyi.jiang@connect.polyu.hk} (SJ), \EMAIL{jqcheng@arizona.edu} (JC), \\ \EMAIL{kai.pan@polyu.edu.hk} (KP), \EMAIL{maxshen@berkeley.edu} (Z-JMS)} \vspace{-0.4cm}
}

\ABSTRACT{Moment-based distributionally robust optimization (DRO) provides an optimization framework to integrate statistical information with traditional optimization approaches. 
Under this framework, one assumes that the underlying joint distribution of random parameters runs in a distributional ambiguity set constructed by moment information and makes decisions against the worst-case distribution within the set.
Although most moment-based DRO problems can be reformulated as semidefinite programming (SDP) problems that can be solved in polynomial time, solving high-dimensional SDPs is still time-consuming. 
Unlike existing approximation approaches that first reduce the dimensionality of random parameters and then solve the approximated SDPs, we propose an optimized dimensionality reduction (ODR) approach.
We first show that the ranks of the matrices in the SDP reformulations are small, by which we are then motivated to integrate the dimensionality reduction of random parameters with the subsequent optimization problems.
Such integration enables two outer and one inner approximations of the original problem, all of which are low-dimensional SDPs that can be solved efficiently, providing two lower bounds and one upper bound correspondingly.
More importantly, these approximations can theoretically achieve the optimal value of the original high-dimensional SDPs.
As these approximations are nonconvex SDPs, we develop modified Alternating Direction Method of Multipliers (ADMM) algorithms to solve them efficiently.
We demonstrate the effectiveness of our proposed ODR approach and algorithm in solving 
multiproduct newsvendor and production-transportation problems.
Numerical results show significant advantages of our approach on the computational time and solution quality over the three best possible benchmark approaches.
Our approach can obtain an optimal or near-optimal (mostly within $0.1\%$) solution and reduce the computational time by up to three orders of magnitude.}




\maketitle


\section{Introduction} \label{Sec:Introduction}

Distributionally robust optimization (DRO) is a modeling framework that integrates statistical information with traditional optimization methods \citep{scarf1958min, delage2010distributionally}. 
Under this framework, one assumes that the underlying joint distribution of random parameters runs in a distributional ambiguity set inferred from given data or prior belief
and then optimizes their decisions against the worst-case distribution within the set.


To solve different applications, researchers study the DRO under various distributional ambiguity sets. 
The ambiguity set plays a crucial role in connecting statistical information and optimization modeling, providing a flexible framework for modeling uncertainties and incorporating partial information of random parameters into the model, such as information from historical data and prior belief. 
Moreover, the performance of DRO depends significantly on the distributional ambiguity set \citep{mohajerin2018data}. This paper focuses on moment-based ambiguity sets, which include all distributions satisfying certain moment constraints. Examples of such constraints include restricting all distributions to have the exact mean and covariance matrix \citep{scarf1958min}, bounding the first and second moments \citep{delage2010distributionally}, or placing the first and second moments in a convex set \citep{ghaoui2003worst}. Moment-based DRO has been extensively studied because it may be more tractable than other ambiguity sets and has a wide range of applications in industry, including but not limited to newsvendor problems \citep{gallego1993distribution,yue2006expected,natarajan2018asymmetry}, portfolio optimization problems \citep{ghaoui2003worst,goldfarb2003robust, zymler2013worst,rujeerapaiboon2016robust,li2018closed,lotfi2018robust}, knapsack problems \citep{cheng2014distributionally}, transportation problems \citep{zhang2017robust,ghosal2020distributionally}, reward-risk ratio optimization \citep{liu2017distributionally}, scheduling problems \citep{shehadeh2020distributionally}, and machine learning \citep{lanckriet2002robust, farnia2016minimax}.

Moment-based DRO model can be reformulated as a semi-infinite program \citep{xu2018distributionally}, which is generally intractable. 
Three approaches are mainly used to solve such a reformulation: 
(i) the cutting plane/surface method \citep{gotoh2002bounding,mehrotra2014cutting}, by which a solution is first obtained by considering a subset of the distributional ambiguity set and cuts are then added iteratively until converging to an optimal solution;
(ii) the dual method \citep{delage2010distributionally,bertsimas2019adaptive}, by which the inner optimization problem (e.g., a minimization problem) is dualized and integrated with the outer optimization problem (e.g., a maximization problem);
(iii) the analytical method \citep{scarf1958min, popescu2007robust}, by which the worst-case distribution is obtained and its properties are analyzed.
Among these methods, the dual method is the most popular. 
Most literature focuses on convex reformulations of different moment-based DRO problems, mainly including second-order cone programming (SOCP) \citep{ghaoui2003worst, lotfi2018robust, goldfarb2003robust} and semidefinite programming (SDP) \citep{ghaoui2003worst, delage2010distributionally, cheng2014distributionally}.



While SOCPs can be solved efficiently, theoretically efficient algorithms (e.g., the interior-point methods) to solve SDPs impose substantial demands on computational time and memory resources \citep{burer2003nonlinear}, particularly when addressing high-dimensional SDPs.
Widely adopted commercial solvers (e.g., Mosek) exhibit prohibitively long computational times when solving high-dimensional SDPs, and the computational burden escalates considerably even as the problem dimension increases gradually (see our numerical results in Section \ref{Sec:experiments}).
Thus, it is of great interest to study efficient algorithms for solving SDPs in the context of moment-based DRO. 
Besides the generic methods (e.g., the interior point methods), two types of algorithms can speed up solving SDP reformulations of moment-based DRO: low-rank SDP algorithms and dimensionality reduction methods.
First, 
some studies develop efficient low-rank algorithms by exploiting the low-rank properties of SDP constraints \citep{burer2003nonlinear}. 
These algorithms rarely have theoretical guarantees but are practically efficient. 
Specifically, the existing studies may reformulate convex SDPs as non-convex problems and subsequently develop efficient algorithms to deliver high-quality solutions within significantly reduced time frames \citep{lemon2016low}.
Second, dimensionality reduction techniques stem from the field of statistics to represent the data with the important information while omitting the trivial one.
In the context of moment-based DRO, such techniques can be extended to reduce the dimension of random parameters and approximate the high-dimensional SDP reformulations using low-dimensional SDPs \citep{cheng2018distributionally, cheramin2022computationally}, thereby reducing the computational time significantly.

However, both the general SDP algorithms and existing dimensionality reduction methods may not perform well in the context of moment-based DRO. 
The general SDP algorithm aims to solve more general SDPs and may fail to consider the specific structure of the moment-based DRO models. 
The existing dimensionality reduction methods fail to consider the subsequent optimization problems when reducing the dimensionality space. 
For example, \cite{cheng2018distributionally} and \cite{cheramin2022computationally} first use the PCA to choose the random parameters corresponding to the largest eigenvalues and then solve the low-dimensional SDP problem with the chosen random parameters. 
Such a sequential process may not provide an optimal solution of the original problem because the aim of leveraging data is to reduce the dimensionality space by focusing on only the statistical information, rather than optimizing the subsequent SDP problems. 
Therefore, in this paper, we integrate the dimensionality reduction with subsequent SDP problems, which leads to an \textit{optimized dimensionality reduction (ODR) approach for moment-based DRO}.
This idea echoes the recently emerging framework that integrates machine learning with decision-making \citep{bertsimas2020predictive, elmachtoub2022smart}. 
We summarize our contributions as follows: 
\begin{enumerate}
    \item We prove the low-rank property of SDP reformulations of moment-based DRO problems. Specifically, we show that the ranks of matrices in SDP reformulations are less than the number of SDP constraints plus one.
    \item Different from the PCA approximation approaches \citep{cheng2018distributionally, cheramin2022computationally} that first reduce the dimensionality and then solve approximation problems, we integrate the dimensionality reduction with the subsequent optimization problems and provide an optimized dimensionality reduction approach.
    \item With the ODR approach, we develop two outer and one inner approximations for the original problem, leading to three low-dimensional SDP problems that can be solved efficiently. 
    More importantly, these low-dimensional approximations can achieve the optimal value of the original high-dimensional SDP. 
    \item The low-dimensional SDP problems are nonconvex with bilinear terms and we develop modified Alternating Direction Method of Multipliers (ADMM) algorithms to solve them efficiently. 
    We apply the ODR approach and ADMM algorithms to solve multiproduct newsvendor and production-transportation problems.
    We compare our ODR approach with three benchmark approaches: the Mosek solver, low-rank algorithm \citep{burer2003nonlinear}, and PCA approximations \citep{cheramin2022computationally}. 
    The results demonstrate that our ODR approach significantly outperforms them in terms of computational time and solution quality. 
    Our approach can obtain an optimal or near-optimal (mostly within $0.1\%$) solution and reduce the computational time by up to three orders of magnitude.
    More importantly, our approach is not sensitive to the dimension $m$ of random parameters, while the benchmark approaches perform much worse when $m$ is larger. 
\end{enumerate}

The remainder of this paper is organized as follows. 
In Section \ref{Sec:liteReview}, we review related literature. 
In Section \ref{Sec:introdMomDRO}, we present the SDP reformulation of moment-based DRO problems. 
In Section \ref{Sec:ODRapproach}, we prove the low-rank property of the SDP reformulation and propose the first outer approximation under the ODR approach, leading to the lower bound for the original problem. 
In Sections \ref{sec:upper-bound} and \ref{sec:lower-bound-2}, we provide the upper bound and the second lower bound for the original problem, respectively.
In Section \ref{Sec:experiments}, we perform extensive numerical experiments on multiproduct newsvendor and production-transportation problems.
Section \ref{Sec:conclusion} concludes the paper. 
All proofs are presented in the Appendix if not specified.

\medskip


\noindent \textbf{Notation}\\
We use non-bold symbols to denote scalar values, e.g., $s$ and $\gamma_1$, and bold symbols to denote vectors, e.g., $\mathbf{x}=\left(x_1,\ldots,x_n\right)^{\top}$ and ${\mathbf{q}}$. 
Similarly, matrices are represented by bold capital symbols, e.g., $\mathbf{A}$ and $\boldsymbol{\Sigma}$, and the size of a matrix is indicated by $r \times c$, where $r$ and $c$ indicate the numbers of rows and columns, respectively. Italic subscripts indicate indices, e.g., $S_k$, while non-italic ones represent simplified specifications, e.g., $\mathbf{Q} _{ \textup{r} }$. We use $ \mathbb{E}_{\mathbb{P}} \left[ \cdot \right]$ to represent expectation over distribution $\mathbb{P}$ and use $"\bullet"$ to denote the inner product defined by $\mathbf{A} \bullet \mathbf{B} =
\sum_{i,j} A_{ij}B_{ij}$, where $\mathbf{A}$ and $\mathbf{B}$ are two conformal matrices. 
For any matrix $\mathbf{M}$, we use $\mathbf{M} \succeq 0$ (resp. $\mathbf{M}\succ 0$) to indicate that it is positive semi-definite (PSD) (resp. positive definite). 
Symbols $\left \| \cdot \right \|_1$ and $\left \| \cdot \right \|_2$ denote L1-Norm and L2-Norm, respectively. 
For any integer number $n \geq 1$, we use $[n]$ to denote the set $\left\{1, 2, \ldots, n\right\}$. 
The identity matrix of size $m$ is denoted by $\mathbf{I}_m$. Symbols $\boldsymbol{0}_m$ and $\boldsymbol{0}_{ r \times c }$ represent a zero vector of size $m$ and a zero matrix of size $r \times c$, respectively. 
Symbols $\boldsymbol{1}_m$ and $\boldsymbol{1}_{ r \times c }$ represent a one vector of size $m$ and a one matrix of size $r \times c$, respectively. 
We use $\mathbb{1}( \cdot )$ to denote the indicator function, 
which takes $1$ if all the conditions encompassed in $( \cdot )$ are satisfied and 
takes $0$ otherwise.

\section{Literature Review} \label{Sec:liteReview} 
We review related literature from four streams: moment-based DRO, dimensionality reduction, low-rank SDP algorithms, and the integration of machine learning with decision-making.

\subsection{Moment-based DRO}


Extensive studies provide the theories and applications of moment-based DRO (see \citealt{rahimian2019distributionally} and  \citealt{lin2022distributionally} for detailed review). 
\cite{scarf1958min} and \cite{gallego1993distribution} are among the first to introduce the moment-based DRO framework, under which they consider fixed mean and variance of random parameters and analytically obtain the worst-case distribution in the ambiguity set, thereby obtaining the analytical optimal solution.
\cite{ghaoui2003worst} study distributionally robust portfolio optimization problem, where the value at risk (VaR) is minimized against the worst-case distribution in the ambiguity set.
The problem is reformulated into an SOCP if the mean and covariance matrix are given and into an SDP if they belong to bounded convex ambiguity sets. 
There are a series of follow-up studies \citep{delage2010distributionally, li2018closed, lotfi2018robust, zymler2013worst}.
Among them, \cite{delage2010distributionally} consider a general moment-based DRO framework with an ambiguity set considering the information of support, mean, and covariance matrix, leading to an SDP reformulation. 
We consider the same framework in this paper.

\subsection{Dimensionality Reduction}

To efficiently represent high-dimensional data, dimensionality reduction techniques are proposed in the literature to maintain the important information from the data and omit the trivial information. 
Principal component analysis (PCA) provides a high-quality representation of the data with as much information as possible by maintaining the random variables with the largest eigenvalues \citep{abdi2010principal}. 
Several variants of PCA are further studied in the literature: robust PCA \citep{candes2011robust}, scaled PCA \citep{huang2022scaled}, and sparse PCA \citep{dey2022using}. 
They are largely applied in the field of machine learning: independent component analysis \citep{comon1994independent}, latent semantic analysis \citep{landauer1998introduction}, and locality preserving projections \citep{he2003locality}. 
As the PCA may still lose useful information, sufficient dimensionality reduction is proposed to represent the information from the data using a linear combination of the original random variables. 
A series of specific methods are studied in machine learning, such as the sliced inverse regression \citep{li1991sliced}, sliced average variance estimation \citep{li2007asymptotics}, principal Hessian direction \citep{li1992principal}, and minimum average variance estimation \citep{yin2011sufficient}. 



These dimensionality reduction techniques are rarely used to support decision-making in mathematical optimization. 
\cite{wang2006effects} and \cite{wang2011quasi} use the Brownian bridge, PCA, and Quasi-Monte Carlo method to reduce the dimensionality of high-dimensional integrals in finance problems. 
\cite{cheng2018distributionally} and \cite{cheramin2022computationally} are among the first to reduce the dimensionality space of random variables that are modeled in a moment-based DRO framework. 
They adopt the PCA to maintain the important random variables in the ambiguity set and reduce the size of the subsequent SDP reformulation.
However, they consider only the statistical information and fail to consider the structure information of the subsequent SDP reformulation when performing dimensionality reduction.
We resolve this issue in this paper \citep{jiang2023optimized}, which is
recently followed by \cite{he2023prescriptive}. 
\cite{he2023prescriptive} integrate the PCA with a subsequent stochastic program and provide a distributionally robust bound for the error between the objective values of the original and integrated problems. 
The integrated approach in \cite{he2023prescriptive} involves solving nonconvex and high-dimensional SDPs and may not reduce the error to zero, while our approach solves low-dimensional SDPs and can achieve the optimal value of the original moment-based DRO problem, thereby offering guidance on selecting the reduced dimension for practical applicability.

\subsection{Low-rank SDP Algorithms}

As the moment-based DRO model is reformulated as an SDP \citep{vandenberghe1996semidefinite} in this paper, SDP  algorithms are important to solve it. 
Commercial optimization solvers (e.g., Mosek and Gurobi) use the interior point method to solve SDPs.
Although this method can converge very fast \citep{helmberg2002semidefinite}, its computation is very expensive. 
More importantly, as a general algorithm, it does not exploit useful structural properties of the SDP constraints.
To solve this issue, \cite{burer2003nonlinear} are among the first to propose low-rank algorithms \citep{lemon2016low} to solve general SDPs.
Specifically, they analyze the low-rank property of the SDP constraints and transform the convex SDP into a nonconvex optimization problem with a smaller size, which is further solved by augmented Lagrangian methods.
In addition, \cite{yurtsever2021scalable} consider trace-constrained SDPs and show that the SDP constraints are weakly constrained,
by which a low-rank approximation is proposed and efficiently solved.
Our paper solves an integrated optimization problem that incorporates both the dimensionality reduction of random parameters and SDP formulation, under which we also exploit the low-rank property of our formulation. 
Numerical results show that our proposed ODR approach performs better than the low-rank algorithm in \cite{burer2003nonlinear}.

\subsection{Integration of Machine Learning with Decision-making}

\cite{bertsimas2020predictive} summarize that many optimization problems have three primitives: (i) data on uncertain parameters, (ii) auxiliary data on associated covariates, and (iii) a structured optimization concerning decisions, constraints, and objective functions.
Traditional approaches first build machine learning models to perform parameter estimation and then solve the optimization problem with the estimated parameters, while a good prediction may not lead to a good decision.
Thus, \cite{bertsimas2020predictive}, \cite{bertsimas2022data}, and \cite{elmachtoub2022smart} integrate the parameter estimation with optimization problems. 
Similar ideas are reflected in early attempts in \cite{liyanage2005practical} and \cite{see2010robust} that solve inventory management problems. 
More relevant applications are recently studied.
For instance, \cite{ban2019big} and \cite{zhang2023optimalrobust} integrate feature data within 
the newsvendor problem; 
\cite{liu2021time} integrate travel-time predictors with order-assignment optimization to provide last-mile delivery services;
\cite{kallus2022stochastic} propose a new random forest algorithm that considers the downstream optimization problem;
\cite{zhu2022joint} develop a joint estimation and robustness optimization framework;
\cite{qi2022practical} and \cite{ho2022risk} provide an end-to-end framework to integrate prediction and optimization. 
Our paper integrates dimensionality reduction with optimization.

\section{SDP Reformulation} \label{Sec:introdMomDRO}

Given the distribution $\mathbb{P}$ of a random vector $\boldsymbol{\xi} \in \mathbb{R}^m$, the following stochastic programming (SP) formulation seeks an $\mathbf{x} \in \mathcal{X} \subseteq \mathbb{R}^n$ to minimize the expectation of an objective function $f(\mathbf{x},\boldsymbol{\xi})$: 
\begin{equation} \label{Equ:SP}
\min_{\mathbf{x} \in \mathcal{X}} \ \mathbb{E}_{\mathbb{P}} \left[ f\left(\mathbf{x},\boldsymbol{\xi}\right) \right].
\end{equation}
As the distribution $\mathbb{P}$ is often unknown, we assume that $\mathbb{P}$ belongs to a distributional ambiguity set $ \mathcal{D}_{\text{M0}}$ constructed by statistical information estimated from historical data, and then minimize $f(\mathbf{x},\boldsymbol{\xi})$ against the worst-case distribution instead. It leads to the following DRO formulation:
\begin{equation} \label{Equ:DROM1}
\min_{\mathbf{x} \in \mathcal{X}} \ \max_{\mathbb{P} \in \mathcal{D}_{\text{M0}}} \ \mathbb{E}_{\mathbb{P}} \left[ f \left( \mathbf{x}, \boldsymbol{\xi} \right) \right].
\end{equation}
We consider moment-based statistical information \citep{delage2010distributionally} is included in the set $ \mathcal{D}_{\text{M0}}$ as follows:
\[ 
\mathcal{D}_{\text{M0}} \left( \mathcal{S}, \boldsymbol{\mu}, \boldsymbol{\Sigma}, \gamma_1,\gamma_2 \right) = \left\{ \mathbb{P} \ \middle| \ \begin{array}{l}
	\mathbb{P} \left( \boldsymbol{\xi} \in \mathcal{S} \right) = 1 \\ 
	\left( \mathbb{E}_{\mathbb{P}} \left[ \boldsymbol{\xi} \right] - \boldsymbol{\mu} \right)^{\top} \boldsymbol{\Sigma}^{-1} \left( \mathbb{E}_{\mathbb{P}} \left[ \boldsymbol{\xi} \right] - \boldsymbol{\mu} \right) \leq \gamma_1 \\ 
	\mathbb{E}_{\mathbb{P}} \left[ \left( \boldsymbol{\xi} - \boldsymbol{\mu} \right) \left( \boldsymbol{\xi} - \boldsymbol{\mu} \right)^{\top} \right] \preceq \gamma_2 \boldsymbol{\Sigma} \end{array} \right\}, 
\]
which describes that (i) the support of $\boldsymbol{\xi}$ is $\mathcal{S}$, (ii) the mean of $\boldsymbol{\xi}$ lies in an ellipsoid of size $\gamma_1$ centered at $\boldsymbol{\mu}$, and (iii) the covariance of $\boldsymbol{\xi}$ is bounded from above by $\gamma_2 \boldsymbol{\Sigma}$, with $\gamma_1 \geq 0$, $\gamma_2 \geq 1$, and $ \boldsymbol{\Sigma} \succ 0 $.
We perform eigenvalue decomposition on the covariance matrix $\boldsymbol{\Sigma}$ as follows:
\[
\boldsymbol{\Sigma} = \mathbf{U} \boldsymbol{\Lambda} \mathbf{U}^{\top} = \mathbf{U} \boldsymbol{\Lambda}^{\frac{1}{2}} \left( \mathbf{U} \boldsymbol{\Lambda}^{\frac{1}{2}} \right)^{\top},
\]
where $\mathbf{U}\in \mathbb{R}^{m\times m}$ is an orthogonal matrix and $\boldsymbol{\Lambda}\in \mathbb{R}^{m\times m}$ is a diagonal matrix. 
Without loss of generality, we assume that the diagonal elements of $\boldsymbol{\Lambda}$ are arranged in a nonincreasing order.
By letting $\boldsymbol{\xi} = \mathbf{U}\boldsymbol{\Lambda}^{{\frac{1}{2}}} \boldsymbol{\xi}_{\text{I}}+\boldsymbol{\mu}$, we can reformulate Problem \eqref{Equ:DROM1} as:
\begin{equation} \label{Equ:DROM2}
\Theta_{\textup{M}}(m) :=
 \min_{\mathbf{x} \in \mathcal{X}} \ \max_{\mathbb{P}_\text{I} \in \mathcal{D}_{\text{M}}} \ \mathbb{E}_{\mathbb{P}_\text{I}} \left[ f\left(\mathbf{x},\mathbf{U}\boldsymbol{\Lambda}^{{\frac{1}{2}}} \boldsymbol{\xi}_{\text{I}}+\boldsymbol{\mu}\right) \right], 
\end{equation}
where
\[ 
\mathcal{D}_{\text{M}} \left( \mathcal{S}_{\text{I}},\gamma_1,\gamma_2 \right)=\left\{ \mathbb{P}_\text{I} \ \middle| \ \begin{array}{l}
	\mathbb{P}_\text{I} \left( \boldsymbol{\xi}_{\text{I}} \in \mathcal{S}_{\text{I}} \right) = 1, \ \
	\mathbb{E}_{\mathbb{P}_\text{I}} \left[ \boldsymbol{\xi}_{\text{I}}^{\top} \right] \mathbb{E}_{\mathbb{P}_\text{I}} \Big[ \boldsymbol{\xi}_{\text{I}} \Big] \leq \gamma_1 \\
	\mathbb{E}_{\mathbb{P}_\text{I}} \left[ \boldsymbol{\xi}_{\text{I}} \boldsymbol{\xi}_{\text{I}}^{\top} \right] \preceq \gamma_2\mathbf{I}_{m} \end{array} \right\}, \]
with  $\mathcal{S}_{\text{I}}:= \{\boldsymbol{\xi}_{\text{I}} \in \mathbb{R}^m \ | \ \mathbf{U}\boldsymbol{\Lambda}^{{\frac{1}{2}}}\boldsymbol{\xi}_{\text{I}}+\boldsymbol{\mu} \in \mathcal{S} \}$.
Similar to \cite{cheng2018distributionally} and \cite{cheramin2022computationally}, we make the following assumption throughout the paper.

\begin{assumption} \label{Ass:piecewise}
Function $f(\mathbf{x},\boldsymbol{\xi})$ is piecewise linear convex in $\boldsymbol{\xi}$, i.e., $f ( \mathbf{x}, \boldsymbol{\xi} ) = \max_{k=1}^K  \{ y_k^0 (\mathbf{x}) + y_k(\mathbf{x})^{\top} \boldsymbol{\xi} \}$
with $y_k(\mathbf{x}) = (y_k^1(\mathbf{x}), \ldots, y_k^m(\mathbf{x}) )^{\top}$ and $y_k^0(\mathbf{x})$ affine in $\mathbf{x}$ for any $k \in [K]$, 
and $\mathcal{S}$ is polyhedral, i.e., $\mathcal{S} = \{\boldsymbol{\xi} \ | \ \mathbf{A} \boldsymbol{\xi} \leq \mathbf{b} \}$ with $\mathbf{A} \in \mathbb{R}^{l\times m}$ and $\mathbf{b}\in \mathbb{R}^{l}$, with at least one interior point.
\end{assumption}


\begin{proposition}[\citealt{cheramin2022computationally}] \label{Prop:cheramin2022}
Under Assumption \ref{Ass:piecewise}, Problem \eqref{Equ:DROM2} has the same optimal value as the following SDP formulation:
{\small \begin{subequations}
\begin{eqnarray}
\Theta_{\textup{M}}(m) = &\min\limits_{\mathbf{x}, s, {\hat{\boldsymbol{\lambda}}}, \mathbf{q}, \mathbf{Q}} & s + \gamma_2 \mathbf{I}_{m} \bullet \mathbf{Q} + \sqrt{\gamma_1}\left \| \mathbf{q} \right \|_2 \\
& {\normalfont \text{s.t.}} & \begin{bmatrix}  s-y_k^0(\mathbf{x})-\boldsymbol{\lambda}_k^{\top}\mathbf{b}-y_k(\mathbf{x})^{\top}\boldsymbol{\mu}+\boldsymbol{\lambda}_k^{\top}\mathbf{A}\boldsymbol{\mu} & \hspace{0.1 in}  \frac{1}{2}  \left(\mathbf{q} +\left(\mathbf{U}\boldsymbol{\Lambda}^{{\frac{1}{2}}}\right)^{\top}\left(\mathbf{A}^{\top}\boldsymbol{\lambda}_k-y_k(\mathbf{x})\right)\right)^{\top}\\ 
 \frac{1}{2} \left( \mathbf{q} + \left( \mathbf{U}\boldsymbol{\Lambda}^{{\frac{1}{2}}} \right)^{\top} \left( \mathbf{A}^{\top}\boldsymbol{\lambda}_k-y_k(\mathbf{x}) \right) \right) & \mathbf{Q}\end{bmatrix} \succeq 0, \nonumber \\
&& \hspace{4.25in} \forall k\in[K], \label{Cons:SDP} \\
&& \ \boldsymbol{\lambda}_k  \in \mathbb{R}_+^{l}, \ \forall k \in [K], \ \mathbf{x} \in \mathcal{X},\label{Cons:MainProblemcons2}
\end{eqnarray} \label{Equ:MainProblem}
\end{subequations}}%
where ${\hat{\boldsymbol{\lambda}}=\left\{\boldsymbol{\lambda}_1, \ldots, \boldsymbol{\lambda}_K \right\}}$, $\mathbf{q} \in \mathbb{R}^m$, and $\mathbf{Q} \in \mathbb{R}^{m \times m}$.
\end{proposition}

Although Problem \eqref{Equ:MainProblem} is a convex program when $\mathbf{x}$ is given, it can be difficult to solve because a large $m$ leads to high-dimensional SDP constraints at size $m+1$.
As such SDP constraints originate from the covariance matrix $\boldsymbol{\Sigma}$, early attempts in \cite{cheng2018distributionally} and \cite{cheramin2022computationally} exploit the statistical information $\boldsymbol{\Sigma}$ to address the computational challenge while maintaining solution quality.
Specifically, they adopt the PCA, a dimensionality reduction method commonly used in statistical learning, to capture the dominant variability of $\mathbf{U} \boldsymbol{\Lambda}^{{\frac{1}{2}}} \boldsymbol{\xi}_{\text{I}}$ by maintaining the first $m_1 (\leq m)$ components of $\boldsymbol{\xi}_{\text{I}}$ and fixing its other components at $0$; that is,
\begin{align}
   \boldsymbol{\xi} \approx \mathbf{U} \boldsymbol{\Lambda}^{{\frac{1}{2}}}
\left[\boldsymbol{\xi}_{\textup{r}}; \boldsymbol{0}_{m-m_1} \right] + \boldsymbol{\mu} = 
\mathbf{U}_{m \times m_1} \boldsymbol{\Lambda}^{{\frac{1}{2}}}_{m_1} \boldsymbol{\xi}_{\textup{r}} + \boldsymbol{\mu}, \label{eqn:dimen-reduction-1} 
\end{align}
 where $\boldsymbol{\xi}_{ \textup{r} } \in \mathbb{R}^{m_1}$ and $\mathbf{U}_{m \times m_1} \in \mathbb{R}^{m\times m_1}$ and $\boldsymbol{\Lambda}^{{\frac{1}{2}}}_{m_1} \in \mathbb{R}^{m_1\times m_1}$ are upper-left submatrices of $\mathbf{U}$ and $\boldsymbol{\Lambda}$, respectively.
That is, the $m_1$ components of $\boldsymbol{\xi}_{\text{I}}$ corresponding to the largest eigenvalues are maintained as uncertain and the other components are fixed at their means.
With a lower-dimensional random vector $\boldsymbol{\xi}_{\textup{r}}$, 
we can have a relaxation of Problem \eqref{Equ:DROM2}:
\begin{subequations}\label{DRSP-RA} 
\begin{align}
&\Theta_{\textup{M}}(m_1) := \min\limits_{ \mathbf{x} \in \mathcal{X}} \ \max\limits_{\mathbb{P}_{ \textup{r} }  \in \mathcal{D}_{\text{L}}} \ \mathbb{E}_{\mathbb{P}_{ \textup{r} } } \left[ f\left( \mathbf{x}, \mathbf{U}_{m\times m_1} \boldsymbol{\Lambda}^{{\frac{1}{2}}}_{m_1}\boldsymbol{\xi}_{ \textup{r} } + \boldsymbol{\mu} \right) \right],\label{DRSP-RA:obj}\\
& \hspace{-1.1 in} \textup{where} \nonumber\\	
&\mathcal{D}_{\text{L}} \left( \mathcal{S}_{ \textup{r} } ,\gamma_1,\gamma_2 \right) = \left\{ \mathbb{P}_{ \textup{r} } \ \middle| \			
	\begin{array}{l}
	\mathbb{P}_{ \textup{r} } \left( \boldsymbol{\xi}_{ \textup{r} } \in \mathcal{S}_{ \textup{r} } \right) = 1, \ \
	\mathbb{E}_{\mathbb{P}_{ \textup{r} } } \left[\boldsymbol{\xi}_{ \textup{r} }^{\top} \right] \mathbb{E}_{\mathbb{P}_{ \textup{r} } } \Big[\boldsymbol{\xi}_{ \textup{r} } \Big] \leq \gamma_1 \\
	\mathbb{E}_{\mathbb{P}_{ \textup{r} } } \left[ \boldsymbol{\xi}_{ \textup{r} } \boldsymbol{\xi}_{ \textup{r} } ^{\top} \right] \preceq \gamma_2 \mathbf{I}_{m_1}
	\end{array}
	\right\} \label{DRSP-RA:LB-amb} \\
	&\hspace{-1.1 in}\textup{with}\nonumber\\
	&\mathcal{S}_{ \textup{r} }  := \left\{ \boldsymbol{\xi}_{ \textup{r} } \in \mathbb{R}^{m_1} \ \middle| \ \mathbf{U}_{m\times m_1} \boldsymbol{\Lambda}^{{\frac{1}{2}}}_{m_1}\boldsymbol{\xi}_{ \textup{r} } +\boldsymbol{\mu}
	\in \mathcal{S} \right\}.\label{DRSP-RA:support}
\end{align}
\end{subequations}

\noindent
Meanwhile, the corresponding SDP formulation of Problem \eqref{DRSP-RA} has SDP constraints with smaller size at $m_1+1$ and can be solved more efficiently than Problem \eqref{Equ:MainProblem}, leading to an efficient  ``\textit{PCA approximation}.''
Specifically, \cite{cheramin2022computationally} show that the following PCA approximation
{\small \begin{subequations}\label{Equ:PCASDP}
\begin{eqnarray}
 \Theta_{\textup{M}}(m_1) =  & \min\limits_{ \substack{ \mathbf{x},s,{\hat{\boldsymbol{\lambda}}}, \\ \mathbf{q}_{ \textup{r} },\mathbf{Q}_{ \textup{r} } } } & s + \gamma_2 \mathbf{I}_{m_1} \bullet \mathbf{Q}_{ \textup{r} }  + \sqrt{\gamma_1}\left \| \mathbf{q}_{ \textup{r} } \right \|_2 \\
& \hspace{-4cm}  \textnormal{s.t.} &  \hspace{-2cm}   \small{ \begin{bmatrix}  s - y_k^0(\mathbf{x}) - \boldsymbol{\lambda}_k^{\top}\mathbf{b} - y_k(\mathbf{x})^{\top} \boldsymbol{\mu} +\boldsymbol{\lambda}_k^{\top}\mathbf{A}\boldsymbol{\mu} & \hspace{0.1 in} \frac{1}{2}  \left(\mathbf{q}_{\textup{r}} +\left( \mathbf{U}_{m \times m_1} \boldsymbol{\Lambda}^{{\frac{1}{2}}}_{m_1}  \right)^{\top} \left( \mathbf{A}^{\top}\boldsymbol{\lambda}_k-y_k(\mathbf{x}) \right) \right)^{\top} \\ 
	 \frac{1}{2} \left( \mathbf{q}_{\textup{r}} + \left( \mathbf{U}_{m \times m_1} \boldsymbol{\Lambda}^{{\frac{1}{2}}}_{m_1}  \right)^{\top} \left(\mathbf{A}^{\top}\boldsymbol{\lambda}_k-y_k(\mathbf{x})\right) \right) & \mathbf{Q}_{\textup{r}} \end{bmatrix} \succeq 0, } \nonumber \\
&& \hspace{4.0 in} \forall k\in[K],\\
&& \hspace{-2cm} \boldsymbol{\lambda}_k \in \mathbb{R}_+^{l}, \ \forall k\in[K],\mathbf{x} \in \mathcal{X},
\end{eqnarray}%
\end{subequations}}%
where $ \hat{\boldsymbol{\lambda}} =  \{\boldsymbol{\lambda}_1, \ldots, \boldsymbol{\lambda}_K \} $, $ \mathbf{q}_{\textup{r}} \in \mathbb{R}^{m_1} $, and $ \mathbf{Q}_{\textup{r}} \in \mathbb{R}^{m_1 \times m_1}$, provides a \textit{lower bound} for the optimal value of Problem \eqref{Equ:DROM2} (i.e., Problem \eqref{Equ:MainProblem}).
The PCA approximation that leads to an \textit{upper bound} for the optimal value of Problem \eqref{Equ:DROM2} can be similarly derived. Hereafter, we call the problem whose optimal value is a lower bound of the original Problem \eqref{Equ:DROM2} as an \textit{outer approximation}. 
In contrast, the problem generating an upper bound is called an \textit{inner approximation} of Problem \eqref{Equ:DROM2}.




However, relying on only the statistical information (i.e., dominant variability) to choose the components and reducing the high-dimensional uncertainty space may not lead to the best approximation performance.
Although \cite{cheramin2022computationally} provide a performance guarantee to bound the gap between the original and approximated objective values, it is difficult to close the gap when reducing the dimensionality of $\boldsymbol{\xi}_{\text{I}}$. 
Such a difficulty is not surprising because maintaining only the largest statistical variability in the PCA approximations does not capture the optimality conditions of the original problems (e.g., Problem \eqref{Equ:DROM2}).
We provide an example as follows to illustrate that choosing the components of $\boldsymbol{\xi}_{\text{I}}$ corresponding to the largest eigenvalues can be even worse than choosing the components corresponding to the least eigenvalues.

\begin{example} \label{exam:cvar-pca}  
Given $\mathbf{x}\in \mathcal{X}$, we consider the CVaR$_{1-\alpha}$ of a cost function $g(\mathbf{x},\boldsymbol{\xi})$ formulated as the following optimization problem \citep{rockafellar2000optimization}: 
\begin{equation}
\min_{t\in \mathbb{R}} \ t + \frac{1}{\alpha} \mathbb{E}_{\mathbb{P}} \left[ g(\mathbf{x}, \boldsymbol{\xi}) - t \right]^+, \nonumber
\end{equation}
where $\alpha \in (0,1)$ is a risk tolerance level and function $[\cdot]^+:=\max\{0, \cdot\}$. 
For brevity, we let $g(\mathbf{x}, \boldsymbol{\xi}) = \mathbf{x}^{\top} \boldsymbol{\xi}$, $\mathcal{X} = \{\mathbf{x} \in \mathbb{R}^m_+ \ | \ \sum^m_{i=1} x_i = 1 \}$, 
$\mathcal{D} = \{\mathbb{P} \ | \  \mathbb{P} (\boldsymbol{\xi} \in \mathcal{S} ) = 1,  \ \mathbb{E}_{\mathbb{P}} [\boldsymbol{\xi} ] = \boldsymbol{\mu}, \ \mathbb{E}_{\mathbb{P}} [ ( \boldsymbol{\xi} - \boldsymbol{\mu} ) ( \boldsymbol{\xi} - \boldsymbol{\mu} )^{\top} ] \preceq \boldsymbol{\Sigma}  \}  $,
$\mathcal{S}$ is compact, and $\boldsymbol{\mu}$ is in the interior of $\mathcal{S}$.
The distributionally robust counterpart of the above CVaR problem can be formulated as
\begin{align}
    & \min_{\mathbf{x} \in \mathcal{X}} \ \max_{\mathbb{P}\in \mathcal{D}} \ \min_{t \in \mathbb{R} } \ t + \frac{1}{\alpha} \mathbb{E}_{\mathbb{P}}  \left[ g \left( \mathbf{x}, \boldsymbol{\xi} \right) - t \right]^+ \nonumber  \\
    = & \min_{\mathbf{x}\in \mathcal{X}, \ t\in \mathbb{R}} \ \max_{\mathbb{P}\in \mathcal{D}} \ t + \frac{1}{\alpha} \mathbb{E}_{\mathbb{P}}  \left[ g \left( \mathbf{x}, \boldsymbol{\xi} \right) - t \right]^+, \label{Equ:CVaR1}
\end{align}
where the equality holds by the Sion's minimax theorem \citep{sion1958general} because $t + (1/\alpha) \mathbb{E}_{\mathbb{P}}  [ g ( \mathbf{x},\boldsymbol{\xi} ) - t ]^+$ is convex in $t$, concave (specifically, linear) in $\mathbb{P}$, and $\mathcal{D}$ is compact. 
By Proposition \ref{Prop:cheramin2022}, Problem \eqref{Equ:CVaR1} has the same optimal value with the following SDP formulation: 
{\small \begin{subequations} \label{Equ:example-1-sdp}
\begin{eqnarray}
& \min\limits_{ \substack{ \mathbf{x},s,t,\boldsymbol{\lambda}_1, \\ \boldsymbol{\lambda}_2,\mathbf{q},\mathbf{Q} } } & s + \mathbf{I}_{m} \bullet \mathbf{Q}  \\
& \text{s.t.} & 
\begin{bmatrix}  
s - t - \boldsymbol{\lambda}_1^{\top}\mathbf{b}+\boldsymbol{\lambda}_1^{\top}\mathbf{A}\boldsymbol{\mu} &   \frac{1}{2}  \left(\mathbf{q} + \left(\mathbf{U}\boldsymbol{\Lambda}^{{\frac{1}{2}}}\right)^{\top}\mathbf{A}^{\top}\boldsymbol{\lambda}_1\right)^{\top}\\ 
 \frac{1}{2}  \left( \mathbf{q} + \left( \mathbf{U} \boldsymbol{\Lambda}^{{\frac{1}{2}}} \right)^{\top} \mathbf{A}^{\top} \boldsymbol{\lambda}_1\right) & \mathbf{Q} 
 \end{bmatrix} \succeq 0,  \\
 && 
 \begin{bmatrix}  
 s - \left( 1 - \frac{1}{\alpha} \right)  t - \boldsymbol{\lambda}_2^{\top} \mathbf{b} - \left( \frac{1}{\alpha} \mathbf{x} \right)^{\top} \boldsymbol{\mu} + \boldsymbol{\lambda}_2^{\top} \mathbf{A} \boldsymbol{\mu} &  \frac{1}{2}  \left(\mathbf{q} + \left( \mathbf{U} \boldsymbol{\Lambda}^{{\frac{1}{2}}} \right)^{\top} \left( \mathbf{A}^{\top} \boldsymbol{\lambda}_2 - \frac{1}{\alpha} \mathbf{x} \right) \right)^{\top} \\ 
 \frac{1}{2}  \left( \mathbf{q} + \left( \mathbf{U} \boldsymbol{\Lambda}^{{\frac{1}{2}}} \right)^{\top} \left( \mathbf{A}^{\top} \boldsymbol{\lambda}_2 - \frac{1}{\alpha} \mathbf{x} \right) \right) & \mathbf{Q}
 \end{bmatrix} \succeq 0, \\
&& \mathbf{x} \in \mathcal{X}, \ t \in \mathbb{R}, \ \boldsymbol{\lambda}_1  \in \mathbb{R}_+^{l}, \ \boldsymbol{\lambda}_2  \in \mathbb{R}_+^{l}. \nonumber
\end{eqnarray} 
\end{subequations}}

\noindent 
Let $\alpha=0.05$, 
$\mathcal{S} = \{ \boldsymbol{\xi} \in \mathbb{R}^3 \ | \ 0 \leq \xi_1 \leq 8 , 1 \leq \xi_2 \leq 12, 2 \leq \xi_3 \leq 16 \}$, 
$\boldsymbol{\mu}=[1,2,3]$, 
$\boldsymbol{\Sigma} = 
{\scriptsize \begin{bmatrix}
1 & 0.2 & 0.1\\
0.2 & 3 & 0.15\\
0.1 & 0.15 & 2
\end{bmatrix}}$ 
with eigenvalues $3.0441$, $1.9827$, and $0.9731$.
Solving Problem \eqref{Equ:example-1-sdp} gives the optimal value $5.0214$ with $x_1 = 0.7194$, $x_2 = 0.1354$, $x_3 = 0.1452$, and $t=3.1289$.
When we follow \cite{cheng2018distributionally} and \cite{cheramin2022computationally} to perform PCA approximation over Problem \eqref{Equ:CVaR1} by capturing only one component of the three components in $\boldsymbol{\xi}$, we observe the following:
\begin{itemize}
    \item Choosing the component corresponding to the largest eigenvalue $3.0441$, 
    the PCA approximation gives the optimal value at $1.7877$ 
    with $x_1 = 1$, $x_2 = 0$, $x_3 = 0$, and $t=1.3731$.
    \item Choosing the component corresponding to the second largest eigenvalue 
    $1.9827$, the PCA approximation gives the optimal value at $1.2999$ 
    with $x_1 = 0.7001$, $x_2 = 0.2999$, $x_3 = 0$, and $t=1.2999$.
    \item Choosing the component corresponding to the smallest eigenvalue 
    $0.9731$, the PCA approximation gives the optimal value at 
    $1.9154$ with $x_1 = 0.0846$, $x_2 = 0.9154$, $x_3 = 0$, and $t=1.9154$.
\end{itemize}
\end{example}


Example \ref{exam:cvar-pca} shows that performing dimensionality reduction (i.e., from $ \boldsymbol{\xi} $ to $ \boldsymbol{\xi}_{ \textup{r} } $) using the components with the largest variability may not produce a good optimal value from the \textit{subsequent} PCA approximation (i.e., an SDP) 
and it can be even worse than using the component with the smallest variability. 
To solve this issue, we integrate the dimensionality reduction with the subsequent approximation in the following sections, leading to an \textit{optimized dimensionality reduction (ODR) approach}. 
Correspondingly, we obtain efficient lower and upper bounds in the following Sections \ref{Sec:ODRapproach}--\ref{sec:lower-bound-2} and more importantly, the bounds can achieve the optimal value of the original Problem \eqref{Equ:DROM2}.

\section{Lower Bound} \label{Sec:ODRapproach} 

We extend the dimensionality reduction method (i.e., PCA) in \eqref{eqn:dimen-reduction-1} by introducing a decision variable $\mathbf{B} \in  \mathcal{B}_{m_1} := \{\mathbf{B}\in \mathbb{R}^{m\times m_1} \ | \ \mathbf{B}^{\top} \mathbf{B} = \mathbf{I}_{m_1} \}  \subseteq \mathbb{R}^{m \times m_1}$ such that 
\begin{align}
    \boldsymbol{\xi} = \mathbf{U}\boldsymbol{\Lambda}^{{\frac{1}{2}}} \boldsymbol{\xi}_{\text{I}} + \boldsymbol{\mu} \approx \mathbf{U} \boldsymbol{\Lambda}^{{\frac{1}{2}}} \mathbf{B}  \boldsymbol{\xi}_{ \textup{r} } +\boldsymbol{\mu}, \label{eqn:dimen-reduction-B}
\end{align}
where $\mathbf{B}$ will be optimized in the subsequent PCA approximation, i.e., optimized dimensionality reduction. 
By \eqref{eqn:dimen-reduction-B}, we project $\boldsymbol{\xi}_{ \textup{I} }$ onto a subspace of $\mathbb{R}^{m\times m}$ 
spanned by $\mathbf{B} \in  \mathcal{B}_{m_1} $ and approximate $\boldsymbol{\xi}_{ \textup{I} }$ by the projection $\mathbf{B} \boldsymbol{\xi}_{ \textup{r} }$. 
Note that when $\mathbf{B} = 
{\scriptsize \begin{bmatrix}
\mathbf{I}_{m_1} \\
\mathbf{0}_{(m-m_1) \times m_1}
\end{bmatrix} }$, 
\eqref{eqn:dimen-reduction-B} reduces to \eqref{eqn:dimen-reduction-1}.
By adopting \eqref{eqn:dimen-reduction-B}, when we still reduce the dimensionality space from $m$ to $m_1$, we allow $ \mathbf{B}  \boldsymbol{\xi}_{ \textup{r} } $ to take linear combinations of the original components of $ \boldsymbol{\xi}_{\text{I}}$, instead of taking only the components corresponding to the largest eigenvalues. 
Therefore, we would like to choose a good (even an optimal) $ \mathbf{B} $ to obtain a better lower bound for Problem \eqref{Equ:DROM2} than Problem \eqref{Equ:PCASDP}.


Given any $m_1 \in [m]$ and $\mathbf{B} \in \mathcal{B}_{m_1}$, we obtain a relaxation of Problem \eqref{Equ:DROM2} by extending Problem \eqref{DRSP-RA}. 
If the relaxation provides a lower bound for the optimal value of Problem \eqref{Equ:DROM2}, then we may choose the best $ \mathbf{B} \in \mathcal{B}_{m_1}$ such that we obtain the largest possible lower bound.
Thus, we build the following \textit{integrated dimensionality reduction and optimization} problem:
\begin{equation} \label{Equ:DROM3}  
\Theta_{\textup{L}}(m_1) =  \max\limits_{\mathbf{B} \in \mathcal{B}_{m_1}} \ \min\limits_{\mathbf{x} \in \mathcal{X}} \ \max\limits_{\mathbb{P}_{ \textup{r} }  \in \mathcal{D}_{\text{L}} } \ \mathbb{E}_{\mathbb{P}_{\textup{r}} } \left[ f\left(\mathbf{x},\mathbf{U}\boldsymbol{\Lambda}^{{\frac{1}{2}}} \mathbf{B}  \boldsymbol{\xi}_{ \textup{r} } +\boldsymbol{\mu}\right) \right],
\end{equation}
where $\mathcal{D}_{\text{L}}$ is defined in \eqref{DRSP-RA:LB-amb} with
\begin{align}
    \mathcal{S}_{ \textup{r} }  := \left\{ \boldsymbol{\xi}_{ \textup{r} } \in \mathbb{R}^{m_1} \ \middle| \ \mathbf{U}\boldsymbol{\Lambda}^{{\frac{1}{2}}}\mathbf{B} \boldsymbol{\xi}_{ \textup{r} } +\boldsymbol{\mu} \in \mathcal{S} \right\}. \label{eqn:support-sr-b}
\end{align}

\noindent 
We will show that Problem \eqref{Equ:DROM3} provides a lower bound for Problem \eqref{Equ:DROM2} (see Theorem \ref{Theo:lowerbound}). 
Before presenting this theorem, we prepare the following two lemmas.

\begin{lemma} \label{Lem:3cons}
When $\mathbf{B} \in \mathbb{R}^{m\times m_1}$, the following three constraints are equivalent:
(i) ${\scriptsize \begin{bmatrix}
 \mathbf{I}_m & \mathbf{B} \\ \mathbf{B}^{\top} & \mathbf{I}_{m_1}
\end{bmatrix} } \succeq 0$, (ii) $\mathbf{B} \mathbf{B}^{\top} \preceq \mathbf{I}_m$, and (iii) $\mathbf{B}^{\top} \mathbf{B} \preceq \mathbf{I}_{m_1}$.
\end{lemma}

Lemma \ref{Lem:3cons} shows that both $\mathbf{B} \mathbf{B}^{\top} \preceq \mathbf{I}_m$ and $\mathbf{B}^{\top} \mathbf{B} \preceq \mathbf{I}_{m_1}$ can be reformulated as an SDP constraint ${\scriptsize \begin{bmatrix}
 \mathbf{I}_m & \mathbf{B} \\ \mathbf{B}^{\top} & \mathbf{I}_{m_1}
\end{bmatrix}} \succeq 0$.
Although this SDP constraint has a high dimension at $m+m_1$, it is very sparse and usually does not create additional computational challenges.


\begin{lemma} \label{Lem:VXV}
For any matrix $\mathbf{V} \in \mathbb{R}^{m\times n}$ and symmetric matrices $\mathbf{X} \in \mathbb{R}^{m\times m}$ and $\mathbf{Y} \in \mathbb{R}^{m \times m}$, we have: (i) If $\mathbf{X} \succeq \mathbf{Y}$, then $\mathbf{V}^{\top} \mathbf{X} \mathbf{V} \succeq \mathbf{V}^{\top} \mathbf{Y} \mathbf{V}$; (ii) If $n=m$ and $\mathbf{V}$ is invertible, then $\mathbf{X} \succeq \mathbf{Y}$ is equivalent to $\mathbf{V}^{\top} \mathbf{X} \mathbf{V} \succeq \mathbf{V}^{\top} \mathbf{Y} \mathbf{V}$.
\end{lemma}

Lemma \ref{Lem:VXV} shows that a PSD matrix (e.g., $\mathbf{X} - \mathbf{Y}$) remains PSD if it is pre-multiplied by an arbitrary matrix with appropriate dimensions (e.g., $ \mathbf{V}^{\top}$) and post-multiplied by this arbitrary matrix's transpose (e.g., $ \mathbf{V}$).
Furthermore, if this arbitrary matrix is invertible, then the original PSD matrix is equivalent to the matrix after the pre-multiplication and post-multiplication.
With Lemmas \ref{Lem:3cons} and \ref{Lem:VXV}, we are now ready to present the following theorem.

\begin{theorem} \label{Theo:lowerbound}
The following three conclusions hold:
(i) Problem \eqref{Equ:DROM3} provides a lower bound for the optimal value of Problem \eqref{Equ:DROM2}, i.e., $\Theta_{\textup{L}}(m_1) \leq \Theta_{\textup{M}}(m)$ for any $m_1 \leq m$; 
(ii) the optimal value of Problem \eqref{Equ:DROM3} is nondecreasing in $m_1$, i.e., $\Theta_{\textup{L}}(m_1) \leq \Theta_{\textup{L}}(m_2)$ for any $m_1 < m_2 \leq m$; 
and (iii) when $m_1 = m$, Problem \eqref{Equ:DROM2} and Problem \eqref{Equ:DROM3} have the same optimal value, i.e., $\Theta_{\textup{L}}(m) = \Theta_{\textup{M}}(m)$. 
\end{theorem}

Theorem \ref{Theo:lowerbound} shows that we obtain a lower bound for the optimal value of Problem \eqref{Equ:DROM2} when reducing the dimensionality space of $ \boldsymbol{\xi}_{\text{I}}$ \textit{while optimizing the choice of $\mathbf{B} \in \mathcal{B}_{m_1}$} in Problem \eqref{Equ:DROM3}.
When the reduced dimensionality (i.e., $m_1$) is higher, we obtain a better lower bound. 
We maintain the optimal value of Problem \eqref{Equ:DROM2} if the dimensionality space is not reduced (i.e., $m_1 = m$).
Note that the conclusions in Theorem \ref{Theo:lowerbound} are similar to Theorem 2 in \cite{cheramin2022computationally}, both demonstrating the validity of dimensionality reduction in solving the moment-based DRO problems.
However, here by optimizing the choice of $\mathbf{B} \in \mathcal{B}_{m_1}$, Problem \eqref{Equ:DROM3} provides a better lower bound than Problem \eqref{DRSP-RA} (i.e., the PCA approximation in \citealt{cheramin2022computationally}) does because the latter problem is a special case of the former problem.
More importantly, we may expect to close the gap between $\Theta_{\textup{L}}(m_1) $ and $ \Theta_{\textup{M}}(m)$ when we choose a small $m_1$.
To that end, we follow the PCA approximation \eqref{Equ:PCASDP} to reformulate Problem \eqref{Equ:DROM3} as the following SDP formulation:
\begin{equation}
\Theta_{\textup{L}} (m_1) = \max\limits_{\mathbf{B} \in \mathcal{B}_{m_1}} \ \underline{\Theta}(m_1, \mathbf{B}), \label{Equ:appro}
\end{equation}
where
{\small \begin{subequations} 
\begin{eqnarray}
\underline{\Theta}(m_1,\mathbf{B}) :=  & \min\limits_{ \substack{\mathbf{x}, s, {\hat{\boldsymbol{\lambda}}}, \\ \mathbf{q}_{ \textup{r} }, \mathbf{Q}_{ \textup{r} } } } & s + \gamma_2 \mathbf{I}_{m_1} \bullet \mathbf{Q}_{ \textup{r} }  + \sqrt{\gamma_1}\left \| \mathbf{q}_{ \textup{r} } \right \|_2 \label{Approximation1:obj}\\
&\textnormal{s.t.} &   \small{ \begin{bmatrix}  s - y_k^0(\mathbf{x}) - \boldsymbol{\lambda}_k^{\top} b - y_k(\mathbf{x})^{\top} \boldsymbol{\mu} + \boldsymbol{\lambda}_k^{\top} \mathbf{A} \boldsymbol{\mu} & \hspace{0.1 in} \frac{1}{2}  \left(\mathbf{q}_{ \textup{r} } +\left( \mathbf{U} \boldsymbol{\Lambda}^{{\frac{1}{2}}} \mathbf{B} \right)^{\top} \left( \mathbf{A}^{\top} \boldsymbol{\lambda}_k - y_k(\mathbf{x}) \right) \right)^{\top} \\
	 \frac{1}{2} \left( \mathbf{q}_{ \textup{r} } +\left(\mathbf{U} \boldsymbol{\Lambda}^{{\frac{1}{2}}}\mathbf{B} \right)^{\top}\left(\mathbf{A}^{\top}\boldsymbol{\lambda}_k-y_k(\mathbf{x})\right) \right) & \mathbf{Q}_{ \textup{r} } \end{bmatrix} \succeq 0,} \nonumber \\
&& \hspace{4.0 in} \forall k\in[K], \label{Cons:appro1} \\
&& \mathbf{x} \in \mathcal{X}; \ \hat{\boldsymbol{\lambda}} = \left\{\boldsymbol{\lambda}_1, \dots, \boldsymbol{\lambda}_K \right\}, \ \boldsymbol{\lambda}_k \in \mathbb{R}_+^{l}, \ \forall k\in[K]; \ \mathbf{q}_{ \textup{r} } \in \mathbb{R}^{m_1}; \ \mathbf{Q}_{ \textup{r} } \in \mathbb{R}^{ m_1 \times m_1 }. \label{Cons:appro2}
\end{eqnarray}
\label{Equ:appro-1}
\end{subequations}}%

Now we would like to find an $m_1 < m$ such that $\Theta_{\textup{L}}(m_1) $ in Problem \eqref{Equ:appro} is close (even equal) to $ \Theta_{\textup{M}}(m)$ in Problem \eqref{Equ:MainProblem}. 
Note that, if $\Theta_{\textup{L}}(m_1) = \Theta_{\textup{M}}(m)$, then comparing the SDP constraints between \eqref{Equ:MainProblem} and \eqref{Equ:appro} shows that the rank of $ \mathbf{Q} $ in the optimal solution of Problem \eqref{Equ:MainProblem} can be smaller than $m$.
Specifically, 
we have
the following significant conclusion holds.

\begin{theorem}  \label{Prop:lowrank}
Consider $K < m$ and any optimal solution $(\mathbf{x}^*, s^*, \hat{\boldsymbol{\lambda}}^*, \mathbf{q}^*, \mathbf{Q}^*)$ of Problem \eqref{Equ:MainProblem} with $S_k = s^* - y_k^0 (\mathbf{x}^*) - \boldsymbol{\lambda}_k^{* \top} \mathbf{b} - y_k(\mathbf{x}^*)^{\top} \boldsymbol{\mu} + \boldsymbol{\lambda}_k^{* \top} \mathbf{A} \boldsymbol{\mu}$ for any $k \in [K]$.
We can always construct another optimal solution $(\mathbf{x}^*, s^*, \hat{\boldsymbol{\lambda}}^*, \mathbf{q}^{\prime}, \mathbf{Q}^{\prime})$ of Problem \eqref{Equ:MainProblem} such that
${\rm{rank}}( \mathbf{Q}^{\prime} )  \leq K$, 
$ \mathbf{q}^{\prime} = \mathbf{V} \boldsymbol{\delta} $, 
$\mathbf{Q}^{\prime} = \mathbf{V} \mathbf{Y}_{11} \mathbf{V}^{\top}$, 
and $(\mathbf{U} \boldsymbol{\Lambda}^{{\frac{1}{2}}} )^{\top}  ( \mathbf{A}^{\top} \boldsymbol{\lambda}^*_k - y_k(\mathbf{x}^*) ) = \mathbf{V} \boldsymbol{\nu}_k$ 
for any $k \in [K]$, 
where $\mathbf{Y}_{11} \in \mathbb{R}^{K \times K} $, 
$ \mathbf{Y}_{11} \succeq 0$, 
$ \mathbf{V} = [ \mathbf{v}_k, \ \forall k \in [K] ] \in \mathbb{R}^{m \times K}$ with orthonormal vectors $\mathbf{v}_k\in \mathbb{R}^m$, 
$\boldsymbol{\delta} \in \mathbb{R}^{K}$, 
and $\boldsymbol{\nu}_k \in \mathbb{R}^{K}$ for any $k \in [K]$.
\end{theorem}
\proof{Proof.}
Note that the optimal solution  $(\mathbf{x}^*, s^*, \hat{\boldsymbol{\lambda}}^*, \mathbf{q}^*, \mathbf{Q}^*)$ of Problem \eqref{Equ:MainProblem} leads to the optimal value $ s^* +\gamma_2\mathbf{I}_{m} \bullet \mathbf{Q}^* + \sqrt{\gamma_1}\left \| \mathbf{q}^*\right \|_2$. 
Based on this optimal solution, we construct a feasible solution of Problem \eqref{Equ:MainProblem}, denoted by 
$(\mathbf{x}^{\prime}, s^{\prime}, \hat{\boldsymbol{\lambda}}^{\prime}, \mathbf{q}^{\prime}, \mathbf{Q}^{\prime})$ such that $\mathbf{x}^{\prime} = \mathbf{x}^*$, $s^{\prime} = s^*$, and $\hat{\boldsymbol{\lambda}}^{\prime} = \hat{\boldsymbol{\lambda}}^*$.

Now we construct the values of $\mathbf{q}^{\prime}$ and $\mathbf{Q}^{\prime}$. 
By constraints \eqref{Cons:SDP}, we have
\begin{equation} 
\small{\begin{bmatrix}  S_k & \frac{1}{2} \left(\mathbf{q}^* + \left(\mathbf{U} \boldsymbol{\Lambda}^{{\frac{1}{2}}} \right)^{\top}  \left(\mathbf{A}^{\top} \boldsymbol{\lambda}^*_k - y_k \left( \mathbf{x}^* \right) \right) \right)^{\top}\\ 
 \frac{1}{2} \left(\mathbf{q}^* + \left(\mathbf{U} \boldsymbol{\Lambda}^{{\frac{1}{2}}} \right)^{\top}  \left(\mathbf{A}^{\top} \boldsymbol{\lambda}^*_k - y_k \left( \mathbf{x}^* \right) \right) \right) & \mathbf{Q}^* \end{bmatrix}} \succeq 0, \ \ \forall k\in[K]. \label{eqn:prop-2-sdp-1}
\end{equation}
By Schur complement, we can equivalently rewrite \eqref{eqn:prop-2-sdp-1} as 
\begin{align}
4S_k \mathbf{Q}^* \succeq \left(\mathbf{q}^* + \left(\mathbf{U} \boldsymbol{\Lambda}^{{\frac{1}{2}}} \right)^{\top}  \left(\mathbf{A}^{\top} \boldsymbol{\lambda}^*_k-y_k(\mathbf{x}^*) \right) \right) \left(\mathbf{q}^* + \left(\mathbf{U} \boldsymbol{\Lambda}^{{\frac{1}{2}}} \right)^{\top}  \left(\mathbf{A}^{\top} \boldsymbol{\lambda}^*_k-y_k(\mathbf{x}^*) \right) \right)^{\top}, \forall k \in [K]. \label{Cons:Prop2-1}
\end{align}
Note that $K < m$. Thus, through the Gram–Schmidt process, we can always construct $K$ orthonormal vectors $\mathbf{v}_k \in \mathbb{R}^m, \ \forall k \in [K],$ and $K$ real vectors $\mathbf{a}_k \in \mathbb{R}^{K}, \ \forall k \in [K]$, such that
\begin{align}
     \left(\mathbf{U} \boldsymbol{\Lambda}^{{\frac{1}{2}}} \right)^{\top}  \left( \mathbf{A}^{\top} \boldsymbol{\lambda}^*_k - y_k(\mathbf{x}^*) \right) = \mathbf{V} \mathbf{a}_k, \ \forall k \in [K], \label{Cons:Prop2-2}
\end{align}
$\mathbf{V} = [\mathbf{v}_k,\ \forall k \in [K] ] \in \mathbb{R}^{m\times K}$. 
We further extend $\mathbf{V}$ to $[\mathbf{V},\bar{\mathbf{V}} ] \in \mathbb{R}^{m\times m}$ with $\bar{\mathbf{V}} \in \mathbb{R}^{m \times (m-K)}$ such that all the column vectors of $[\mathbf{V}, \bar{\mathbf{V}} ]$ can span the space of $\mathbb{R}^m$.
As $\mathbf{q}^* \in \mathbb{R}^m$, we can find $\mathbf{a}_0 \in \mathbf{R}^{K}$ and $\bar{\mathbf{a}}_0 \in \mathbf{R}^{m-K}$ such that 
\begin{align}
    \mathbf{q}^* = \mathbf{V} \mathbf{a}_0 + \bar{\mathbf{V}} \bar{\mathbf{a}}_0. \label{Cons:Prop2-3}
\end{align}
As $\mathbf{Q}^* \in \mathbb{R}^{m \times m}$, we can then decompose $\mathbf{Q}^*$ as
\begin{align}
    \mathbf{Q}^* & = \begin{bmatrix}  \mathbf{V} & \bar{\mathbf{V}}\end{bmatrix} \begin{bmatrix}  \mathbf{Y}_{11} & \mathbf{Y}_{12}\\ \mathbf{Y}_{21} & \mathbf{Y}_{22} \end{bmatrix} \begin{bmatrix}  \mathbf{V}^{\top} \\ \bar{\mathbf{V}}^{\top}\end{bmatrix} = \mathbf{V}\mathbf{Y}_{11} \mathbf{V}^{\top} + \bar{\mathbf{V}} \mathbf{Y}_{21} \mathbf{V}^{\top} + \mathbf{V} \mathbf{Y}_{12} \bar{\mathbf{V}}^{\top} + \bar{\mathbf{V}} \mathbf{Y}_{22} \bar{\mathbf{V}}^{\top},  \label{Cons:Prop2-4}
\end{align}
where $\mathbf{Y}_{11} \in  \mathbb{R}^{K \times K}$, $\mathbf{Y}_{12} \in  \mathbb{R}^{K \times (m - K)}$, 
$\mathbf{Y}_{21} \in  \mathbb{R}^{ (m - K) \times K}$,
 and $\mathbf{Y}_{22} \in  \mathbb{R}^{ (m - K) \times (m - K)}$. 
As $\mathbf{Q}^* \succeq 0$, we have  
$ {\scriptsize  \begin{bmatrix}  \mathbf{Y}_{11} & \mathbf{Y}_{12}\\ \mathbf{Y}_{21} & \mathbf{Y}_{22} \end{bmatrix} } = 
 {\scriptsize \begin{bmatrix}  \mathbf{V} & \bar{\mathbf{V}} \end{bmatrix}^{-1} \mathbf{Q}^* \begin{bmatrix}  \mathbf{V}^{\top} \\ \bar{\mathbf{V}}^{\top}\end{bmatrix}^{-1} } \succeq 0$ by Lemma \ref{Lem:VXV}. 
By \eqref{Cons:Prop2-1}, \eqref{Cons:Prop2-2}, and \eqref{Cons:Prop2-3}, we have 
\begin{align} 
4S_k\mathbf{Q}^* \succeq \left( \mathbf{V} \mathbf{a}_0 + \bar{\mathbf{V}} \bar{\mathbf{a}}_0+\mathbf{V} \mathbf{a}_k \right) \left( \mathbf{V} \mathbf{a}_0 + \bar{\mathbf{V}} \bar{\mathbf{a}}_0 + \mathbf{V} \mathbf{a}_k \right)^{\top}, \ \forall k \in [K]. \label{Cons:Prop2-5}
\end{align}
By \eqref{Cons:Prop2-4} and \eqref{Cons:Prop2-5}, we have
\begin{align*}
4 S_k \left(\mathbf{V} \mathbf{Y}_{11} \mathbf{V}^{\top} + \bar{\mathbf{V}} \mathbf{Y}_{21} \mathbf{V}^{\top} + \mathbf{V} \mathbf{Y}_{12} \bar{\mathbf{V}}^{\top} + \bar{\mathbf{V}} \mathbf{Y}_{22} \bar{\mathbf{V}}^{\top}\right) \succeq \left(\mathbf{V} \mathbf{a}_0 + \bar{\mathbf{V}} \bar{\mathbf{a}}_0+\mathbf{V} \mathbf{a}_k \right)\left(\mathbf{V} \mathbf{a}_0 + \bar{\mathbf{V}} \bar{\mathbf{a}}_0 + \mathbf{V} \mathbf{a}_k \right)^{\top}, \forall k \in [K].
\end{align*} 
By Lemma \ref{Lem:VXV}, we further have
\begin{align}
&4S_k \mathbf{V}^{\top}(\mathbf{V}\mathbf{Y}_{11} \mathbf{V}^{\top} + \bar{\mathbf{V}} \mathbf{Y}_{21} \mathbf{V}^{\top} + \mathbf{V} \mathbf{Y}_{12} \bar{\mathbf{V}}^{\top} + \bar{\mathbf{V}} \mathbf{Y}_{22} \bar{\mathbf{V}}^{\top}) \mathbf{V}  \nonumber \\
& \hspace{2cm} \succeq \mathbf{V}^{\top} (\mathbf{V} \mathbf{a}_0 + \bar{\mathbf{V}} \bar{\mathbf{a}}_0+\mathbf{V} \mathbf{a}_k )(\mathbf{V} \mathbf{a}_0 + \bar{\mathbf{V}} \bar{\mathbf{a}}_0 + \mathbf{V} \mathbf{a}_k )^{\top} \mathbf{V}, \ \forall k \in [K]. \label{Cons:Prop2-6} 
\end{align}
Because $ \mathbf{V}^{\top} \bar{\mathbf{V}} = \mathbf{0} $, $ \bar{\mathbf{V}}^{\top} \mathbf{V} = \mathbf{0}$, and $\mathbf{V}^{\top} \mathbf{V}  = \mathbf{I}_{K} $, constraints \eqref{Cons:Prop2-6} become
\begin{align}
 & 4 S_k \mathbf{Y}_{11} \succeq (\mathbf{a}_0 + \mathbf{a}_k) (\mathbf{a}_0 +\mathbf{a}_k)^{\top}, \ \forall k \in [K]. \label{Cons:Prop2-7}
\end{align}
Now we let $\mathbf{q}^{\prime} = \mathbf{V} \mathbf{a}_0$ and $\mathbf{Q}^{\prime} = \mathbf{V} \mathbf{Y}_{11} \mathbf{V}^{\top}$. 
By \eqref{Cons:Prop2-7} and Lemma \ref{Lem:VXV}, we have
\begin{align}
    4S_k\mathbf{Q}^{\prime}  = & 4S_k\mathbf{V} \mathbf{Y}_{11} \mathbf{V}^{\top}  \succeq (\mathbf{V}\mathbf{a}_0 + \mathbf{V}\mathbf{a}_k) (\mathbf{V}\mathbf{a}_0 + \mathbf{V}\mathbf{a}_k)^{\top} \nonumber \\
    = & \left( \mathbf{q}^{\prime} + \left( \mathbf{U} \boldsymbol{\Lambda}^{{\frac{1}{2}}} \right)^{\top} \left( \mathbf{A}^{\top} \boldsymbol{\lambda}^*_k - y_k (\mathbf{x}^*) \right) \right) \left( \mathbf{q}^{\prime} + \left( \mathbf{U} \boldsymbol{\Lambda}^{{\frac{1}{2}}} \right)^{\top}  \left( \mathbf{A}^{\top} \boldsymbol{\lambda}^*_k - y_k(\mathbf{x}^*) \right) \right)^{\top}, \ \forall k \in [K]. \label{Cons:Prop2-8}
\end{align}
Comparing \eqref{Cons:SDP} and \eqref{Cons:Prop2-8}, 
we have $(\mathbf{x}^{\prime}, s^{\prime}, \hat{\boldsymbol{\lambda}}^{\prime}, \mathbf{q}^{\prime}, \mathbf{Q}^{\prime})$ is a feasible solution of Problem \eqref{Equ:MainProblem}
and the corresponding objective value is
\begin{align}
    s^{\prime} + \gamma_2 \mathbf{I}_{m} \bullet \mathbf{Q}^{\prime} + \sqrt{\gamma_1} \left \| \mathbf{q}^{\prime} \right \|_2 
    \geq s^* +\gamma_2 \mathbf{I}_{m} \bullet \mathbf{Q}^* + \sqrt{\gamma_1} \left \| \mathbf{q}^* \right \|_2, \label{Cons:Prop2-9}
\end{align}
where the inequality holds because $(\mathbf{x}^{\prime}, s^{\prime}, \hat{\boldsymbol{\lambda}}^{\prime}, \mathbf{q}^{\prime}, \mathbf{Q}^{\prime})$ is a feasible solution of Problem \eqref{Equ:MainProblem} and Problem \eqref{Equ:MainProblem} is a minimization problem.
Note that 
\begin{align*}
   &  \mathbf{I}_{m} \bullet \mathbf{Q}^* = tr(\mathbf{Q}^*) = tr \left(\begin{bmatrix}  \mathbf{V} & \bar{\mathbf{V}}\end{bmatrix} \begin{bmatrix}  \mathbf{Y}_{11} & \mathbf{Y}_{12}\\ \mathbf{Y}_{21} & \mathbf{Y}_{22}\end{bmatrix} \begin{bmatrix}  \mathbf{V}^{\top} \\ \bar{\mathbf{V}}^{\top}\end{bmatrix} \right) = tr \left(  \begin{bmatrix}  \mathbf{Y}_{11} & \mathbf{Y}_{12}\\ \mathbf{Y}_{21} & \mathbf{Y}_{22}\end{bmatrix} \begin{bmatrix}  \mathbf{V}^{\top} \\ \bar{\mathbf{V}}^{\top}\end{bmatrix} \begin{bmatrix}  \mathbf{V} & \bar{\mathbf{V}}\end{bmatrix} \right) = tr \left(  \begin{bmatrix}  \mathbf{Y}_{11} & \mathbf{Y}_{12}\\ \mathbf{Y}_{21} & \mathbf{Y}_{22}\end{bmatrix} \right) \\
   & = \mathbf{I}_{K} \bullet \mathbf{Y}_{11} + \mathbf{I}_{m-K} \bullet \mathbf{Y}_{22} \geq \mathbf{I}_{K} \bullet \mathbf{Y}_{11} = tr (\mathbf{Y}_{11} ) = tr ( \mathbf{Y}_{11} \mathbf{V}^{\top} \mathbf{V} ) = tr ( \mathbf{V} \mathbf{Y}_{11} \mathbf{V}^{\top}  ) = tr (\mathbf{Q}^{\prime}) = \mathbf{I}_{m} \bullet \mathbf{Q}^{\prime},
\end{align*}
where the first equality holds by the definition of a matrix's trace, 
the second equality holds by \eqref{Cons:Prop2-4},
the third equality holds by the cyclic property of a matrix's trace,
the fourth equality holds because 
$  {\scriptsize  \begin{bmatrix}  \mathbf{V}^{\top} \\ \bar{\mathbf{V}}^{\top}\end{bmatrix} \begin{bmatrix}  \mathbf{V} & \bar{\mathbf{V}}\end{bmatrix}} = \mathbf{I}_{m} $, and
the first inequality holds because 
$  {\scriptsize  \begin{bmatrix}  \mathbf{Y}_{11} & \mathbf{Y}_{12}\\ \mathbf{Y}_{21} & \mathbf{Y}_{22}\end{bmatrix} } \succeq 0$
and accordingly $\mathbf{I}_{m-K} \bullet \mathbf{Y}_{22} \geq 0$.
Meanwhile, 
\begin{align*}
  \left \| \mathbf{q}^* \right \|_2^2 = & \left( \mathbf{q}^{*} \right)^{\top} \mathbf{q}^* =  \left( \mathbf{V} \mathbf{a}_0 + \bar{\mathbf{V}} \bar{\mathbf{a}}_0 \right)^{\top} \left( \mathbf{V} \mathbf{a}_0 + \bar{\mathbf{V}} \bar{\mathbf{a}}_0 \right) =  \left( \mathbf{a}_0^{\top} \mathbf{a}_0 + \bar{\mathbf{a}}_0^{\top} \bar{\mathbf{a}}_0 \right)  \\
     \geq &  \mathbf{a}_0^{\top} \mathbf{a}_0 = \left( \mathbf{V} \mathbf{a}_0 \right)^{\top} \left( \mathbf{V} \mathbf{a}_0 \right) = \left( \mathbf{q}^{\prime} \right)^{\top} \mathbf{q}^{\prime} = \left \| \mathbf{q}^{\prime} \right \|_2^2,
\end{align*}
where the second equality holds by \eqref{Cons:Prop2-3}, the third equality holds because $ \mathbf{V}^{\top} \bar{\mathbf{V}} = \mathbf{0}, \ \bar{\mathbf{V}}^{\top} \mathbf{V} = \mathbf{0},\  \mathbf{V}^{\top} \mathbf{V}  = \mathbf{I}_{K}, \ \bar{\mathbf{V}}^{\top} \bar{\mathbf{V}}  = \mathbf{I}_{m-K}$, and the first inequality holds because $\bar{\mathbf{a}}_0^{\top} \bar{\mathbf{a}}_0 \geq 0$. Thus, we have
\begin{align}
s^{\prime} +\gamma_2 \mathbf{I}_{m} \bullet \mathbf{Q}^{\prime} + \sqrt{\gamma_1} \left \| \mathbf{q}^{\prime} \right \|_2 
    \leq s^* +\gamma_2 \mathbf{I}_{m} \bullet \mathbf{Q}^* + \sqrt{\gamma_1} \left \| \mathbf{q}^* \right \|_2. \label{Cons:Prop2-10}
\end{align}
Combining \eqref{Cons:Prop2-9} and \eqref{Cons:Prop2-10} leads to 
\begin{align*}
    s^{\prime} + \gamma_2 \mathbf{I}_{m} \bullet \mathbf{Q}^{\prime} + \sqrt{\gamma_1} \left \| \mathbf{q}^{\prime} \right \|_2 
    = s^* +\gamma_2 \mathbf{I}_{m} \bullet \mathbf{Q}^* + \sqrt{\gamma_1} \left \| \mathbf{q}^* \right \|_2,
\end{align*}
which indicates that $(\mathbf{x}^{\prime}, s^{\prime}, \hat{\boldsymbol{\lambda}}^{\prime}, \mathbf{q}^{\prime}, \mathbf{Q}^{\prime})$ 
is also an optimal solution of Problem \eqref{Equ:MainProblem}. 
Meanwhile, note that ${\rm{rank}}(\mathbf{Q}^{\prime}) = {\rm{rank}}(\mathbf{V} \mathbf{Y}_{11} \mathbf{V}^{\top}) \leq \min\{{\rm{rank}}(\mathbf{V}), {\rm{rank}}(\mathbf{Y}_{11})\} \leq K$, $\boldsymbol{\delta} = \mathbf{a}_0$, and $\boldsymbol{\nu}_k = \mathbf{a}_k$ for any $k\in [K]$.
Thus, the proof is complete. 
\Halmos

\medskip

As $K$ is the number of pieces formulating the piecewise linear function $f(\mathbf{x}, \boldsymbol{\xi})$, 
Theorem \ref{Prop:lowrank} shows that the rank of $ \mathbf{Q}^{\prime} $ that optimizes Problem \eqref{Equ:MainProblem} can be small.
Note that when $K \geq m$, we have ${\rm{rank}} ( \mathbf{Q}^{\prime} ) \leq m \leq K$, thereby no need to consider this case in Theorem \ref{Prop:lowrank}. 
Therefore, given any $m_1 \in [m]$ and $\mathbf{B} \in \mathcal{B}_{m_1}$, the rank of the optimal $\mathbf{Q}_{ \textup{r} }$ in Problem \eqref{Equ:appro-1}
may also be small.
With an optimized $\mathbf{B} \in \mathbb{R}^{m \times m_1}$ in Problem \eqref{Equ:appro}, we then would like to choose a small $m_1$ and find a $\mathbf{B} \in \mathcal{B}_{m_1}$ such that $\Theta_{\textup{L}}(m_1) $ can be close to $ \Theta_{\textup{M}}(m)$.
Specifically, we have the following conclusion holds.

\begin{theorem} \label{Theo:sameoptimal}
Consider the optimal solution $(\mathbf{x}^*, s^*, \hat{\boldsymbol{\lambda}}^*, \mathbf{q}^{\prime}, \mathbf{Q}^{\prime})$ of Problem \eqref{Equ:MainProblem}, 
$S_k (\forall k \in [K])$, 
$\mathbf{V}$, 
$\boldsymbol{\delta}$, 
$\boldsymbol{\nu}_k (\forall k \in [K])$, 
and $\mathbf{Y}_{11}$ that are defined in Theorem \ref{Prop:lowrank}. 
When $m_1 \geq K$, there exists a feasible solution 
$\mathbf{B}^{\dag} = [\mathbf{V}, \mathbf{C}]$ in Problem \eqref{Equ:appro} with $\mathbf{C} \in \mathbb{R}^{m \times (m_1-K)}$ and $[\mathbf{V}, \mathbf{C}]^{\top} [\mathbf{V}, \mathbf{C}] = \mathbf{I}_{m_1}$ and given this $\mathbf{B}^{\dag}$, there exists a feasible solution
$(\mathbf{x}^{\dag} = \mathbf{x}^*, 
s^{\dag} = s^*, 
\hat{\boldsymbol{\lambda}}^{\dag} = \hat{\boldsymbol{\lambda}}^*,
\mathbf{q}_{\textup{r}}^{\dag} = (\boldsymbol{\delta}^{\top}, \mathbf{0}_{m_1-K}^{\top})^{\top}, 
\mathbf{Q}_{\textup{r}}^{\dag} = 
{\scriptsize
\begin{bmatrix}
     \mathbf{Y}_{11} & \mathbf{0}_{K \times (m_1 - K)} \\ 
     \mathbf{0}_{(m_1 - K) \times K} & \mathbf{0}_{(m_1 - K) \times (m_1 - K)} 
 \end{bmatrix}})$ 
 in Problem \eqref{Equ:appro-1} such that the corresponding objective value equals the optimal value of Problem \eqref{Equ:MainProblem}, $\Theta_{\textup{M}}(m)$.
\end{theorem}

Theorem \ref{Theo:sameoptimal} shows that we may reduce the dimensionality space of the random parameters from $m$ to $K$ while maintaining high-quality solutions. 
Meanwhile, when $m_1 \geq K$, we can always find a feasible solution of Problems \eqref{Equ:appro} and \eqref{Equ:appro-1} such that the corresponding objective value is equal to the optimal value of the original Problem \eqref{Equ:MainProblem}. 
More importantly, the SDP constraints in Problem  \eqref{Equ:appro} have smaller sizes (i.e., $m_1+1$) than those in Problem \eqref{Equ:MainProblem} (i.e., $m+1$), potentially reducing computational challenges because $K$ is usually small (e.g., $K=2$ in the distributionally robust CVaR problem in Example \ref{exam:cvar-pca}). 
We used to conjecture that this constructed feasible solution is an optimal solution of Problems \eqref{Equ:appro} and \eqref{Equ:appro-1} such that the optimal value of Problem \eqref{Equ:appro} is equal to the optimal value of Problem \eqref{Equ:MainProblem} when $m_1 \geq K$. 
Most numerical experiments (see Section \ref{Sec:experiments}) show this conjecture may be correct, but we find a counter-example.

Now we provide an example as follows to illustrate that the optimal value of Problem \eqref{Equ:appro-1} with $\mathbf{B} = \mathbf{V}$ is strictly less than the optimal value of Problem \eqref{Equ:MainProblem}, which means that the constructed feasible solution $(\mathbf{B} = \mathbf{V})$ is not optimal.

\begin{example} \label{exam:not-optimal} 
We consider an instance of Problem \eqref{Equ:MainProblem}, where
$m=n=4$, $ \gamma_1 = 1$, 
$\gamma_2 = 2$, 
$ \mathbf{A} = \mathbf{0}_{l \times m}$, 
$ \mathbf{b} = \mathbf{0}_{l}$, 
$ \boldsymbol{\mu} = \mathbf{1}_{m}$, 
$ \boldsymbol{\Sigma} = \mathbf{I}_m$,
$K = 3$, 
$y_k^0(\mathbf{x}) = 0 \ (\forall k \in [K])$, 
$y_k(\mathbf{x}) = \mathbf{W}_k \mathbf{x}\ (\forall k \in [K])$ 
with $\mathbf{W}_1 = 
{\scriptsize
\begin{bmatrix} 
1 & 0 & 0 & 0 \\
0 & 0 & 0 & 0 \\
0 & 0 & 0 & 0 \\
0 & 0 & 0 & 1 
\end{bmatrix}}$, 
$\mathbf{W}_2 = 
{\scriptsize \begin{bmatrix}
0 & 0 & 0 & 0 \\
0 & 1 & 0 & 0 \\
0 & 0 & 0 & 0 \\
0 & 0 & 0 & 2 
\end{bmatrix}}$, and 
$\mathbf{W}_3 = 
{\scriptsize \begin{bmatrix}
0 & 0 & 0 & 0 \\
0 & 0 & 0 & 0 \\ 
0 & 0 & 1 & 0 \\ 
0 & 0 & 0 & 1  
\end{bmatrix}}$, 
and $\mathcal{X} = \{\mathbf{x}\in \mathbb{R}^4 \ | \  x_1 = x_3 = x_4 = 1, x_2 \in \{ -7, 1\} \}$, then Problem \eqref{Equ:MainProblem} becomes
{\small \begin{align}
& \min\limits_{\mathbf{x} \in \mathcal{X}, s, \mathbf{q}, \mathbf{Q}} 
\left\{ s + 2 \mathbf{I}_{m} \bullet \mathbf{Q} +  \left \| \mathbf{q} \right \|_2 \ \middle| \ \begin{bmatrix}  
s -\mathbf{x}^{\top}\mathbf{W}_k^{\top}  \mathbf{1}_m &   \frac{1}{2}  \left( \mathbf{q} - \mathbf{W}_k \mathbf{x} \right)^{\top} \\ 
 \frac{1}{2} \left( \mathbf{q} -  \mathbf{W}_k \mathbf{x} \right) & \mathbf{Q}\end{bmatrix} \succeq 0, \ \forall k \in [K] \right\}. \label{Equ:Example2-primary}
\end{align}}%
Solving Problem \eqref{Equ:Example2-primary} gives the optimal value $5.9882$ with 
$\mathbf{x} = [1,1,1,1]^{\top}$, 
$\mathbf{Q} = {\scriptsize  \begin{bmatrix}
0.0911 &  -0.0558 &  -0.0354 &  -0.0558\\
   -0.0558  &  0.1115  & -0.0558   & 0.1115\\
   -0.0354  & -0.0558  &  0.0911 &  -0.0558\\
   -0.0558  &  0.1115  & -0.0558 &   0.1115
\end{bmatrix}}$, and $\rm{rank}(\mathbf{Q}) = 2$. By Theorem \ref{Prop:lowrank}, we can correspondingly obtain a feasible $\mathbf{V} = {\scriptsize  \begin{bmatrix}
    0.7071  & -0.5774 &  -0.1543\\
         0  &  0.5774 &  -0.3086\\
         0  &       0 &  0.9258\\
    0.7071  &  0.5774 &   0.1543
\end{bmatrix} }$.
Now given $m_1 = K = 3$ and $\mathbf{B} = \mathbf{V} $, 
Problem \eqref{Equ:appro-1} becomes 
{\small \begin{align}
& \min\limits_{\mathbf{x} \in \mathcal{X}, s, \mathbf{q}_{\textup{r}}, \mathbf{Q}_{\textup{r}}}
\left\{ s + 2 \mathbf{I}_{m_1} \bullet \mathbf{Q}_{\textup{r}} +  \left \| \mathbf{q}_{\textup{r}} \right \|_2  
\ \middle| \ \begin{bmatrix}  s -\mathbf{x}^{\top}\mathbf{W}_k^{\top}  \mathbf{1}_m &   \frac{1}{2}  \left(\mathbf{q}_{\textup{r}} - \mathbf{V}^{\top} \mathbf{W}_k \mathbf{x} \right)^{\top}\\ 
 \frac{1}{2} \left(\mathbf{q}_{\textup{r}} -  \mathbf{V}^{\top} \mathbf{W}_k \mathbf{x} \right) & \mathbf{Q}_{\textup{r}}\end{bmatrix} \succeq 0, \ \forall k\in[K]
\right\}. \label{Equ:Example2-appro}
\end{align}}%
Solving Problem \eqref{Equ:Example2-appro} gives the optimal value $5.1139$ with $\mathbf{x} = [1,-7,1,1]^{\top}$.
That is, the optimal value of Problem \eqref{Equ:appro-1} with $\mathbf{B} = \mathbf{V}$ is strictly less than the optimal value of Problem \eqref{Equ:MainProblem}.
\end{example}



Theorem \ref{Theo:sameoptimal} and Example \ref{exam:not-optimal} show that while the optimized dimensionality reduction maintains very high-quality solutions (mostly the optimal solutions as shown in our later numerical experiments in Section \ref{Sec:experiments}), we may still potentially lose some useful information 
that achieves the optimal solution of the original problem.
To resolve this issue, we will also derive an upper bound and a new lower bound for the optimal value of the original problem in the later sections.

Note that Problem \eqref{Equ:appro} is a nonconvex optimization problem due to the max-min operator. 
That is, we develop a low-dimensional nonconvex optimization technique 
to solve the original high-dimensional SDP problem, which can be significantly difficult to solve because of the large sizes of SDP matrices.
To further efficiently solve Problem \eqref{Equ:appro}, we first reformulate it into a bilinear SDP problem (see Proposition \ref{prop:bilinear-sdp}) under the following assumption and then propose efficient algorithms (see Section \ref{Sec:Algorithm}) to solve it.

\begin{assumption} \label{assump:convex-X}
The set $\mathcal{X}$ is convex with at least one interior point. More specifically, we consider the convex set $\mathcal{X} $ in a generic SDP form: $\mathcal{X} = \left\{\mathbf{x}\in \mathbb{R}^n \ | \  \sum_{i=1}^n (\boldsymbol{\Delta}_i x_i) + \boldsymbol{\Delta}_0 \succeq 0 \right\}$, 
where $\boldsymbol{\Delta}_i \in \mathbb{R}^{\tau \times \tau}$ is symmetric for any $i \in \{0, 1, \ldots, n \}$ and some $\tau \geq 1$.
\end{assumption}

We use $\mathbf{a}_{ij}\mathbf{x} + a_{ij}^0$ ($ \forall i \in [\tau], j \in [\tau]$) to 
denote the elements of the matrix $\sum_{i=1}^n (\boldsymbol{\Delta}_i x_i) + \boldsymbol{\Delta}_0$, 
where $ \mathbf{a}_{ij}^{\top} \in \mathbb{R}^{ n}$.
We let $y_k^0(\mathbf{x}) = \mathbf{w}_k^0 \mathbf{x} + d_k^0$ 
and $y_k(\mathbf{x}) = (\mathbf{w}_k^1 \mathbf{x} + d_k^1, \ldots, \mathbf{w}_k^m \mathbf{x} + d_k^m)^{\top} = \mathbf{W}_k \mathbf{x} + \mathbf{d}_k$ for any $k\in[K]$, 
where  $ (\mathbf{w}_k^i)^{\top} \in \mathbb{R}^{ n}$ for any $ i \in  \{0, 1, \ldots, m\}$ and $k\in[K]$, $\mathbf{W}_k \in \mathbb{R}^{m\times n}$ for any $k\in[K]$,
and $\mathbf{d}_k \in \mathbb{R}^m$ for any $k\in[K]$.
The following proposition holds.

\begin{proposition} \label{prop:bilinear-sdp} 
Under Assumption \ref{assump:convex-X}, Problem \eqref{Equ:appro} has the same optimal value as the following bilinear SDP formulation:
{\small \begin{subeqnarray} \label{Equ:bilinear}
\Theta_{\textup{L}} (m_1) =  &\max\limits_{\substack{t_k, \mathbf{p}_k, \mathbf{P}_k, \forall k \in [K],\\ \mathbf{Z},\mathbf{B}  }}   & \sum_{k=1}^K \left( t_k d_k^0 + \left( t_k \boldsymbol{\mu}^{\top}  +  \mathbf{p}_k^{\top} \left( \mathbf{U} \boldsymbol{\Lambda}^{{\frac{1}{2}}}\mathbf{B} \right)^{\top} \right) \mathbf{d}_k \right) - \sum_{i=1}^\tau \sum_{j=1}^\tau z_{ij} a_{ij}^0\\
&\rm{s.t.}
& 1-\sum_{k=1}^K t_k = 0,\ \sqrt{\gamma_1} - \left\|\sum_{k=1}^K \mathbf{p}_k \right\|_2 \geq 0,\ \gamma_2 \mathbf{I}_{m_1} - \sum_{k=1}^K \mathbf{P}_k \succeq 0, \\
&& t_k \left(\mathbf{A}\boldsymbol{\mu}-\mathbf{b}\right)^{\top} + \mathbf{p}_k^{\top}\left(\mathbf{U} \boldsymbol{\Lambda}^{{\frac{1}{2}}}\mathbf{B} \right)^{\top} \mathbf{A}^{\top}   \leq 0, \ \forall k \in [K], \\
&& \sum_{k=1}^K \left(  t_k \mathbf{w}_k^0 + \left(t_k \boldsymbol{\mu}^{\top}  +  \mathbf{p}_k^{\top} \left(\mathbf{U} \boldsymbol{\Lambda}^{{\frac{1}{2}}}\mathbf{B} \right)^{\top} \right)  \mathbf{W}_k \right) - \sum_{i=1}^\tau \sum_{j=1}^\tau z_{ij}\mathbf{a}_{ij} = 0 , \\
&& \begin{bmatrix}  t_k & \mathbf{p}_k^{\top} \\ \mathbf{p}_k & \mathbf{P}_k \end{bmatrix} \succeq 0, \ \forall k\in [K],\ \mathbf{B}^{\top} \mathbf{B} = \mathbf{I}_{m_1}, \ \mathbf{Z} \succeq 0,
\end{subeqnarray}}%
where $\mathbf{p}_k \in \mathbb{R}^{m_1}$ ($ k \in [K]$), $\mathbf{P}_k \in \mathbb{R}^{m_1\times m_1}$ ($ k \in [K]$), $\mathbf{Z} \in \mathbb{R}^{\tau \times \tau}$, $\mathbf{B}\in \mathbb{R}^{m\times m_1}$, and $z_{ij}$ is the element of the matrix $\mathbf{Z}$. In addition, $\mathbf{Z}$ is the dual variable of the constraint $\sum_{i=1}^n (\boldsymbol{\Delta}_i x_i) + \boldsymbol{\Delta}_0 \succeq 0$ in $\mathcal{X}$ and $ {\scriptsize \begin{bmatrix}  t_k & \mathbf{p}_k^{\top} \\ \mathbf{p}_k & \mathbf{P}_k \end{bmatrix} }$ ($ \forall k \in [K]$) are the dual variables of constraints \eqref{Cons:appro1}. 
\end{proposition}

As solving Problem \eqref{Equ:bilinear} may not achieve the optimal value of Problem \eqref{Equ:MainProblem}, we will further develop an upper bound for the optimal value of 
Problem \eqref{Equ:MainProblem} while \textit{closing the gap} between them in the next section.


\section{Upper Bound} \label{sec:upper-bound}

We develop an inner approximation for Problem \eqref{Equ:DROM2} by relaxing the second-moment constraint in $\mathcal{D}_{\text{M}}$ while optimizing the choice of components in $\boldsymbol{\xi}_{\text{I}}$ to be relaxed, leading to the best possible upper bound for the optimal value of Problem \eqref{Equ:DROM2}.
Specifically, given $m_1 \in [m]$, we 
build the following optimized inner approximation of Problem \eqref{Equ:DROM2}:
\begin{equation}
\Theta_{\textup{U}} (m_1) := \min\limits_{\mathbf{B} \in \mathcal{B}_{m_1} } \  \overline{\Theta} (m_1,\mathbf{B}), \label{Equ:upperbound-setform}
\end{equation}
where
\begin{equation} \label{Equ:upperbound-setform-2}
\overline{\Theta} (m_1,\mathbf{B}) = 
\min_{\mathbf{x} \in \mathcal{X}} \ \max_{\mathbb{P}_{\text{I}} \in \mathcal{D}_{\textup{U}}} \ \mathbb{E}_{\mathbb{P}_{\text{I}}} \left[ f\left(\mathbf{x},\mathbf{U}\boldsymbol{\Lambda}^{{\frac{1}{2}}} \boldsymbol{\xi}_{\textup{I}}+\boldsymbol{\mu}\right) \right] 
\end{equation}
with
\begin{align}
    \mathcal{D}_{\textup{U}}(\mathcal{S}_\textup{I},\gamma_1,\gamma_2)=\left\{ \mathbb{P}_{\text{I}} \ \middle| \ \begin{array}{l}
	\mathbb{P}_{\text{I}} \left( \boldsymbol{\xi}_{\textup{I}} \in \mathcal{S}_\textup{I} \right) = 1, \ \ 	\mathbb{E}_{\mathbb{P}_{\text{I}}} \left[ \boldsymbol{\xi}_{\textup{I}}^{\top} \right] \mathbb{E}_{\mathbb{P}_{\text{I}}} \Big[ \boldsymbol{\xi}_{\textup{I}} \Big] \leq\gamma_1 \\
	\mathbb{E}_{\mathbb{P}_{\text{I}}} \left[ \mathbf{B}^{\top}\boldsymbol{\xi}_{\textup{I}} \boldsymbol{\xi}_{\textup{I}}^{\top}\mathbf{B} \right] \preceq \gamma_2\mathbf{I}_{m_1} \end{array}\right\}. \label{eqn:upperbound-ambiguity-set}
\end{align}
The second-moment constraint in $\mathcal{D}_{\textup{U}}$ is relaxed from $\mathbb{E}_{\mathbb{P}_\text{I}} [ \boldsymbol{\xi}_{\text{I}} \boldsymbol{\xi}_{\text{I}}^{\top} ] \preceq \gamma_2\mathbf{I}_{m}$.
Intuitively, the feasible region defined by the ambiguity set $ \mathcal{D}_{\textup{U}}$
is larger than that by $ \mathcal{D}_{\textup{M}}$. 
Therefore, we have several conclusions based on this new ambiguity set 
$ \mathcal{D}_{\textup{U}}$.


\begin{theorem} \label{Theo:upperbound}
The following three conclusions hold:
(i) Problem \eqref{Equ:upperbound-setform} provides an upper bound for the optimal value of Problem \eqref{Equ:DROM2}, i.e., $\Theta_{\textup{U}}(m_1) \geq \Theta_{\textup{M}}(m)$ for any $m_1 \leq m$; 
(ii) the optimal value of Problem \eqref{Equ:upperbound-setform} is nonincreasing in $m_1$, i.e., $\Theta_{\textup{U}}(m_1) \geq \Theta_{\textup{U}}(m_2)$ for any $m_1 < m_2 \leq m$; 
and (iii) when $m_1 = m$, Problem \eqref{Equ:upperbound-setform} and Problem \eqref{Equ:DROM2} have the same optimal value, i.e., $\Theta_{\textup{U}}(m) = \Theta_{\textup{M}}(m)$. 
\end{theorem}

Theorem \ref{Theo:upperbound} shows that Problem \eqref{Equ:upperbound-setform} provides a valid upper bound for the optimal value of Problem \eqref{Equ:DROM2}, $\Theta_{\textup{M}}(m)$, and the upper bound is closer to $\Theta_{\textup{M}}(m)$ if less information is relaxed in $ \mathcal{D}_{\textup{U}} $.

\begin{proposition} \label{prop:ub-sdp-equivalent}
Under Assumption \ref{Ass:piecewise}, Problem \eqref{Equ:upperbound-setform-2} has the same optimal value as the following SDP formulation:
{\small \begin{subeqnarray}  \label{Equ:upperbound-sdpform}
\overline{\Theta} (m_1,\mathbf{B}) =  & \min\limits_{ \substack{\mathbf{x}, s, {\hat{\boldsymbol{\lambda}}}, \\ \mathbf{q}, \mathbf{Q}_{ \textup{r} }, \hat{\mathbf{u}}} } & s + \gamma_2 \mathbf{I}_{m_1} \bullet \mathbf{Q}_{ \textup{r} }  + \sqrt{\gamma_1}\left \| \mathbf{q} \right \|_2 \slabel{Obj:upperbound-sdpform}\\
&\textnormal{s.t.} &   \begin{bmatrix}  s-y_k^0(\mathbf{x})-\boldsymbol{\lambda}_k^{\top}\mathbf{b}-y_k(\mathbf{x})^{\top}\boldsymbol{\mu} +\boldsymbol{\lambda}_k^{\top}\mathbf{A}\boldsymbol{\mu} & \hspace{0.1 in} \frac{1}{2} \mathbf{u}_k^{\top} \\ 
	 \frac{1}{2} \mathbf{u}_k & \mathbf{Q}_{ \textup{r} } \end{bmatrix} \succeq 0, \ \forall k\in[K], \slabel{Cons:upperbound-sdpform-1} \\
&& \mathbf{q} +  \left(\mathbf{U}\boldsymbol{\Lambda}^{{\frac{1}{2}}}  \right)^{\top}\left(\mathbf{A}^{\top}\boldsymbol{\lambda}_k-y_k(\mathbf{x})\right) = \mathbf{B} \mathbf{u}_k, \ \forall k \in [K], \slabel{Cons:upperbound-sdpform-2} \\
&& \mathbf{x} \in \mathcal{X},  \ \mathbf{q}  \in \mathbb{R}^{m}, \ \mathbf{Q}_{ \textup{r} } \in \mathbb{R}^{ m_1 \times m_1 }, \slabel{Cons:upperbound-sdpform-2.5} \\
&& \hat{\boldsymbol{\lambda}} = \left\{\boldsymbol{\lambda}_1, \dots, \boldsymbol{\lambda}_K \right\}, \ \boldsymbol{\lambda}_k \in \mathbb{R}_+^{l}, \ \hat{\mathbf{u}} = \left\{\mathbf{u}_1, \dots, \mathbf{u}_K \right\}, \ \mathbf{u}_k \in \mathbb{R}^{m_1}, \ \forall k\in[K]. \slabel{Cons:upperbound-sdpform-3}
\end{subeqnarray}}%
\end{proposition}

Proposition \ref{prop:ub-sdp-equivalent} shows that Problem \eqref{Equ:upperbound-setform} can be updated by replacing its inner optimization Problem \eqref{Equ:upperbound-setform-2} with Problem \eqref{Equ:upperbound-sdpform}.
With the updated Problem \eqref{Equ:upperbound-setform}, we have the following conclusion.

\begin{theorem} \label{Theo:upperboundSameoptimal}
Consider the optimal solution $(\mathbf{x}^*, s^*, \hat{\boldsymbol{\lambda}}^*, \mathbf{q}^{\prime}, \mathbf{Q}^{\prime})$ of Problem \eqref{Equ:MainProblem}, $S_k\ (\forall k \in [K])$, $\mathbf{V}$, $\boldsymbol{\delta}$, $\mathbf{Y}_{11}$, and $\boldsymbol{\nu}_k \ (\forall k \in [K])$ that are defined in Theorem \ref{Prop:lowrank}.
If $m_1 \geq K$, then 
$\Theta_{\textup{U}}(m_1) = \Theta_{\textup{M}}(m)$. 
Specifically, when $m_1 = K$, there exist optimal $\mathbf{B} = \mathbf{V}$
and $\mathbf{Q}_{\textup{r}} = \mathbf{Y}_{11}$ in Problem \eqref{Equ:upperbound-setform} such that $\Theta_{\textup{U}}(m_1) = \overline{\Theta} (m_1,\mathbf{V}) $.
\end{theorem}


Theorem \ref{Theo:upperboundSameoptimal} shows that when $m_1 \geq K$, the optimal value of Problem \eqref{Equ:upperbound-setform} is always equal to the optimal value of the original problem, $\Theta_{\textup{M}}(m)$. 
We may interpret the insights as follows. 
Comparing the inner-approximation Problem \eqref{Equ:upperbound-setform} and the original
Problem \eqref{Equ:DROM2}, we can notice that they differ only in the second-moment constraints in their ambiguity sets.
When $m_1 = K$, we relax the original second-moment constraint from 
$\mathbb{E}_{\mathbb{P}_\text{I}} [ \boldsymbol{\xi}_{\text{I}} \boldsymbol{\xi}_{\text{I}}^{\top} ] \preceq \gamma_2\mathbf{I}_{m}$ to 
$ \mathbf{B}^{\top}	\mathbb{E}_{\mathbb{P}_{\text{I}}} [ \boldsymbol{\xi}_{\textup{I}} \boldsymbol{\xi}_{\textup{I}}^{\top}] \mathbf{B}  \preceq \gamma_2 \mathbf{I}_{K} $
with $\mathbf{B} \in \mathcal{B}_{K}$.
That is, this relaxation eventually does not lead to a different optimal value. 
Specifically, under a worst-case distribution $\mathbb{P}_{\text{I}}^*$ generated by solving 
Problem \eqref{Equ:DROM2}, 
we have
$\mathbb{E}_{\mathbb{P}_\text{I}^*} [ \boldsymbol{\xi}_{\text{I}} \boldsymbol{\xi}_{\text{I}}^{\top} ] \preceq \gamma_2\mathbf{I}_{m}$ may be equivalent to 
$ \mathbf{B}^{\top}	\mathbb{E}_{\mathbb{P}_{\text{I}}^*} [ \boldsymbol{\xi}_{\textup{I}} \boldsymbol{\xi}_{\textup{I}}^{\top}] \mathbf{B}  \preceq \gamma_2 \mathbf{I}_{K} $.
Such equivalence largely depends on the property provided by Theorem \ref{Prop:lowrank}, which states that the rank of an optimal solution of $ \mathbf{Q} $ of Problem \eqref{Equ:MainProblem} is 
not larger than $K$.
Note that the variable $ \mathbf{Q} $ in  Problem \eqref{Equ:MainProblem} is a dual variable with respect to the second-moment constraint $\mathbb{E}_{\mathbb{P}_\text{I}} [ \boldsymbol{\xi}_{\text{I}} \boldsymbol{\xi}_{\text{I}}^{\top} ] \preceq \gamma_2\mathbf{I}_{m}$, indicating that the rank of $\mathbb{E}_{\mathbb{P}_\text{I}^*} [ \boldsymbol{\xi}_{\text{I}} \boldsymbol{\xi}_{\text{I}}^{\top} ]$ may not be large.
Specifically, we have the following proposition holds.


\begin{proposition} \label{Lem:VXVwithRank}
For any PSD matrix $\mathbf{X} \in \mathbb{R}^{m \times m}$ such that ${\rm{rank}}(\mathbf{X}) \leq m_1 \leq m$, we have the following equivalence holds:
\[ 
\mathbf{X} \preceq \mathbf{I}_{m} \iff \mathbf{B}^{\top} \mathbf{X} \mathbf{B} \preceq \mathbf{I}_{m_1}, \ \forall \mathbf{B} \in \mathcal{B}_{m_1}.
\] 
\end{proposition}

\begin{corollary} \label{Lem:VXVwithRank2} 
For any PSD matrix $\mathbf{X} \in \mathbb{R}^{m\times m}$ and ${\rm{rank}}(\mathbf{X})\leq m_1 \leq m$, there exists a $\mathbf{B} \in \mathcal{B}_{m_1}$ such that $\mathbf{X} \preceq \mathbf{I}_{m}$ is equivalent to $\mathbf{B}^{\top} \mathbf{X} \mathbf{B} \preceq \mathbf{I}_{m_1}$.
\end{corollary}

In the context of solving Problem \eqref{Equ:DROM2} and its inner-approximation Problem \eqref{Equ:upperbound-setform}, 
Proposition \ref{Lem:VXVwithRank} and Corollary \ref{Lem:VXVwithRank2} show that 
there exist a worst-case distribution $\mathbb{P}_{\text{I}}^* \in \mathcal{D}_{\text{M}}$ such that the rank of $\mathbb{E}_{\mathbb{P}_\text{I}^*} [ \boldsymbol{\xi}_{\text{I}} \boldsymbol{\xi}_{\text{I}}^{\top} ]$ is not larger than $K$
and an optimal solution $\mathbf{B}^* \in \mathcal{B}_{m_1}$ such that 
$\mathbb{E}_{\mathbb{P}_\text{I}^*} [ \boldsymbol{\xi}_{\text{I}} \boldsymbol{\xi}_{\text{I}}^{\top} ] \preceq \gamma_2\mathbf{I}_{m}$ is equivalent to 
$ (\mathbf{B}^*)^{\top}	\mathbb{E}_{\mathbb{P}_{\text{I}}^*} [ \boldsymbol{\xi}_{\textup{I}} \boldsymbol{\xi}_{\textup{I}}^{\top}] \mathbf{B}^*  \preceq \gamma_2 \mathbf{I}_{K} $.
As such, even when we use a relaxed second-moment constraint, Problem \eqref{Equ:upperbound-setform} with $m_1 \geq K$ does not lose the optimality.

\section{Lower Bound Revisited} \label{sec:lower-bound-2}

Given that Problem \eqref{Equ:upperbound-setform} with $m_1 = K$ maintains the optimal value of the original Problem \eqref{Equ:DROM2}, 
we can further perform dimensionality reduction based on Problem \eqref{Equ:upperbound-setform} as we do in Section \ref{Sec:ODRapproach}, 
thereby obtaining a new lower bound for the optimal value of Problem \eqref{Equ:DROM2}.
Specifically, we consider $K \leq m$ and recall that $\mathcal{B}_K = \{\mathbf{B} \in \mathbb{R}^{m \times K} \ | \ \mathbf{B}^{\top} \mathbf{B} = \mathbf{I}_{K} \}$.
Given any $m_1 \in [K]$, we consider
\begin{equation} \label{Equ:lowerbound-setform-2}
\min\limits_{\mathbf{B} \in \mathcal{B}_K } \ \min_{\mathbf{x} \in \mathcal{X}} \ \max_{\mathbb{P}_{\text{I}} \in \mathcal{D}_{\textup{L2}}} \ \mathbb{E}_{\mathbb{P}_{\text{I}}} \left[ f\left(\mathbf{x},\mathbf{U}\boldsymbol{\Lambda}^{{\frac{1}{2}}} \boldsymbol{\xi}_{\textup{I}}+\boldsymbol{\mu}\right) \right] 
\end{equation}
with
\begin{align*}
    \mathcal{D}_{\textup{L2}}(\mathcal{S}_\textup{I},\gamma_1,\gamma_2)=\left\{ \mathbb{P}_{\text{I}} \ \middle| \ \begin{array}{l}
	\mathbb{P}_{\text{I}} \left( \boldsymbol{\xi}_{\textup{I}} \in \mathcal{S}_\textup{I} \right) = 1, \ \ 	\mathbb{E}_{\mathbb{P}_{\text{I}}} \left[ \boldsymbol{\xi}_{\textup{I}}^{\top} \right] \mathbb{E}_{\mathbb{P}_{\text{I}}} \Big[ \boldsymbol{\xi}_{\textup{I}} \Big] \leq \gamma_1 \\
	\mathbb{E}_{\mathbb{P}_{\text{I}}} \left[ \mathbf{B}_1^{\top}\boldsymbol{\xi}_{\textup{I}} \boldsymbol{\xi}_{\textup{I}}^{\top}\mathbf{B}_1 \right] \preceq \gamma_2\mathbf{I}_{m_1}, \ \
 \mathbb{E}_{\mathbb{P}_{\text{I}}} \left[ \mathbf{B}_2^{\top}\boldsymbol{\xi}_{\textup{I}} \boldsymbol{\xi}_{\textup{I}}^{\top}\mathbf{B}_2 \right] \preceq \mathbf{0}_{(K-m_1)\times (K-m_1)}
 \end{array}
 \right\}, 
\end{align*}
$\mathbf{B} = [ \mathbf{B}_1, \mathbf{B}_2 ]$, $\mathbf{B}_1 \in \mathbb{R}^{m \times m_1}$, and $\mathbf{B}_2 \in \mathbb{R}^{m \times (K-m_1)}$. 
To obtain the above ambiguity set $\mathcal{D}_{\textup{L2}}$, we shrink the ambiguity set $\mathcal{D}_{\textup{U}}$ of Problem \eqref{Equ:upperbound-setform} by replacing the second-moment constraint
$ \mathbb{E}_{\mathbb{P}_{\text{I}}} [ \mathbf{B}^{\top}	 \boldsymbol{\xi}_{\textup{I}} \boldsymbol{\xi}_{\textup{I}}^{\top} \mathbf{B} ]  \preceq \gamma_2 \mathbf{I}_{K} $
with $\mathbb{E}_{\mathbb{P}_{\text{I}}} [ \mathbf{B}_1^{\top}\boldsymbol{\xi}_{\textup{I}} \boldsymbol{\xi}_{\textup{I}}^{\top}\mathbf{B}_1 ] \preceq \gamma_2\mathbf{I}_{m_1}$ and 
$\mathbb{E}_{\mathbb{P}_{\text{I}}} [ \mathbf{B}_2^{\top}\boldsymbol{\xi}_{\textup{I}} \boldsymbol{\xi}_{\textup{I}}^{\top}\mathbf{B}_2 ] \preceq \mathbf{0}_{(K-m_1)\times (K-m_1)}$.
The constraint $\mathbb{E}_{\mathbb{P}_{\text{I}}} [ \mathbf{B}_2^{\top} \boldsymbol{\xi}_{\textup{I}} \boldsymbol{\xi}_{\textup{I}}^{\top} \mathbf{B}_2 ] \preceq \mathbf{0}_{(K-m_1)\times (K-m_1)}$ implies that we project the random vector $\boldsymbol{\xi}_{\textup{I}}$ to the space spanned by $\mathbf{B}_2$ and the second-moment value of the projected random vector is fixed at $0$.
By doing so, we may slightly lose some information to characterize the distribution $\mathbb{P}_{\text{I}}$, 
but we can obtain a formulation with an even smaller size than Problem \eqref{Equ:upperbound-setform} and maintain high-quality solutions.
Specifically, the following theorem holds.

\begin{theorem} \label{thm-new-lower-bound}
Under Assumption \ref{Ass:piecewise}, by dualizing the inner maximization problem of Problem \eqref{Equ:lowerbound-setform-2}, we obtain the following SDP formulation:
{\small \begin{subeqnarray}  \label{Equ:.lowerbound-sdpform-rewrite}
\Theta_{\textup{L2}}(m_1) = & \min\limits_{ \substack{\mathbf{x}, s, {\hat{\boldsymbol{\lambda}}}, \\ \mathbf{q}, \mathbf{Q}_{ \textup{r} }^{\prime}, \hat{\mathbf{u}}^{\prime}, \hat{\mathbf{u}}^{\prime \prime}, \\ \mathbf{B}_1, \mathbf{B}_2 } } & s + \gamma_2 \mathbf{I}_{m_1} \bullet \mathbf{Q}_{ \textup{r} }^{\prime}  + \sqrt{\gamma_1}\left \| \mathbf{q} \right \|_2 \\
&\textnormal{s.t.} & \small{  \begin{bmatrix} s-y_k^0(\mathbf{x})-\boldsymbol{\lambda}_k^{\top}\mathbf{b}-y_k(\mathbf{x})^{\top}\boldsymbol{\mu} +\boldsymbol{\lambda}_k^{\top}\mathbf{A}\boldsymbol{\mu} & \hspace{0.1 in} \frac{1}{2}   (\mathbf{u}_k^{\prime})^{\top}  \\ 
\frac{1}{2}  \mathbf{u}_k^{\prime}  & \mathbf{Q}_{ \textup{r} }^{\prime} \end{bmatrix}} \succeq 0, \ \forall k\in[K],   \\
&& \mathbf{q} +  \left(\mathbf{U}\boldsymbol{\Lambda}^{{\frac{1}{2}}}  \right)^{\top}\left(\mathbf{A}^{\top}\boldsymbol{\lambda}_k-y_k(\mathbf{x})\right) = \mathbf{B}_1 \mathbf{u}_k^{\prime} +  \mathbf{B}_2 \mathbf{u}_k^{\prime \prime} ,\ \forall k \in [K] ,   \\
&& \mathbf{x} \in \mathcal{X}, \ \mathbf{q} \in \mathbb{R}^{m}, \ \mathbf{Q}_{ \textup{r} }^{\prime} \in \mathbb{R}^{ m_1 \times m_1 },   \\ 
&& \mathbf{B}_1 \in \mathbb{R}^{m \times m_1},  \mathbf{B}_2  \in \mathbb{R}^{m \times (K-m_1)}, \ [  \mathbf{B}_1,  \mathbf{B}_2 ]^{\top} [\mathbf{B}_1,  \mathbf{B}_2] = \mathbf{I}_K, \\
&& \hat{\boldsymbol{\lambda}} = \left\{\boldsymbol{\lambda}_1, \dots, \boldsymbol{\lambda}_K \right\}, \ \boldsymbol{\lambda}_k \in \mathbb{R}_+^{l}, \ \forall k \in [K], \\
&& \hat{\mathbf{u}}^{\prime} = \left\{\mathbf{u}_1^{\prime}, \dots, \mathbf{u}_K^{\prime} \right\}, \ \mathbf{u}_k^{\prime} \in \mathbb{R}^{m_1}, \ \forall k\in[K],  \\ 
&&\hat{\mathbf{u}}^{\prime \prime} = \left\{\mathbf{u}_1^{\prime \prime}, \dots, \mathbf{u}_K^{\prime \prime} \right\}, \ \mathbf{u}_k^{\prime \prime} \in \mathbb{R}^{K-m_1}, \ \forall k \in [K].
\end{subeqnarray}}%
In addition, the following three conclusions hold:
(i) Problem \eqref{Equ:.lowerbound-sdpform-rewrite} provides a lower bound for the optimal value of Problem \eqref{Equ:MainProblem}, i.e., $\Theta_{\textup{L2}}(m_1) \leq \Theta_{\textup{M}}(m)$ for any $m_1 \leq K$; 
(ii) the optimal value of Problem \eqref{Equ:.lowerbound-sdpform-rewrite} is nondecreasing in $m_1$, i.e., $\Theta_{\textup{L2}}(m_1) \leq \Theta_{\textup{L2}}(m_2)$ for any $m_1 < m_2 \leq K$; 
and (iii) when $m_1 = K$, Problem \eqref{Equ:.lowerbound-sdpform-rewrite} and Problem \eqref{Equ:MainProblem} have the same optimal value, i.e., $\Theta_{\textup{L2}}(K) = \Theta_{\textup{M}}(m)$. 
\end{theorem}

Recall that the lower bound provided by Problem \eqref{Equ:appro} may not achieve the optimal value of the original Problem \eqref{Equ:MainProblem} when reducing the dimensionality to $K$. 
However, the new lower bound provided by Problem \eqref{Equ:.lowerbound-sdpform-rewrite} achieves the optimal value of the original problem when $m_1 = K$.

\section{Efficient Algorithm} \label{Sec:Algorithm}

In Sections \ref{Sec:ODRapproach}--\ref{sec:lower-bound-2}, we provide two outer approximations (i.e., Problems \eqref{Equ:bilinear} and \eqref{Equ:.lowerbound-sdpform-rewrite}) leading to lower bounds for the optimal value of Problem \eqref{Equ:DROM2} and 
an inner approximation (i.e., Problem \eqref{Equ:upperbound-setform}) leading to an upper bound.
Although these approximations have matrices with smaller sizes than the original problem (i.e., Problem \eqref{Equ:MainProblem}), they are nonconvex with bilinear terms.
We derive Alternating Direction Method of Multipliers (ADMM) algorithms to solve them efficiently.

Note that the detailed ADMM algorithms for Problems \eqref{Equ:bilinear}, 
\eqref{Equ:upperbound-setform}, and \eqref{Equ:.lowerbound-sdpform-rewrite} are similar and 
the formulation for Problem \eqref{Equ:upperbound-setform} is simpler than those for Problems \eqref{Equ:bilinear} and \eqref{Equ:.lowerbound-sdpform-rewrite}.
We hence only introduce the algorithmic details for solving Problem \eqref{Equ:upperbound-setform} thereafter. 
Recall that Problem \eqref{Equ:upperbound-setform} can be formulated as follows:
{\small \begin{subeqnarray}  \label{Equ:upperbound-nonbilinear}
\Theta_{\textup{U}}(m_1) =  & \min\limits_{ \substack{\mathbf{B}, \mathbf{x}, s , \mathbf{q}, \mathbf{Q}_{ \textup{r} }, \\ \boldsymbol{\lambda}_k, \tilde{\mathbf{u}}_k, \mathbf{u}_k, \forall k \in [K]} } & s + \gamma_2 \mathbf{I}_{m_1} \bullet \mathbf{Q}_{ \textup{r} }  + \sqrt{\gamma_1}\left \| \mathbf{q} \right \|_2 \slabel{Obj:upperbound-nonbilinear}\\
&\textnormal{s.t.} &   \begin{bmatrix}  s-y_k^0(\mathbf{x})-\boldsymbol{\lambda}_k^{\top}\mathbf{b}-y_k(\mathbf{x})^{\top}\boldsymbol{\mu} +\boldsymbol{\lambda}_k^{\top}\mathbf{A}\boldsymbol{\mu} & \hspace{0.1 in} \frac{1}{2} \mathbf{u}_k^{\top} \\ 
	 \frac{1}{2} \mathbf{u}_k & \mathbf{Q}_{ \textup{r} } \end{bmatrix} \succeq 0, \ \forall k\in[K], \slabel{Cons:upperbound-nonbilinear-1} \\
&& \mathbf{q} +  \left(\mathbf{U}\boldsymbol{\Lambda}^{{\frac{1}{2}}}  \right)^{\top}\left(\mathbf{A}^{\top}\boldsymbol{\lambda}_k-y_k(\mathbf{x})\right) = \tilde{\mathbf{u}}_k, \ \forall k \in [K] , \slabel{Cons:upperbound-nonbilinear-2} \\
&& \mathbf{x} \in \mathcal{X}, \  \mathbf{B}^{\top} \mathbf{B} = \mathbf{I}_{m_1}, \ \mathbf{B} \in \mathbb{R}^{m\times m_1}, \  \mathbf{q}  \in \mathbb{R}^{m}, \ \mathbf{Q}_{ \textup{r} } \in \mathbb{R}^{ m_1 \times m_1 }, \slabel{Cons:upperbound-nonbilinear-3}\\ 
&& \boldsymbol{\lambda}_k \in \mathbb{R}_+^{l}, \ \mathbf{u}_k \in \mathbb{R}^{m_1}, \ \tilde{\mathbf{u}}_k \in \mathbb{R}^m, \ \forall k\in[K], \slabel{Cons:upperbound-nonbilinear-4} \\ 
&& \tilde{\mathbf{u}}_k = \mathbf{B} \mathbf{u}_k,\ \forall k \in [K]. \slabel{Cons:upperbound-nonbilinear-5}
\end{subeqnarray}}%
where constraints \eqref{Cons:upperbound-nonbilinear-5} and $\mathbf{B}^{\top} \mathbf{B} = \mathbf{I}_{m_1}$ are bilinear constraints. 
We consider the following augmented Lagrangian problem for Problem \eqref{Equ:upperbound-nonbilinear}:
{\small \begin{align}
\min\limits_{\substack{ \mathbf{x}, s,  \mathbf{q}, \mathbf{Q}_{ \textup{r} }, \\ \boldsymbol{\lambda}_k, \tilde{\mathbf{u}}_k, \mathbf{u}_k, \forall k \in [K];\\ \mathbf{B}; \ \boldsymbol{\beta}_k, \forall k \in [K] }} \left\{
s + \gamma_2 \mathbf{I}_{m_1} \bullet \mathbf{Q}_{ \textup{r} }  + \sqrt{\gamma_1}\left \| \mathbf{q} \right \|_2 + \sum_{k=1}^K  \boldsymbol{\beta}_k^{\top} \left( \tilde{\mathbf{u}}_k - \mathbf{B} \mathbf{u}_k \right)   + \sum_{k=1}^K \frac{\rho}{2} \left\|\tilde{\mathbf{u}}_k - \mathbf{B} \mathbf{u}_k \right\|_2^2 
\ \middle| \
\eqref{Cons:upperbound-nonbilinear-1}-\eqref{Cons:upperbound-nonbilinear-4}
\right\}, \label{Equ:ALProblem}
\end{align} }%
where $\boldsymbol{\beta}_k \in \mathbb{R}^m\ (\forall k\in[K])$ are Lagrangian multipliers and $\rho > 0$ is the penalty parameter.
Thus, we design Algorithm \ref{Alg:upperbound} to solve Problem \eqref{Equ:upperbound-nonbilinear}.

\begin{algorithm}[hbt!]
\SingleSpacedXI
\caption{\small ADMM for Problem \eqref{Equ:upperbound-nonbilinear}}\label{Alg:upperbound}
\scriptsize
\begin{algorithmic} 
\Initialize{$\mathbf{B}^0, \boldsymbol{\beta}_k^0, \forall k \in [K]$}
\Repeati{update $\left(\mathbf{x}, s,    \mathbf{q}, \mathbf{Q}_{ \textup{r} },\boldsymbol{\lambda}_k, \tilde{\mathbf{u}}_k, \mathbf{u}_k, \forall k \in [K] \right), \mathbf{B}$ and $\boldsymbol{\beta}_k (\forall k \in [K])$ alternatingly by} \\
~\\
Given $\mathbf{B}^i$ and $\boldsymbol{\beta}_k^i$ for any $k \in [K]$, solve Problem \eqref{Equ:ALProblem} to obtain the optimal solution $\left(\mathbf{x}, s,    \mathbf{q}, \mathbf{Q}_{ \textup{r} }, \boldsymbol{\lambda}_k, \tilde{\mathbf{u}}_k, \mathbf{u}_k, \forall k \in [K] \right)^{i+1}$; \\
~\\
Given $\left(\mathbf{x}, s, \mathbf{q}, \mathbf{Q}_{ \textup{r} }, \boldsymbol{\lambda}_k, \tilde{\mathbf{u}}_k, \mathbf{u}_k, \forall k \in [K]\right)^{i+1}$ and $\boldsymbol{\beta}_k^i$ for any $ k \in [K]$, solve Problem \eqref{Equ:ALProblem} to obtain the optimal solution $\mathbf{B}^{i+1}$; \\
~\\
$\boldsymbol{\beta}_k^{i+1} = \boldsymbol{\beta}_k^i + \rho \left(\tilde{\mathbf{u}}_k^{i+1} - \mathbf{B}^{i+1} \mathbf{u}_k^{i+1} \right),\ \forall k\in[K]$;
\Utlcon
\end{algorithmic}
\end{algorithm}

In Algorithm \ref{Alg:upperbound}, we initialize $\mathbf{B}^0 = {\scriptsize \begin{bmatrix}
 \mathbf{I}_{m_1} \\ 
 \boldsymbol{0}_{(m-m_1) \times m_1}
\end{bmatrix} }$ based on the PCA approximation in \citep{cheramin2022computationally}
and set $\boldsymbol{\beta}^0_k = \mathbf{0}$ for any $k \in [K]$. 
We terminate the iteration when the improvement (regarding the relative gap) of the objective value 
is less than $10^{-4}$. In this algorithm, given $\mathbf{B}$ and $\boldsymbol{\beta}_k$ for any 
$k \in [K]$, Problem \eqref{Equ:ALProblem} becomes a low-dimensional ($m_1+1$) SDP problem. 
Given $(\mathbf{x}, s, \mathbf{q}, \mathbf{Q}_{ \textup{r} }, \boldsymbol{\lambda}_k, \tilde{\mathbf{u}}_k, \mathbf{u}_k, \boldsymbol{\beta}_k, \forall k \in [K])$, 
Problem \eqref{Equ:ALProblem} becomes a nonconvex optimization problem, while the following 
proposition shows that it has an analytical optimal solution.
Thus, Algorithm \ref{Alg:upperbound} can be performed efficiently.


\begin{proposition} \label{Prop:algorithm}
Given $(\mathbf{x}, s, \mathbf{q}, \mathbf{Q}_{ \textup{r} }, \boldsymbol{\lambda}_k, 
\tilde{\mathbf{u}}_k, \mathbf{u}_k, \boldsymbol{\beta}_k, \forall k \in [K])$, 
Problem \eqref{Equ:ALProblem} has an optimal solution $\mathbf{B}^{*} = 
\tilde{\mathbf{U}}   \tilde{\mathbf{V}}^{\top}$, 
where $\sum_{k=1}^K ( \boldsymbol{\beta}_k \mathbf{u}_k^{\top} + \rho \tilde{\mathbf{u}}_k \mathbf{u}_k^{\top} ) 
= \tilde{\mathbf{U}}  \tilde{\boldsymbol{\Sigma}}  \tilde{\mathbf{V}}^{\top}$ for 
$\tilde{\mathbf{U}} \in \mathbb{R}^{m\times m_1}$, $\tilde{\boldsymbol{\Sigma}} \in \mathbb{R}^{m_1\times m_1}$, and $\tilde{\mathbf{V}} \in \mathbb{R}^{m_1\times m_1}$ 
by the singular value decomposition (SVD).
\end{proposition}

Note that we can develop similar ADMM algorithms to solve the augmented Lagrangian problems for Problems \eqref{Equ:bilinear} and \eqref{Equ:.lowerbound-sdpform-rewrite}. 
As an ADMM algorithm returns a near-optimal (and feasible) 
$\mathbf{B}^*$, we use this $\mathbf{B}^*$ to recover the lower or upper bound for the optimal value of Problem \eqref{Equ:DROM2} as follows. 
(i) For Problem \eqref{Equ:bilinear} (outer approximation), given the $\mathbf{B}^*$ returned by the ADMM algorithm,
we can solve a low-dimensional SDP problem \eqref{Equ:appro-1} to retrieve the lower bound.
(ii) For Problem \eqref{Equ:upperbound-setform} (inner approximation), given the $\mathbf{B}^*$ returned by the ADMM algorithm, 
we can solve a low-dimensional SDP problem \eqref{Equ:upperbound-setform} to retrieve the upper bound. 
(iii) For Problem \eqref{Equ:.lowerbound-sdpform-rewrite} (outer approximation), given the $[\mathbf{B}_1, \mathbf{B}_2]^*$ returned by the ADMM algorithm, 
we can solve a low-dimensional SDP problem \eqref{Equ:appro-1} to retrieve the lower bound. 
In addition, we can measure the quality of the $\mathbf{B}^*$ by deriving a theoretical interval in which the optimal value of Problem \eqref{Equ:DROM2} is located via the following proposition.

\begin{proposition} \label{Prop:upperbound}
Given any $m_1 \in [m]$ and $ \mathbf{B}^{\prime} \in \mathcal{B}_{m_1}$, 
we use $(\mathbf{x}^*, s^*,  \hat{\boldsymbol{\lambda}}^*, \mathbf{q}^*_{\textup{r}}, \mathbf{Q}^*_{\textup{r}})$ to denote an optimal solution of Problem \eqref{Equ:appro-1}.
Let $P = \sum_{k=1}^K (\gamma_2/4) \mathbf{I}_m \bullet  \mathbf{M}_k$ and $S = \min\{S_k, \ \forall k \in [K] \}$, where 
{\small \begin{align*}
& \mathbf{M}_k = \left(\mathbf{B}^{\prime} \mathbf{q}^*_{\textup{r}} +\left(\mathbf{U}\boldsymbol{\Lambda}^{{\frac{1}{2}}}\right)^{\top}\left( \mathbf{A}^{\top} \boldsymbol{\lambda}^*_k - y_k(\mathbf{x}^*) \right) 
 \right) 
\left( \mathbf{B}^{\prime} \mathbf{q}^*_{\textup{r}} + \left(\mathbf{U} \boldsymbol{\Lambda}^{{\frac{1}{2}}} \right)^{\top} \left( \mathbf{A}^{\top} \boldsymbol{\lambda}^*_k - y_k(\mathbf{x}^*) \right) \right)^{\top}, \ \forall k \in [K], \\
& S_k = s^* - y_k^0(\mathbf{x}^*) - \boldsymbol{\lambda}_k^{* \top} \mathbf{b} -y_k(\mathbf{x}^*)^{\top}\boldsymbol{\mu}+\boldsymbol{\lambda}_k^{* \top}\mathbf{A}\boldsymbol{\mu}, \ \forall k \in [K].
\end{align*}}%
We have
{\small \begin{equation}
    \Theta_{\textup{M}}(m) - \underline{\Theta}(m_1,\mathbf{B}^{\prime}) \leq  
    \frac{P}{S} \mathbb{1}(\sqrt{P} - S < 0) + (2 \sqrt{P} - S) \mathbb{1}(\sqrt{P} - S \geq 0). \nonumber
\end{equation}}%
\end{proposition}

Given any feasible $\mathbf{B}^{\prime}$, Problem \eqref{Equ:appro-1} is a low-dimensional SDP and 
can be solved efficiently. Thus, by solving an easy problem, we can obtain a theoretical lower bound 
$\underline{\Theta}(m_1,\mathbf{B}^{\prime})$ and upper bound 
$\underline{\Theta}(m_1,\mathbf{B}^{\prime}) + (P/S) \mathbb{1}(\sqrt{P} - S < 0) + (2 \sqrt{P} - S) 
\mathbb{1}(\sqrt{P} - S \geq 0)$ for the original Problem \eqref{Equ:DROM2} and quantify the gap between 
them in Proposition \ref{Prop:upperbound}. More specifically, we can obtain their theoretical relative gap (denoted by ``Theoretical Gap'')
as follows:
{\small \begin{align}
\text{Theoretical Gap} = \frac{ \frac{P}{S} \mathbb{1}(\sqrt{P} - S < 0) + (2 \sqrt{P} - S) \mathbb{1}(\sqrt{P} - S \geq 0) }{| \underline{\Theta}(m_1,\mathbf{B}^{\prime}) |} \times 100\%. \label{eqn:theoretical-gap}
\end{align}}%

\section{Numerical Experiments}  
\label{Sec:experiments}

We perform extensive numerical experiments to demonstrate the effectiveness of our ODR approach in solving two moment-based DRO problems: multiproduct newsvendor and production-transportation problems.
The mathematical models are implemented in MATLAB R2022a (ver. 9.12) by the modeling language CVX (ver. 2.2) and solved by the Mosek solver (ver. 9.3.20) on a PC with 64-bit Windows Operating System, an Intel(R) Xeon(R) W-2195 CPU @ 2.30GHz processor, and a 128 GB of memory. 
The time limit for each run is set at two hours. 
In Section \ref{Sec:newsvendor}, we specify the proposed inner and outer approximations under the ODR approach in the context of the multiproduct newsvendor and production-transportation problems. 
In Section \ref{Sec:results}, we report and analyze all the numerical results.

\subsection{Numerical Setup} \label{Sec:newsvendor}

\color{red}
\subsubsection{Multiproduct Newsvendor Problem}
\color{black}

In the deterministic multiproduct newsvendor problem, we consider $m$ products and the demand for each product $i \in [m]$ is $\xi_i$. 
Given the wholesale, retail, and salvage prices: 
$\mathbf{c} \in \mathbb{R}_+^m$, 
$\mathbf{v} \in \mathbb{R}_+^m$, and 
$\mathbf{g} \in \mathbb{R}_+^m$, respectively, 
we decide an ordering amount $ \mathbf{x} \in \mathbb{R}_+^m $ to minimize the total cost 
{\small \begin{align*}
f\left( \mathbf{x},\boldsymbol{\xi} \right) = & \mathbf{c}^{\top}\mathbf{x} - \mathbf{v}^{\top} \min \{\mathbf{x},\boldsymbol{\xi}\}-\mathbf{g}^{\top}\left[\mathbf{x} - \boldsymbol{\xi}\right]^+   = \left(\mathbf{c} - \mathbf{v}\right)^{\top}\mathbf{x} + \left(\mathbf{v} - \mathbf{g}\right)^{\top}\left[\mathbf{x} - \boldsymbol{\xi}\right]^+ \\
 =& \max\left\{\left(\mathbf{c}-\mathbf{v}\right)^{\top}\mathbf{x},\left(\mathbf{c}-\mathbf{g}\right)^{\top}\mathbf{x} + \left(\mathbf{g}-\mathbf{v}\right)^{\top}\boldsymbol{\xi}\right\}.
\end{align*}}%
Note that this piecewise linear function $f ( \mathbf{x},\boldsymbol{\xi} )$ has only two pieces, i.e., $K=2$.
When the demand $\boldsymbol{\xi}$ is uncertain and its probability distribution belongs to a distributional ambiguity set $\mathcal{D}_{\text{M0}}$ as defined in Section \ref{Sec:introdMomDRO},
we obtain the following DRO counterpart:
{\small \begin{equation}\label{Equ:Newsvendor}
\min\limits_{\mathbf{x} \geq 0}  \max\limits_{\mathbb{P} \in \mathcal{D}_{\text{M0}}} \  \mathbb{E}_{\mathbb{P}}\left[ \max\left\{ \left( \mathbf{c} - \mathbf{v} \right)^{\top} \mathbf{x}, \left(\mathbf{c} - \mathbf{g} \right)^{\top} \mathbf{x} + \left( \mathbf{g} - \mathbf{v} \right)^{\top} \boldsymbol{\xi} \right\} \right].
\end{equation} }%

\noindent
We can apply Proposition \ref{Prop:cheramin2022} to reformulate 
Problem \eqref{Equ:Newsvendor} and
the proposed inner and outer  approximations (i.e., Problems \eqref{Equ:bilinear}, \eqref{Equ:upperbound-setform}, and \eqref{Equ:.lowerbound-sdpform-rewrite}) to approximate it (see the details in Appendix \ref{subsec-apx:newsvendor}).

When the dimension $m$ of $\boldsymbol{\xi}$ is large, 
the original SDP reformulation of  Problem \eqref{Equ:Newsvendor} (i.e., Problem \eqref{Equ:NewsvendorSDP}) becomes very difficult to solve because of the large-scale SDP constraints. 
Nevertheless, our approximations under the ODR approach (i.e., Problems \eqref{Equ:NewsvendorDual}--\eqref{Equ:NewsvendorRLB}) have SDP matrices with very small sizes (e.g., $K+1=3$), largely reducing the computational burden while maintaining the original problem's optimal value and optimal solution $\mathbf{x}$.


\subsubsection{Production-Transportation Problem}

In the deterministic production-transportation problem, we make production and transportation decisions to 
satisfy all customer demands while minimizing the total cost.
We consider $m$ suppliers, each with normalized capacity one, and $n$ customers, each with demand $d_j$ ($ \forall j \in [n]$) such that $\sum_{j=1}^n d_j \leq m$. 
Supply $i \in [m]$ produces $x_i$ goods with unit production cost $ c_i $ and delivers $z_{ij}$ goods to customer $j \in [n]$ with unit transportation cost $\xi_{ij}$.
We can formulate this problem as follows:
{\small  \begin{subeqnarray} \label{Equ:determinedTransportation}
\min\limits_{\mathbf{x}, \mathbf{z}} && \mathbf{c}^{\top} \mathbf{x} + \boldsymbol{\xi}^{\top} \mathbf{z} \\
\text{s.t.}
&& 0 \leq x_i \leq 1, \ \forall i \in [m], \slabel{Cons:determined_transportation_x} \\
&&  \sum_{i=1}^m z_{ij} = d_j, \ \forall j \in [n], 
\ \ \  \sum_{j=1}^n z_{ij} = x_i, \ \forall i \in [m], 
\ \ \  z_{ij} \geq 0, \ \forall i \in [m], \ j \in [n], \slabel{Cons:determined_transportation_z1}
\end{subeqnarray} }%
where $\boldsymbol{\xi} = (\xi_{ij}, \ \forall i \in [m], \ j \in [n])^{\top}$ and $\mathbf{z}_{ij} = (z_{ij}, \ \forall i \in [m], \ j \in [n])^{\top}$. 
We consider $\boldsymbol{\xi}$ is random (e.g., see \citealt{bertsimas2010models} and \citealt{cheramin2022computationally}) and its probability distribution $\mathbb{P}$ belongs to $\mathcal{D}_{\textup{M0}}$. 
Thus, we can derive the two-stage DRO counterpart:
{\small \begin{align}
\min\limits_{\mathbf{x}} \left\{ \mathbf{c}^{\top} \mathbf{x} + \max\limits_{\mathbb{P} \in \mathcal{D}_{\textup{M0}}} \mathbb{E}_{\mathbb{P}} \left[ \mathcal{U}\left(\mathcal{Q}\left(\mathbf{x}, \boldsymbol{\xi}\right)\right) \right] \ \middle| \ \eqref{Cons:determined_transportation_x} \right\}, \label{Equ:two-stage-dro-transportation}
\end{align} }%
where $\mathcal{U}(\mathcal{Q}(\mathbf{x}, \boldsymbol{\xi})) 
= \max_{k \in [K]} \{ \alpha_k \mathcal{Q}(\mathbf{x}, \boldsymbol{\xi}) + \beta_k \}$ 
is a convex nondecreasing disutility function that incorporates
risk attitudes into the second-stage cost and $\mathcal{Q}(\mathbf{x}, \boldsymbol{\xi}) = \min_{\mathbf{z}} \{ \boldsymbol{\xi}^{\top} \mathbf{z} \ | \ \eqref{Cons:determined_transportation_z1} \} $. 
More specifically, we can reformulate \eqref{Equ:two-stage-dro-transportation} as follows:
{\small \begin{align}
\min\limits_{\mathbf{0} \leq \mathbf{x} \leq \mathbf{1}} \max\limits_{\mathbb{P}_{\text{I}} \in \mathcal{D}_{\textup{M}}} \mathbb{E}_{\mathbb{P}} \left[ \mathbf{c}^{\top} \mathbf{x} +  \mathcal{U} \left( \mathcal{Q} \left(\mathbf{x}, \mathbf{U} \boldsymbol{\Lambda}^{{\frac{1}{2}}} \boldsymbol{\xi}_{\text{I}} + \boldsymbol{\mu} \right) \right) \right]. \label{Equ:two-stage-dro-transportation-1}
\end{align} }%


We can further apply our inner and outer  approximations to approximate 
Problem \eqref{Equ:two-stage-dro-transportation-1} (see the details in Appendix \ref{subsec-apx:production-transportation}). 
When the dimension $mn$ of $\boldsymbol{\xi}$ is large, the original SDP reformulation of Problem \eqref{Equ:two-stage-dro-transportation-1} (i.e., Problem \eqref{Equ:one-stage-dro-transportation-sdp}) becomes very difficult to solve due to the high-dimensional SDP constraints, while our approximations (i.e., Problems \eqref{Equ:transportation-dual}--\eqref{Equ:transportation-LBR}) have SDP matrices with very small sizes (e.g., $K+1$) and can be solved efficiently.

\subsection{Numerical Results} \label{Sec:results}

We compare the performance of our ODR approach (that provides two lower bounds and one upper bound) with three benchmark approaches:
(i) the Mosek solver with default settings, which can provide the optimal value of the original problem; 
(ii) The low-rank algorithm proposed by \cite{burer2003nonlinear} to solve the SDP reformulation of the original problem, i.e., Problem \eqref{Equ:MainProblem}, generating a lower bound for the optimal value of Problem \eqref{Equ:MainProblem};
and 
(iii) The existing PCA approximation proposed by \cite{cheramin2022computationally}, generating PCA-based lower and upper bounds.
For the third benchmark, we consider the reduced dimension $m_1 \in \{100\% \times {\rm{dim}}(\boldsymbol{\xi}), 80\%\times {\rm{dim}}(\boldsymbol{\xi}), 60\%\times {\rm{dim}}(\boldsymbol{\xi}), 40\%\times {\rm{dim}}(\boldsymbol{\xi}), 20\%\times {\rm{dim}}(\boldsymbol{\xi}), K\}$, where ${\rm{dim}}(\boldsymbol{\xi})$ represents
the dimension of the random vector $\boldsymbol{\xi}$. 
Note that $K = 2$ in the multiproduct newsvendor problem. For the production-transportation problem, we consider $K \in \{5, 10, 15\}$.
Our proposed inner and outer approximations under the ODR approach are solved using Algorithm \ref{Alg:upperbound}, providing near-optimal $\mathbf{B}^*$ and recovering valid lower and upper bounds (see Section \ref{Sec:Algorithm}).

\subsubsection{Instance Generation and Table Header Description}

We consider various instances of the multiproduct newsvendor and production-transportation problems.
In the multiproduct newsvendor problem, the mean and standard deviation of $\boldsymbol{\xi}$ are randomly generated from the intervals $[0,10]$ and $[1,2]$, respectively.
We further generate a correlation matrix randomly using the MATLAB function ``gallery(`randcorr',n)" and then convert it to a covariance matrix.
We follow \cite{xu2018distributionally} to set the wholesale, retail, and savage prices as $c_i = 0.1 ( 5+i-1 )$, $v_i = 0.15 (5+i-1 )$, and $g_i = 0.05 (5+i-1 )$ for any $i \in [m]$, respectively. 
Meanwhile, we consider $ m \in \{100,200,400,800,1200,1600,2000\}$ in this problem.

In the production-transportation problem, we follow \cite{bertsimas2010models} to randomly generate the locations of $m$ suppliers and $n$ customers from a unit square and use $\xi_{ij}$ to measure the distance between supplier $i \in [m]$ and customer $j \in [n]$. 
For any $i \in [m]$ and $j \in [n]$, we generate $10,000$ samples of $\xi_{ij}$ from the interval $[0.5\bar{\xi}_{ij}, 1.5\bar{\xi}_{ij}]$ and use them to estimate the mean, standard deviation, and covariance matrix of $\boldsymbol{\xi}$. 
Using $\gamma_0$ denoting the average transportation cost, we then
generate production cost $c_i$ and demand $d_j$ uniformly on the intervals $[0.5\gamma_0, 1.5\gamma_0]$ and $[0.5m/n, m/n]$, respectively, for any $i \in [m]$ and $j \in [n]$. 
The disutility function $\mathcal{U}(x) = 0.25(e^{2x}-1)$ and is approximated by a piecewise-linear function with $K$ equidistant segments on the interval $[0,1]$. 
Meanwhile, we consider $(m,n) \in \{ (4,25), (5,20), (5,40), (8,25), (10,40), (20,30), (20,40) \}$.

For any given value of $m$ in the multiproduct newsvendor problem
or $(m,n)$ in the production-transportation problem, 
we randomly generate five instances and report the average performance in Tables \ref{Table:newsvendor-averageData}--\ref{Table:production-transportation-averageData-K15}. 
Here we describe several table headers that are shared by these tables.
We use ``Mosek'' and ``Low-rank'' to represent the performance of the Mosek solver and the low-rank algorithm \citep{burer2003nonlinear}, respectively.
The abbreviations ``LB,'' ``UB,'' and ``RLB'' represent lower, upper, and revisited lower bounds, respectively.
Specifically, the labels ``ODR-LB,'' ``ODR-UB,'' and ``ODR-RLB'' denote 
the lower-bound performance after solving the first outer approximation \eqref{Equ:bilinear} with $m_1 = K$, 
the upper-bound performance after solving the inner approximation \eqref{Equ:upperbound-setform} with $m_1 = K$, 
and the other lower-bound performance after solving the second outer approximation  \eqref{Equ:.lowerbound-sdpform-rewrite} with $m_1 = K$, respectively. 
We use ``PCA-100\%,'' ``PCA-80\%,'' ``PCA-60\%,''
``PCA-40\%,'' ``PCA-20\%,'' and ``PCA-$K$'' to denote the performance of the PCA approximation by \cite{cheramin2022computationally}
when the reduced dimension $m_1$ equals $100\% \times {\rm{dim}}(\boldsymbol{\xi})$, 
$80\%\times {\rm{dim}}(\boldsymbol{\xi})$, 
$60\%\times {\rm{dim}}(\boldsymbol{\xi})$, 
$40\%\times {\rm{dim}}(\boldsymbol{\xi})$, $20\%\times {\rm{dim}}(\boldsymbol{\xi})$, and $ K $, respectively. 

In all the tables, we use ``Size'' to represent the value of $m$ in the multiproduct newsvendor problem
or $(m,n)$ in the production-transportation problem and ``Time'' to represent the computational time in seconds required to solve each instance. 
We use ``Gap1'' (resp. ``Gap2'') to represent the relative gap in percentage between a lower (resp. an upper) bound and the optimal value provided by the Mosek solver. 
That is,
{\small \begin{align*}
\text{Gap1} = \frac{\text{optimal value}-\text{lower bound}}{|\text{optimal value}|} \times 100\%, \ \text{Gap2} = \frac{\text{upper bound}-\text{optimal value}}{|\text{optimal value}|} \times 100\%. 
\end{align*}}%
We further use ``Interval Gap'' to represent the relative gap in percentage between a lower bound and an upper bound, i.e.,
{\small \begin{align}
\text{Interval Gap} = \frac{\text{upper bound}-\text{lower bound}}{|\text{upper bound}|} \times 100\%. \label{eqn:interval-gap}
\end{align}}%
Specifically, for both the ODR approach and the low-rank algorithm, we take the objective value of ``ODR-UB'' as the value of ``upper bound'' in \eqref{eqn:interval-gap}.
For the PCA approximation approach, the value of ``upper bound'' in \eqref{eqn:interval-gap} is provided by this approach itself. 
The value of the ``Theoretical Gap'' for each instance is defined in \eqref{eqn:theoretical-gap}. 

Finally, we use ``-'' to represent that no result can be obtained within the time limit (i.e., two hours). 
For instance, the Mosek solver cannot solve the original problem to the optimality within two hours for the newsvendor problem with $m \geq 400$.
Hence, we cannot obtain the value of ``Gap1'' for the ``Mosek,'' ``ODR-LB,'' and ``ODR-RLB'' approaches.

\begin{table}[!htbp]
\centering
\caption{Average Performance on the Newsvendor Problem}\label{Table:newsvendor-averageData}
\centering
\resizebox{0.85\textwidth}{!}{ 
\begin{tabular}{|cc|c|c|c|c|c|c|c|}
\hline\rule{0pt}{2.5ex}
& Size ($m$) & 100 & 200 & 400 & 800 & 1200 & 1600 & 2000 \\ [0.2em]
\hline\rule{0pt}{2.5ex}
Mosek & Time (secs) & 13.02 & 363.54 & - & - & - & - & -\\ [0.2em]
\hline\rule{0pt}{2.5ex}
 \multirow{3}{*}{Low-rank} &  Gap1 ($\%$) &  2.52  & 1.79  & - & - & - & - & -\\[0.2em]
& Time (secs)  &0.26 & 0.80 & 5.46 & 47.34 & 110.33 & 309.00 & 825.62 \\[0.2em]
& Interval Gap ($\%$) & 4.27 & 3.66 & 2.67 & 2.32 & 2.28 & 2.36 &  2.52 \\ [0.2em]
\hline\rule{0pt}{2.5ex}
\multirow{3}{*}{ODR-LB} & Gap1 ($\%$) &0.09  & 0.00  & - & - & - & - & - \\[0.2em]
& Time (secs) &0.77  & 0.78  & 0.83  & 0.85  & 1.13  & 2.01  & 2.54  \\ [0.2em]
& Interval Gap ($\%$) &1.81  & 1.83  & 1.44  & 1.44  & 1.56  & 1.73  & 1.96  \\ [0.2em]
& Theoretical Gap ($\%$) & 2.73 & 3.27 & 2.47 & 3.15 & 4.28 & 2.91 & 3.58  \\ [0.2em]
\hline\rule{0pt}{2.5ex}
\multirow{3}{*}{ODR-RLB} & Gap1 ($\%$) &0.03  & 0.03  & - & - & - & - & - \\[0.2em]
& Time (secs) & 1.95  & 2.60  & 4.33  & 9.75  & 20.83  & 38.36  & 56.68 \\ [0.2em]
& Interval Gap ($\%$) &1.74  & 1.86  & 1.46  & 1.46  & 1.56  & 1.78  & 1.98  \\ [0.2em]
& Theoretical Gap ($\%$) & 1.92 & 2.36 & 1.27 & 2.44 & 1.90 & 3.01 & 2.58  \\ [0.2em]
\hline\rule{0pt}{2.5ex}
\multirow{2}{*}{ODR-UB} & Gap2 ($\%$) & 1.68  & 1.80  & - & - & - & - & -\\[0.2em]
& Time (secs) & 1.95  & 2.60  & 4.33  & 9.75  & 20.83  & 38.36  & 56.68 \\ [0.2em]
& Theoretical Gap ($\%$) & 1.92 & 2.36 & 1.27 & 2.44 & 1.90 & 3.01 & 2.58  \\ [0.2em]
\hline\rule{0pt}{2.5ex}
\multirow{5}{*}{PCA-100\%} & Gap1 ($\%$) & 0.00  & 0.00  & - & - & - & - & -\\[0.2em]
& Time (secs) & 13.04 & 361.54 & - & - & - & - & -\\[0.2em]
& Gap2 ($\%$) & 0.00  & 0.00  & - & - & - & - & -\\[0.2em]
& Time (secs) & 12.99 & 361.91 & - & - & - & - & -\\[0.2em]
& Interval Gap ($\%$) &0.00  & 0.00  & - & - & - & - & -\\[0.2em]
\hline\rule{0pt}{2.5ex}
\multirow{5}{*}{PCA-80\%} & Gap1 ($\%$) & 0.50  & 0.31  & - & - & - & - & -\\[0.2em]
& Time (secs) & 5.05 & 120.72 & 3348.00 & - & - & - & -\\[0.2em]
& Gap2 ($\%$) & 12.23  & 11.13  & - & - & - & - & -\\ [0.2em]
& Time (secs) & 7.77 & 155.39 & 4793.72 & - & - & - & -\\ [0.2em]
& Interval Gap ($\%$) &14.54  & 12.90  & 13.57  & - & - & - & -\\[0.2em]
\hline\rule{0pt}{2.5ex}
\multirow{5}{*}{PCA-60\%} & Gap1 ($\%$) & 0.98  & 0.73  & - & - & - & - & -\\[0.2em]
& Time (secs) &1.44 & 28.73 & 785.63 &  - & -  &  - &  -\\[0.2em]
& Gap2 ($\%$) & 23.33  & 24.07  & - & - & - & - & -\\[0.2em]
& Time (secs) &2.29 & 44.28 & 1196.98 & - & - & - & -\\[0.2em]
& Interval Gap ($\%$) &31.76  & 32.79  & 31.87  & - & - & - & -\\[0.2em]
\hline\rule{0pt}{2.5ex}
\multirow{5}{*}{PCA-40\%} & Gap1 ($\%$) & 1.69  & 1.20  & - & - & - & - & -\\[0.2em]
& Time (secs) & 0.39 & 5.21 & 125.43 & 3351.00 & - & - & -\\[0.2em]
& Gap2 ($\%$) & 35.79  & 35.94  & - & - & - & - & -\\ [0.2em]
& Time (secs) & 0.57 & 8.03 & 177.96 & 5237.40 & - & - & -\\[0.2em]
& Interval Gap ($\%$) & 58.50  & 58.26  & 56.45  & 57.18  & - & - & -\\[0.2em]
\hline\rule{0pt}{2.5ex}
\multirow{5}{*}{PCA-20\%} & Gap1 ($\%$) &2.71  & 1.90  & - & - & - & - & -\\[0.2em]
& Time (secs) & 0.15 & 0.43 & 6.49 & 136.97 & 971.60 & 3546.30 & -\\ [0.2em]
& Gap2 ($\%$) & 47.92  & 48.19  & - & - & - & - & -\\ [0.2em]
& Time (secs) & 0.17 & 0.60 & 9.25 & 203.97 & 1340.46 & 4940.28 & -\\[0.2em]
& Interval Gap ($\%$) & 97.74  & 84.05  & 90.65  & 93.17  & 92.38  & 94.27  & -\\ [0.2em]
\hline\rule{0pt}{2.5ex}
\multirow{5}{*}{PCA-2} & Gap1 ($\%$) & 4.26  & 3.24  & - & - & - & - & -\\[0.2em]
& Time (secs) & 0.11 & 0.12 & 0.13 & 0.20 & 0.26 & 0.36 & 0.50 \\ [0.2em]
& Gap2 ($\%$) &57.60  & 59.25  & - & - & - & - & -\\[0.2em]
& Time (secs) & 0.13 & 0.14 & 0.16 & 0.22 & 0.32 & 0.46 & 0.60 \\[0.2em]
& Interval Gap ($\%$) & 147.12  & 154.40  & 141.39  & 149.18  & 146.92  & 150.96  & 153.57 \\[0.2em]
\hline
\end{tabular} 
}
\end{table}

\begin{table}[!htbp] 
\centering
\caption{Average Performance on the Production-Transportation Problem with $K=5$}\label{Table:production-transportation-averageData-K5}
\centering
\resizebox{0.87\textwidth}{!}{
\begin{tabular}{|cc|c|c|c|c|c|c|c|}
\hline\rule{0pt}{2.5ex}
& Size ($(m,n)$) & (4,25) & (5,20) & (5,40) & (8,25) & (10,40) & (20,30) & (20,40)\\ [0.2em]
\hline\rule{0pt}{2.5ex}
Mosek & Time (secs) &  56.04   & 70.10   & 1524.45  & 1612.68  & - & - & - \\[0.2em]
\hline\rule{0pt}{2.5ex}
 \multirow{3}{*}{Low-rank} &  Gap1 ($\%$) & 3.69  & 4.17  & 2.94 & 3.42 & - & - & - \\[0.2em]
& Time (secs)  & 14.87  & 15.28  & 44.82 & 41.34 & 192.21 & 710.28 & 1678.21 \\[0.2em]
& Interval Gap ($\%$) & 3.69  & 4.17  & 2.94 & 3.42 & 5.44 & 4.12 & 5.93 \\[0.2em]
\hline\rule{0pt}{2.5ex}
\multirow{4}{*}{ODR-LB} & Gap1 ($\%$) &0.14  & 0.34  & 0.22 & 0.29 & - & - & - \\[0.2em] 
& Time (secs) &4.03  & 3.63  & 6.01  & 4.82  & 12.41  & 25.03  & 57.79  \\[0.2em]
& Interval Gap ($\%$) &0.15  & 0.34  & 0.43  & 0.29  & 0.52  & 0.55  & 0.51   \\[0.2em]
& Theoretical Gap ($\%$) & 1.24 & 1.07 & 4.91 & 0.77 & 2.71 & 1.25 & 1.09  \\ [0.2em]
\hline\rule{0pt}{2.5ex}
\multirow{4}{*}{ODR-RLB} & Gap1 ($\%$) &0.01  & 0.02  & 0.01 & 0.00 & - & - & -  \\[0.2em]
& Time (secs) & 5.42  & 5.36  & 22.40  & 20.86  & 123.64  & 330.29  & 665.43   \\[0.2em]
& Interval Gap ($\%$) &0.02  & 0.02  & 0.01  & 0.01  & 0.00  & 0.00  & 0.00   \\[0.2em]
& Theoretical Gap ($\%$) & 1.13 & 1.17 & 1.80 & 0.81 & 1.22 & 1.92 & 1.32  \\ [0.2em]
\hline\rule{0pt}{2.5ex} 
\multirow{3}{*}{ODR-UB} & Gap2 ($\%$) &0.01  & 0.01  & 0.00 & 0.00 & - & - & -  \\[0.2em]
& Time (secs) & 5.48  & 5.32  & 22.41  & 20.86  & 123.48  & 329.95  & 665.35   \\[0.2em]
& Theoretical Gap ($\%$) & 1.13 & 1.17 & 1.80 & 0.81 & 1.22 & 1.92 & 1.32  \\ [0.2em]
\hline\rule{0pt}{2.5ex}
\multirow{5}{*}{PCA-100\%} & Gap1 ($\%$) & 0.00  & 0.00  & 0.00 & 0.00 & - & - & -  \\[0.2em]
& Time (secs) & 56.74   & 68.76   & 1521.32  & 1612.94  & - & - & -  \\[0.2em]
& Gap2 ($\%$) & 0.00  & 0.00  & 0.00 & 0.00 & - & - & -  \\[0.2em]
& Time (secs) & 55.49   & 68.49   & 1521.04  & 1611.67  & - & - & -  \\[0.2em]
& Interval Gap ($\%$) &0.00  & 0.00  & 0.00 & 0.00 & - & - & -  \\[0.2em]
\hline\rule{0pt}{2.5ex}
\multirow{5}{*}{PCA-80\%} & Gap1 ($\%$) & 0.52  & 0.99  & 0.33 & 1.50 & - & - & -  \\[0.2em]
& Time (secs) & 23.41  & 27.68  & 589.74  & 625.65 & - & - & -  \\[0.2em]
& Gap2 ($\%$) &3.31  & 2.18  & 1.15 & 2.76 & - & - & -  \\[0.2em]
& Time (secs) & 29.33  & 31.20  & 657.78  & 639.09 & - & - & -  \\[0.2em]
& Interval Gap ($\%$) &3.63  & 3.08  & 1.44  & 4.12 & - & - & -  \\[0.2em]
\hline\rule{0pt}{2.5ex}
\multirow{5}{*}{PCA-60\%} & Gap1 ($\%$) & 1.60  & 2.61  & 1.26 & 3.41 & - & - & -  \\[0.2em]
& Time (secs) &8.54  & 8.41  & 143.60  & 115.54 & 3087.98  & - & -  \\[0.2em]
& Gap2 ($\%$) & 9.02  & 7.33  & 3.45 & 8.62 & - & - & -  \\[0.2em]
& Time (secs) &9.04  & 8.88  & 158.95  & 152.21 & 3484.45  & - & -  \\[0.2em]
& Interval Gap ($\%$) &9.50  & 9.17  & 4.36  & 11.04 & 4.85 & - & -  \\[0.2em]
\hline\rule{0pt}{2.5ex}
\multirow{5}{*}{PCA-40\%} & Gap1 ($\%$) & 3.75  & 4.23  & 2.33 & 4.00 & - & - & -  \\[0.2em]
& Time (secs) & 2.18  & 2.11  & 30.80  & 23.44  & 512.53 & 3004.06 & - \\[0.2em]
& Gap2 ($\%$) & 14.05  & 11.58  & 6.57 & 9.96 & - & - & -  \\ [0.2em]
& Time (secs) & 2.61  & 2.37  & 37.08  & 30.60  & 633.00 & 3149.11 & -  \\[0.2em]
& Interval Gap ($\%$) & 15.42  & 14.04  & 7.94  & 12.66  & 8.52 & 5.05 & -  \\[0.2em]
\hline\rule{0pt}{2.5ex}
\multirow{5}{*}{PCA-20\%} & Gap1 ($\%$) &4.45  & 5.41  & 3.59 & 4.16 & - & - & -  \\[0.2em]
& Time (secs) & 0.74  & 0.69  & 6.21  & 5.59  & 38.03  & 205.25  & 741.71  \\[0.2em]
& Gap2 ($\%$) & 17.07  & 12.66  & 9.25 & 10.08 & - & - & -  \\[0.2em]
& Time (secs) & 1.11  & 1.01  & 7.31  & 6.26  & 64.54  & 276.63  & 883.55  \\[0.2em]
& Interval Gap ($\%$) & 18.35  & 15.95  & 11.32  & 12.91  & 9.39  & 5.22  & 4.85  \\[0.2em]
\hline\rule{0pt}{2.5ex}
\multirow{5}{*}{PCA-5} & Gap1 ($\%$) & 5.35  & 5.73  & 4.55 & 4.39 & - & - & -  \\[0.2em]
& Time (secs) & 0.38  & 0.36  & 0.52  & 0.45  & 0.66  & 0.79  & 1.00   \\[0.2em]
& Gap2 ($\%$) &17.99  & 15.35  & 12.71 & 10.08 & - & - & -  \\[0.2em]
& Time (secs) & 0.78  & 0.72  & 2.35  & 2.08  & 10.41  & 27.98  & 68.39  \\[0.2em]
& Interval Gap ($\%$) & 19.76  & 18.26  & 15.15  & 13.12  & 9.69  & 5.48  & 5.14  \\  [0.2em] 
\hline 
\end{tabular}}  
\end{table}

\begin{table}[!htbp] 
\centering
\caption{Average Performance on the Production-Transportation Problem with $K=10$}\label{Table:production-transportation-averageData-K10}
\centering
\resizebox{0.87\textwidth}{!}{
\begin{tabular}{|cc|c|c|c|c|c|c|c|}
\hline\rule{0pt}{2.5ex}
& Size ($(m,n)$) & (4,25) & (5,20) & (5,40) & (8,25) & (10,40) & (20,30) & (20,40)\\ [0.2em]
\hline\rule{0pt}{2.5ex}
Mosek & Time (secs) &  115.52   & 103.52   & 2866.54  & 2857.47  & - & - & - \\[0.2em]
\hline\rule{0pt}{2.5ex}
 \multirow{3}{*}{Low-rank} &  Gap1 ($\%$) & 3.41  & 3.92  & 3.21 & 4.08 & - & - & - \\[0.2em]
& Time (secs)  & 17.67  & 18.41  & 56.22 & 54.17 & 227.71 & 929.00 & 1907.28 \\[0.2em]
& Interval Gap ($\%$) & 3.41  & 3.92  & 3.21 & 4.08 & 5.27 & 6.2 & 4.62 \\[0.2em]
\hline\rule{0pt}{2.5ex}
\multirow{4}{*}{ODR-LB} & Gap1 ($\%$) & 0.12  & 0.19  & 0.13 & 0.25 & - & - & - \\[0.2em] 
& Time (secs) & 6.65  & 6.41  & 11.11 & 10.31 & 31.42 & 61.41 & 131.10 \\[0.2em]
& Interval Gap ($\%$) & 0.12  & 0.19  & 0.13 & 0.25 & 0.22 & 0.19 & 0.25 \\[0.2em]
& Theoretical Gap ($\%$) & 1.78  & 1.97  & 2.96 & 1.19 & 2.10 & 1.49 & 1.22 \\[0.2em]
\hline\rule{0pt}{2.5ex}
\multirow{4}{*}{ODR-RLB} & Gap1 ($\%$) & 0.00  & 0.00  & 0.00 & 0.00 & - & - & - \\[0.2em]
& Time (secs) & 11.22  & 10.65  & 40.34 & 32.95 & 122.54 & 342.80 & 748.11 \\[0.2em]
& Interval Gap ($\%$) & 0.00  & 0.00  & 0.00 & 0.00 & 0.00 & 0.00 & 0.00 \\[0.2em]
& Theoretical Gap ($\%$) & 1.86  & 2.27  & 1.74 & 1.22 & 2.10 & 1.39 & 2.30 \\[0.2em]
\hline\rule{0pt}{2.5ex} 
\multirow{3}{*}{ODR-UB} & Gap2 ($\%$) & 0.00  & 0.00  & 0.00 & 0.00 & - & - & - \\[0.2em]
& Time (secs) & 11.14  & 10.69  & 40.28 & 32.92 & 122.48 & 344.10 & 747.90 \\[0.2em]
& Theoretical Gap ($\%$) & 1.86  & 2.27  & 1.74 & 1.22 & 2.10 & 1.39 & 2.30 \\[0.2em]
\hline\rule{0pt}{2.5ex}
\multirow{5}{*}{PCA-100\%} & Gap1 ($\%$) & 0.00  & 0.00  & 0.00 & 0.00 & - & - & -  \\[0.2em]
& Time (secs) & 115.74   & 102.90   & 2852.43  & 2851.09  & - & - & -  \\[0.2em]
& Gap2 ($\%$) & 0.00  & 0.00  & 0.00 & 0.00 & - & - & -  \\[0.2em]
& Time (secs) & 116.32  & 103.70   & 2853.10  & 2859.01  & - & - & -  \\[0.2em]
& Interval Gap ($\%$) &0.00  & 0.00  & 0.00 & 0.00 & - & - & -  \\[0.2em] 
\hline\rule{0pt}{2.5ex}
\multirow{5}{*}{PCA-80\%} & Gap1 ($\%$) & 0.47  & 0.39  & 0.38 & 0.58 & - & - & -  \\[0.2em]
& Time (secs) & 52.31  & 50.48  & 1101.61  & 979.90 & - & - & -  \\[0.2em]
& Gap2 ($\%$) &0.04  & 0.04  & 0.02 & 0.02 & - & - & -  \\[0.2em]
& Time (secs) & 58.12  & 58.88  & 1392.67  & 1364.82 & - & - & -  \\[0.2em]
& Interval Gap ($\%$) &0.51  & 0.42  & 0.40  & 0.61 & - & - & -  \\[0.2em]
\hline\rule{0pt}{2.5ex}
\multirow{5}{*}{PCA-60\%} & Gap1 ($\%$) & 1.54  & 1.39  & 0.98 & 1.57 & - & - & -  \\[0.2em]
& Time (secs) &14.48  & 15.54  & 270.80  & 240.40 & 4182.33  & - & -  \\[0.2em]
& Gap2 ($\%$) & 0.96  & 0.34  & 0.34 & 0.40 & - & - & -  \\[0.2em]
& Time (secs) &19.64  & 20.83  & 455.00  & 440.12 & 4582.31  & - & -  \\[0.2em]
& Interval Gap ($\%$) &2.46  & 1.72  & 1.32  & 1.95 & 2.61 & - & -  \\[0.2em]
\hline\rule{0pt}{2.5ex}
\multirow{5}{*}{PCA-40\%} & Gap1 ($\%$) & 2.30  & 2.31  & 1.76 & 2.18 & - & - & -  \\[0.2em]
& Time (secs) & 4.19  & 4.04  & 53.67  & 45.22  & 831.61 & 4312.63 & - \\[0.2em]
& Gap2 ($\%$) & 1.85  & 1.63  & 1.10 & 0.77 & - & - & -  \\ [0.2em]
& Time (secs) & 6.19  & 5.96  & 94.74  & 97.70  & 931.48 & 4792.10 & -  \\[0.2em]
& Interval Gap ($\%$) & 4.05  & 3.87  & 2.81  & 2.92  & 5.01 & 4.95 & -  \\[0.2em]
\hline\rule{0pt}{2.5ex}
\multirow{5}{*}{PCA-20\%} & Gap1 ($\%$) &2.98  & 3.20  & 2.26 & 2.39 & - & - & -  \\[0.2em]
& Time (secs) & 2.26  & 1.91  & 7.72  & 7.02  & 52.34  & 361.81  & 1038.72  \\[0.2em]
& Gap2 ($\%$) & 2.24  & 2.44  & 1.71 & 1.07 & - & - & -  \\[0.2em]
& Time (secs) & 3.83  & 3.52  & 18.70  & 18.07  & 71.38  & 428.93  & 1391.25  \\[0.2em]
& Interval Gap ($\%$) & 5.09  & 5.49  & 3.90  & 3.42  & 6.21  & 5.48  & 4.29  \\[0.2em]
\hline\rule{0pt}{2.5ex}
\multirow{5}{*}{PCA-10} & Gap1 ($\%$) & 3.36  & 3.26  & 2.40 & 2.53 & - & - & -  \\[0.2em]
& Time (secs) & 0.83  & 0.76  & 1.19  & 1.00  & 2.01  & 1.95  & 2.54   \\[0.2em]
& Gap2 ($\%$) &2.51  & 2.48  & 2.09 & 1.11 & - & - & -  \\[0.2em]
& Time (secs) & 1.93  & 1.81  & 7.12  & 5.97  & 21.15  & 43.62  & 79.15  \\[0.2em]
& Interval Gap ($\%$) & 5.72  & 5.59  & 4.40  & 3.60  & 6.36  & 7.28  & 6.04  \\  [0.2em] 
\hline  
\end{tabular}}   
\end{table}

\begin{table}[!htbp] 
\centering
\caption{Average Performance on the Production-Transportation Problem with $K=15$}\label{Table:production-transportation-averageData-K15}
\centering
\resizebox{0.87\textwidth}{!}{
\begin{tabular}{|cc|c|c|c|c|c|c|c|}
\hline\rule{0pt}{2.5ex}
& Size ($(m,n)$) & (4,25) & (5,20) & (5,40) & (8,25) & (10,40) & (20,30) & (20,40)\\  [0.2em]
\hline\rule{0pt}{2.5ex} 
Mosek & Time (secs) &  154.22    & 159.12    & 4004.18   & 4099.14   & - & - & - \\[0.2em]
\hline\rule{0pt}{2.5ex}
 \multirow{3}{*}{Low-rank} &  Gap1 ($\%$) & 4.51  & 3.52  & 4.06 & 5.91 & - & - & - \\[0.2em]
& Time (secs)  & 21.73  & 28.14  & 74.20 & 62.18 & 331.01 & 1124.41 & 2271.82 \\[0.2em]
& Interval Gap ($\%$) & 4.51  & 3.52  & 4.06 & 5.91 & 4.22 & 5.11 & 4.83 \\[0.2em]
\hline\rule{0pt}{2.5ex}
\multirow{4}{*}{ODR-LB} & Gap1 ($\%$) &0.05  & 0.10  & 0.12 & 0.15 & - & - & - \\[0.2em] 
& Time (secs) & 13.92  & 16.56  & 22.11  & 26.16  & 43.02  & 80.51  & 168.90  \\[0.2em]
& Interval Gap ($\%$) &0.05  & 0.10  & 0.12  & 0.15  & 0.11  & 0.13  & 0.09   \\[0.2em]
& Theoretical Gap ($\%$) & 4.69 & 5.22 & 5.57 & 6.38 & 3.29 & 4.41 & 3.95  \\ [0.2em]
\hline   \rule{0pt}{2.5ex}
\multirow{4}{*}{ODR-RLB} & Gap1 ($\%$) &0.00  & 0.00  & 0.00 & 0.00 & - & - & -  \\[0.2em]
& Time (secs) & 22.60  & 21.45  & 77.18  & 63.99  & 241.03  & 689.24  & 1550.21   \\[0.2em]
& Interval Gap ($\%$) &0.00  & 0.00  & 0.00  & 0.00  & 0.00  & 0.00  & 0.00   \\[0.2em]
& Theoretical Gap ($\%$) & 1.38 & 2.61 & 1.79 & 1.73 & 2.14 & 1.94 & 2.25  \\ [0.2em]
\hline   \rule{0pt}{2.5ex}
\multirow{3}{*}{ODR-UB} & Gap2 ($\%$) &0.00  & 0.00  & 0.00 & 0.00 & - & - & -  \\[0.2em]
& Time (secs) & 22.61  & 21.44  & 76.99  & 64.08  & 149.21  & 401.63  & 878.68   \\[0.2em]
& Theoretical Gap ($\%$) & 1.38 & 2.61 & 1.79 & 1.73 & 2.14 & 1.94 & 2.25   \\ [0.2em]
\hline\rule{0pt}{2.5ex}
\multirow{5}{*}{PCA-100\%} & Gap1 ($\%$) & 0.00  & 0.00  & 0.00 & 0.00 & - & - & -  \\[0.2em]
& Time (secs) & 154.64    & 158.67   & 4024.62  & 4095.91  & - & - & -  \\[0.2em]
& Gap2 ($\%$) & 0.00  & 0.00  & 0.00 & 0.00 & - & - & -  \\[0.2em]
& Time (secs) & 156.32    & 155.64   & 4050.04  & 4059.01  & - & - & -  \\[0.2em]
& Interval Gap ($\%$) &0.00  & 0.00  & 0.00 & 0.00 & - & - & -  \\[0.2em]
\hline\rule{0pt}{2.5ex}
\multirow{5}{*}{PCA-80\%} & Gap1 ($\%$) & 0.47  & 0.53  & 0.62 & 0.58 & - & - & -  \\[0.2em]
& Time (secs) & 68.31  & 69.22  &1439.48  & 1428.65 & - & - & -  \\[0.2em]
& Gap2 ($\%$) &0.04  & 0.06  & 0.15 & 0.08 & - & - & -  \\[0.2em]
& Time (secs) & 79.55  & 80.20  & 2120.31  & 2146.05 & - & - & -  \\[0.2em]
& Interval Gap ($\%$) &0.51  & 0.57  & 0.75  & 0.62 & - & - & -  \\[0.2em]
\hline\rule{0pt}{2.5ex} 
\multirow{5}{*}{PCA-60\%} & Gap1 ($\%$) & 0.95  & 1.13  & 1.28 & 1.25 & - & - & -  \\[0.2em]
& Time (secs) &22.50  & 25.11  & 411.27  & 398.64 & -  & - & -  \\[0.2em]
& Gap2 ($\%$) & 0.22  & 0.31  & 0.47 & 0.62 & - & - & -  \\[0.2em]
& Time (secs) &33.28  & 38.11  & 658.16  & 624.21 & -  & - & -  \\[0.2em]
& Interval Gap ($\%$) &1.17  & 1.40  & 1.74  & 1.85 & - & - & -  \\[0.2em]
\hline\rule{0pt}{2.5ex}
\multirow{5}{*}{PCA-40\%} & Gap1 ($\%$) & 1.68  & 1.67  & 1.78 & 1.81 & - & - & -  \\[0.2em]
& Time (secs) & 6.61  & 8.53  & 81.20  & 80.88  & 1553.12 & - & - \\[0.2em]
& Gap2 ($\%$) & 1.65  & 2.58  & 2.54 & 1.96 & - & - & -  \\ [0.2em]
& Time (secs) & 10.21  & 12.37  & 137.98  & 155.84  & 2296.17 & - & -  \\[0.2em]
& Interval Gap ($\%$) & 3.23  & 4.27  & 4.31  & 3.76  & 4.53 & - & -  \\[0.2em]
\hline\rule{0pt}{2.5ex}
\multirow{5}{*}{PCA-20\%} & Gap1 ($\%$) &2.22  & 3.01  & 2.09 & 1.91 & - & - & -  \\[0.2em]
& Time (secs) & 4.74  & 4.98  & 13.21  & 14.20  & 81.28  & 375.81  & 1671.21  \\[0.2em]
& Gap2 ($\%$) & 3.07  & 3.61  & 2.98 & 2.48 & - & - & -  \\[0.2em]
& Time (secs) & 9.24  & 8.19  & 27.26  & 28.79  & 124.04  & 496.38  & 2517.39  \\[0.2em]
& Interval Gap ($\%$) & 5.27  & 6.59  & 5.06  & 4.22  & 5.78  & 6.29  & 5.88  \\[0.2em]
\hline\rule{0pt}{2.5ex}
\multirow{5}{*}{PCA-15} & Gap1 ($\%$) & 2.24  & 3.32  & 2.94 & 3.09 & - & - & -  \\[0.2em]
& Time (secs) & 1.91  & 2.07  & 2.52  & 2.73  & 4.25  & 4.87  & 6.27   \\[0.2em]
& Gap2 ($\%$) &3.10  & 4.12  & 3.57 & 4.01 & - & - & -  \\[0.2em]
& Time (secs) & 4.74  & 5.12  & 12.86  & 12.92  & 23.16  & 48.21  & 95.72  \\[0.2em]
& Interval Gap ($\%$) & 5.31  & 7.40  & 6.41  & 7.03  & 6.42  & 7.19  & 6.83  \\  [0.2em]  
\hline   
\end{tabular}}     
\end{table}

\subsubsection{Numerical Performance} \label{subsec:performance-results}

From Tables \ref{Table:newsvendor-averageData}--\ref{Table:production-transportation-averageData-K15}, we have the following observations.  
First, in the newsvendor problem with $m \in \{100, 200\}$ and the production-transportation problem with $(m,n) \in \{(4,25), (5,20), (5,40), (8,25)\}$, the Mosek solver solves each instance of the original problem to the optimality. Our ODR approach performs much better than the three benchmark approaches. Both the ``ODR-LB'' and ``ODR-RLB'' provide a smaller value of ``Gap1'' than the low-rank algorithm and the PCA approximation, and require a shorter computational time than the three benchmark approaches.
The ``ODR-UB'' also provides a smaller value of ``Gap2'' than the PCA approximation if $m_1 \neq 100\% \times {\rm{dim}}(\boldsymbol{\xi})$ therein and requires shorter computational time.


Specifically, the objective value of our ``ODR-LB'' reaches 
the optimal value of the original problem for some instances of the multiproduct newsvendor problem (see Table \ref{Table:newsvendor-averageData})
and provides near-optimal solutions for the production-transportation problem with ``Gap1'' less than $0.34\%$ (see  Tables \ref{Table:production-transportation-averageData-K5}--\ref{Table:production-transportation-averageData-K15}). 
More importantly, the ``ODR-LB'' reduces the computational time by up to three orders of magnitude compared to the Mosek solver.
In addition, the ``ODR-RLB'' and ``ODR-UB'' \textit{reach the optimal value of the original problem for all instances}
in Tables \ref{Table:production-transportation-averageData-K10}--\ref{Table:production-transportation-averageData-K15} and 
provide objective values that are near-optimal (within $1.8\%$ for all instances and $0.03\%$ for most instances) in Tables \ref{Table:newsvendor-averageData}--\ref{Table:production-transportation-averageData-K5}, 
while reducing the computational time significantly. 
The results also imply that our ADMM algorithms return the optimal $\mathbf{B}^*$ for most instances.

In addition, Tables \ref{Table:newsvendor-averageData}--\ref{Table:production-transportation-averageData-K15} show that our ODR approach (including ``ODR-LB,'' ``ODR-UB,'' and ``ODR-RLB'') provides a better solution in terms of the objective value than the PCA approximation 
if the reduced dimension $m_1 \leq 80\% \times {\rm{dim}}(\boldsymbol{\xi})$ in the latter approach. 
That is, even if we maintain 80\% of the random parameters corresponding to the largest eigenvalues to be uncertain in 
the PCA approximation by focusing on only their statistical information, the performance is worse than our ODR approach, 
where we \textit{optimize} the dimensionality reduction from ${\rm{dim}}(\boldsymbol{\xi})$ to $K$ (i.e., maintaining only $1\%$ of the original dimensionality size when $m=200$ for the multiproduct newsvendor problem). 
More importantly, the inner and outer approximations of our ODR approach can be solved efficiently.

Second, when the problem size is large, i.e.,
$m \geq 400$ in the newsvendor problem and $(m,n) \in \{(10,40), (20,30), (20,40)\}$ in the production-transportation problem, 
the Mosek solver cannot solve any instance of the original problem to the optimality. Our ODR approach also performs better than the three benchmark approaches. 
Tables \ref{Table:newsvendor-averageData}--\ref{Table:production-transportation-averageData-K15}  show that ``ODR-LB'' provides a smaller value of ``Interval Gap'' (within $2\%$) and 
requires a much shorter computational time than both the low-rank algorithm and the PCA approximation. For instance, when $m=1600$ in the multiproduct newsvendor problem, the low-rank algorithm and ``ODR-LB'' take $309$ and $2.01$ seconds to solve an instance of the multiproduct newsvendor problem and provide the value of ``Interval Gap'' at $2.36\%$ and $1.73\%$, respectively. 
The PCA approximation solves this instance only when the reduced dimension $m_1$ is not larger than $20\% \times m$, by which it takes $3546.3$ seconds while the solution quality is very poor, providing the value of ``Interval Gap'' at $94.27\%$. 
Tables \ref{Table:production-transportation-averageData-K5}--\ref{Table:production-transportation-averageData-K15} show that the ``ODR-RLB'' and ``ODR-UB'' \textit{reach the optimal value of the original problem for all instances} (with ``Interval Gap'' at 0); that is, our ADMM algorithms return the optimal $\mathbf{B}^*$. 
More importantly, our ODR approach is not sensitive to the value of ${\rm{dim}}(\boldsymbol{\xi})$, 
while the benchmark approaches perform much worse when ${\rm{dim}}(\boldsymbol{\xi})$ is larger. Thus, when we cannot obtain the optimal value of the original problem, the ``ODR-LB'', ``ODR-RLB'', and ``ODR-UB'' 
can be efficiently solved to provide a narrower interval that includes the optimal value than the benchmark approaches.
That is, our ODR approach can provide a near-optimal solution very efficiently for the moment-based DRO problems where other best possible benchmark approaches are struggling.

Furthermore, although we obtain near-optimal solutions by setting $m_1 = K$ in the ODR approach, the sensitivity analyses of our ODR approach in Tables \ref{Table:production-transportation-sensiAna-K5}--\ref{Table:production-transportation-sensiAna-K15} (see Appendix \ref{subsec:sensi_analysis}) 
with respect to $m_1$ also show valuable results.
Specifically, we consider the production-transportation problem, where $m_1$ takes values from $\{3, 5, 7\}$ when $K=5$, $\{8, 10, 12\}$ when $K=10$, and $\{13, 15, 17\}$ when $K=15$. 
Concerning our ODR approach (i.e., ``ODR-LB,''  ``ODR-UB,'' and ``ODR-RLB''), we note a general trend where the values of ``Gap1,'' ``Gap2,'' and ``Interval Gap'' all tend to decrease as $m_1$ increases. This trend aligns with the theoretical results in Theorems \ref{Theo:lowerbound}, \ref{Theo:upperbound}, and \ref{thm-new-lower-bound}.

\subsubsection{Numerical Insights} \label{subsec:insights}
 
In the multiproduct newsvendor problem, Tables \ref{Table:newsvendor-averageData} show that our ODR approach performs better than the PCA approximation with respect to the objective values for all the cases except that the ``PCA-100\%'' (i.e., the original problem) provides the optimal value when the problem size is small.
Note that the PCA approximation reduces the dimensionality of the random vector $\boldsymbol{\xi}$ by focusing on only the statistical information of $\boldsymbol{\xi}$, while the ODR approach integrates the dimensionality reduction with the optimization of the original problem. 
Here we would like to further demonstrate the benefits of our approach, thereby providing insights into how we can choose the value of $\mathbf{B}$ without solving the models in our ODR approach.

Consider the multiproduct newsvendor problem.
The PCA approximation chooses the random parameters corresponding to the largest eigenvalues by maximizing the expectation of $\boldsymbol{\xi}^{\top} \boldsymbol{\xi}$, i.e., the variability of $\boldsymbol{\xi}$.
Adopting the idea of our ODR approach to integrate the dimensionality reduction with the subsequent optimization problem, we can consider the objective function $f ( \mathbf{x},\boldsymbol{\xi} )$ when choosing the random parameters in $\boldsymbol{\xi}$. 
Specifically, we can maximize the variability of $(\mathbf{g}-\mathbf{v})^{\top} \boldsymbol{\xi}$, which is the only random component in $f ( \mathbf{x},\boldsymbol{\xi} )$.
By \eqref{eqn:dimen-reduction-B}, we solve the following problem to reduce the dimension from $m$ to $m_1$:
{\small \begin{align} 
\max\limits_{\mathbf{B}^{\top} \mathbf{B} = \mathbf{I}_{m_1}}  &\mathbb{E}_{\mathbb{P}} \left[ (\mathbf{g} - \mathbf{v})^{\top}   \boldsymbol{\xi}  \boldsymbol{\xi}^{\top}  (\mathbf{g} - \mathbf{v}) \right] 
\approx  \mathbb{E}_{\mathbb{P}} \left[ (\mathbf{g}-\mathbf{v})^{\top}  \left( \mathbf{U} \boldsymbol{\Lambda}^{\frac{1}{2}} \mathbf{B} \boldsymbol{\xi}_{\textup{r}} + \boldsymbol{\mu} \right) \left( \mathbf{U} \boldsymbol{\Lambda}^{\frac{1}{2}} \mathbf{B} \boldsymbol{\xi}_{\textup{r}} + \boldsymbol{\mu} \right)^{\top} (\mathbf{g}-\mathbf{v}) \right] \nonumber \\
= & \mathbb{E}_{\mathbb{P}} \left[ (\mathbf{g}-\mathbf{v})^{\top} \left( \left( \mathbf{U} \boldsymbol{\Lambda}^{\frac{1}{2}} \mathbf{B} \boldsymbol{\xi}_{\textup{r}} \right) \left( \mathbf{U} \boldsymbol{\Lambda}^{\frac{1}{2}} \mathbf{B} \boldsymbol{\xi}_{\textup{r}} \right)^{\top} + 2\mathbf{U} \boldsymbol{\Lambda}^{\frac{1}{2}} \mathbf{B} \boldsymbol{\xi}_{\textup{r}} \boldsymbol{\mu}^{\top} + \boldsymbol{\mu}\boldsymbol{\mu}^{\top} \right) (\mathbf{g}-\mathbf{v}) \right] \nonumber \\
= & \mathbb{E}_{\mathbb{P}} \left[ (\mathbf{g}-\mathbf{v})^{\top} \left( \left( \mathbf{U} \boldsymbol{\Lambda}^{\frac{1}{2}} \mathbf{B} \boldsymbol{\xi}_{\textup{r}} \right) \left( \mathbf{U} \boldsymbol{\Lambda}^{\frac{1}{2}} \mathbf{B} \boldsymbol{\xi}_{\textup{r}} \right)^{\top}  + \boldsymbol{\mu}\boldsymbol{\mu}^{\top} \right) (\mathbf{g}-\mathbf{v}) \right]. \label{eqn:integrate-obj-dimen-red}
\end{align}}%
By introducing $\mathbf{r} = (\boldsymbol{\Lambda}^{\frac{1}{2}}  \mathbf{U}^{\top} ) (\mathbf{g}-\mathbf{v})$,
Problem \eqref{eqn:integrate-obj-dimen-red} clearly has the same optimal solution as 
\begin{equation}  \label{Equ:insight}
\max\limits_{\mathbf{B}^{\top} \mathbf{B} = \mathbf{I}_{m_1}}  \mathbf{r}^{\top} \mathbf{B} \mathbf{B}^{\top} \mathbf{r}.
\end{equation}      



\begin{proposition} \label{prop:solve-int-dim-opt}
We have $\mathbf{B}^* = \begin{bmatrix}
  \mathbf{r}/ \| \mathbf{r} \|_2, \mathbf{0}_{m \times (m_1-1)}
\end{bmatrix}$ is an optimal solution of Problem \eqref{Equ:insight}.
\end{proposition}

By considering the partial feature of the original optimization problem, the optimal $\mathbf{B}^*$ of Problem \eqref{eqn:integrate-obj-dimen-red} by Proposition \ref{prop:solve-int-dim-opt} performs better than the PCA approximation that only considers statistical information of random parameters.
Note that our proposed inner and outer approximations of the ODR approach consider the complete feature of the original optimization problem and can provide an even better choice of $\mathbf{B}$.
In the multiproduct newsvendor problem with $K=2$, we can compare the $\mathbf{B}^*$ of Problem \eqref{eqn:integrate-obj-dimen-red} with the optimal $\mathbf{B}$ provided by our proposed outer approximation \eqref{Equ:bilinear} with $m_1 = K$.
Specifically, letting $m=10$, we have 
(i) the optimal value given by the PCA approximation (lower bound) with $m_1 = K$ is $-18.62$; 
(ii) the optimal value given by \eqref{eqn:integrate-obj-dimen-red} is $-17.53$ with  $\mathbf{B} = {\scriptsize \begin{bmatrix}
-0.8696 & -0.0478 & 0.3285  & -0.0930 &  -0.2762 &  0.2126&  -0.0456& -0.0034 & 0.0361 & 0.0097   \\
0 & 0 & 0& 0& 0& 0& 0& 0& 0& 0 
\end{bmatrix}^{\top}}$; 
(iii) the optimal value given by \eqref{Equ:bilinear} (lower bound) with $m_1 = K$ is $-17.38$
 with 
$\mathbf{B} =  {\scriptsize \begin{bmatrix}
    -0.8964  & -0.1886 & 0.2094 & -0.0327 & -0.2497 & 0.2215 & -0.0548 & -0.0216 & 0.0289 & 0.0104\\
   0.0143  &  0.0052 & -0.0014 & -0.0004 & 0.0034 & -0.0035 & 0.0010 & 0.0006 & -0.0003 & -0.0002
\end{bmatrix}^{\top}}$.
Clearly, our ODR approach performs the best and the value of $\mathbf{B}$ from solving \eqref{eqn:integrate-obj-dimen-red} is close to that from our ODR approach (the Frobenius norm of the difference between the two matrices is less than $0.1$).
That is, if a decision-maker does not have enough capacity to solve the approximations of our ODR approach, the decision-maker may consider partial feature of the optimization problem when reducing the dimensionality.


\section{Conclusion} \label{Sec:conclusion}     


Moment-based DRO provides a theoretical framework to integrate moment-based information from available data with optimal decision-making. 
Extensive studies have demonstrated the effectiveness of this framework in solving various industrial applications under uncertainties. 
Although moment-based DRO problems can be reformulated as SDPs that can be solved in polynomial time, solving high-dimensional SDPs is significantly challenging. 
More importantly, high-dimensional random parameters are generally involved in industrial applications, demanding efficient approaches to solve the high-dimensional SDPs in the context of moment-based DRO.

Current approaches adopt the PCA to first reduce the dimensionality of random parameters using only the statistical information and then solve the subsequent low-dimensional approximation (SDPs).
We show that performing dimensionality reduction using the components with the largest variability 
may not produce a good optimal value from the subsequent PCA approximation and it can be even worse 
than using the components with the least variability (Example \ref{exam:cvar-pca}).
Thus, we integrate the dimensionality reduction with subsequent SDP problems and hence propose an optimized dimensionality reduction (ODR) approach for the moment-based DRO (Sections \ref{Sec:ODRapproach}--\ref{sec:lower-bound-2}), aiming to drastically reduce the computational time of solving the SDP reformulations while maintaining the optimal solution of the original problem. 

We first derive an outer approximation under the ODR approach to provide a lower bound for the optimal value of the original problem (Theorem \ref{Theo:lowerbound}), where the lower bound is nondecreasing in the reduced dimension $m_1$.
We expect to choose a small $m_1$ to close the gap between the derived lower bound and the original optimal value.
To that end, we show that the rank of each SDP matrix with respect to an optimal solution of the original high-dimensional SDP reformulation is small, guiding us on how to optimize the dimensionality reduction (Theorem \ref{Prop:lowrank}). 
With this low-rank property, we observe that the derived lower bound can be close to the original optimal value (Theorem \ref{Theo:sameoptimal}) but may not reach it (Example \ref{exam:not-optimal}). 
Furthermore, we derive an inner approximation to provide an upper bound for the optimal value of the original problem (Theorem \ref{Theo:upperbound}).
More importantly, this upper bound reaches the original optimal value when the reduced dimension $m_1$ is small (Theorem \ref{Theo:upperboundSameoptimal}).
Building on this significant result, we further derive an outer approximation to provide the second lower bound for the optimal value of the original problem, where the gap between the new lower bound and the original optimal value can be closed when the reduced dimension $m_1$ is small (Theorem \ref{thm-new-lower-bound}).

The two outer and one inner approximations derived for the original problem are all low-dimensional SDPs and nonconvex with bilinear terms (Propositions \ref{prop:bilinear-sdp} and \ref{prop:ub-sdp-equivalent} and Theorem \ref{thm-new-lower-bound}).
We accordingly develop modified ADMM algorithms to solve them efficiently (Section \ref{Sec:Algorithm}).
Finally, we demonstrate the effectiveness of our ODR approach in solving multiproduct newsvendor and production-transportation problems. 
We compare our ODR approach and algorithms with three benchmark approaches: the Mosek solver, the low-rank algorithm by \cite{burer2003nonlinear}, and existing PCA approximations by \cite{cheramin2022computationally}. 
Numerical results show that our ODR approach significantly outperforms these benchmarks in computational time and solution quality. 
Our approach can obtain an optimal or near-optimal (mostly within $0.1\%$) solution and reduce the computational time by up to three orders of magnitude.
More importantly, unlike the existing approaches that become more computationally challenging when the dimension $m$ of random parameters increases, our approach is not sensitive to $m$, demonstrating significant strength in solving large-scale practical problems (Section \ref{subsec:performance-results}).
In addition, we provide insights into why our ODR approach performs better than the existing PCA approximations (Section \ref{subsec:insights}).

Our research can be further extended in various directions. 
First, this paper considers a piece-wise linear cost function in the original problem. Thus, it would be attractive to consider a more general objective function. 
Second, it would be interesting to apply our approach to more application problems to generate practical insights. 
Third, it is also of great interest to integrate the idea of dimensionality reduction into the Wasserstein DRO or two-stage stochastic programming.
Fourth, our ODR approach can be potentially generalized to solve general SDPs with certain structures. 
Thus, it would be appealing to exploit the structures of SDP constraints and apply the ODR approach to solve more general SDPs. 
We leave the above extensions for future research.



{\bibliographystyle{apalike}
\SingleSpacedXI
\setlength\bibsep{6pt}
\bibliography{sample}}

\newpage

\begin{APPENDICES}


\setcounter{table}{0}
\renewcommand{\thetable}{\Alph{section}\arabic{table}}

\section{Table of Notations}

\begin{table}[htbp]
\begin{center}
\caption{Summary of Notations}\label{tab:notations}
\centering
{\small
\begin{tabular}{ c  p{375pt} }
\toprule
Notation &  Description\\
\midrule
$\textbf{Random}$\\
$\textbf{Variables}$:\\
$\boldsymbol{\xi}$ & The random vector $\boldsymbol{\xi} \in \mathbb{R}^m$ \\ 
$\boldsymbol{\xi}_{\textup{I}}$ & The random vector $\boldsymbol{\xi}_{\textup{I}} \in \mathbb{R}^m$ obtained by the linearly transformation of $\boldsymbol{\xi}$ \\ 
$\boldsymbol{\xi}_{\textup{r}}$ & The random vector $\boldsymbol{\xi}_{\textup{r}} \in \mathbb{R}^{m_1}$ obtained by reducing the dimension of $\boldsymbol{\xi}_{\textup{I}}$ \\ 
\midrule
$\textbf{Distributions}$:\\
$\mathbb{P}$ & The probability distribution of the random vector $\boldsymbol{\xi} $ \\ 
$\mathbb{P}_{\textup{I}}$ & The probability distribution of the random vector $\boldsymbol{\xi}_{\textup{I}} $ \\ 
$\mathbb{P}_{\textup{r}}$ & The probability distribution of the random vector $\boldsymbol{\xi}_{\textup{r}} $ \\ 
\midrule
$\textbf{Decision}$\\
$\textbf{Variables}$:\\
$\mathbf{x}$ & The decision variable $\mathbf{x} \in \mathbb{R}^n$ \\ 
$s$, $\boldsymbol{\lambda}_k$, $\mathbf{q}$, $\mathbf{Q}$ & Decision variables in original SDP problem \\ 
$\hat{\boldsymbol{\lambda}}$ & $ \hat{\boldsymbol{\lambda}} :=  \{\boldsymbol{\lambda}_1, \ldots, \boldsymbol{\lambda}_K \}$\\ 
$\mathbf{q}_{\textup{r}}$, $\mathbf{Q}_{\textup{r}}$ & Decision variables in PCA approximation \\ 
$\mathbf{B}$  & The decision variable used in the optimized dimensionality reduction \\ 
$t_k, \mathbf{p}_k, \mathbf{P}_k, \mathbf{Z}$ & Decision variables used in the lower bound  \\
$\mathbf{Q}_{ \textup{r} }^{\prime}, \hat{\mathbf{u}}^{\prime}, \hat{\mathbf{u}}^{\prime \prime}, \mathbf{B}_1, \mathbf{B}_2 $ & Decision variables used in the revisited lower bound \\ 
\midrule
$\textbf{Parameters}$\\
$\textbf{and Sets}$:\\
$\mathcal{X}$ & The feasible set of decision variable $\mathbf{x}$ \\ 
$\mathcal{D}_{\textup{M0}}$ & The distributional ambiguity
set constructed by statistical information \\ 
$\mathcal{D}_{\textup{M}}$ & The distributional ambiguity
set corresponding to $\boldsymbol{\xi}_{\textup{I}}$ \\
$\mathcal{S}$ & The support of $\boldsymbol{\xi}$ \\
$\gamma_1$ & A scalar $\gamma_1 \geq 0$ \\
$\gamma_2$ & A scalar $\gamma_2 \geq 1$ \\
$\boldsymbol{\mu}$ & The estimated mean of $\boldsymbol{\xi} $ \\
$\boldsymbol{\Sigma}$ &  The estimated covariance matrix of $\boldsymbol{\xi} $\\
$\mathbf{U}$, $\mathbf{\Lambda}$ & Two matrices produced by the eigenvalue decomposition on the covariance matrix $\boldsymbol{\Sigma}$ \\
$\mathbf{A}, \mathbf{b}$ & $\mathcal{S} := \{ \boldsymbol{\xi} \ | \ \mathbf{A} \boldsymbol{\xi} \leq \mathbf{b} \}$\\ 
$\mathcal{S}_{\textup{I}}$ & The support of $\boldsymbol{\xi}_{\textup{I}}$\\ 
$\mathcal{S}_{\textup{r}}$ & The support of $\boldsymbol{\xi}_{\textup{r}}$\\ 
$\mathcal{D}_{\textup{L}}$ & The distributional ambiguity
set corresponding to $\boldsymbol{\xi}_{\textup{r}}$ \\
$\mathcal{B}_{m_1}$ & The feasible set of decision variable $\mathbf{B} \in \mathbb{R}^{m \times m_1 }$ \\
$\mathcal{D}_{\textup{U}}$ & The distributional ambiguity
set by relaxing the second-moment constraint in $\mathcal{D}_{\textup{M}}$ \\  
\midrule
$\textbf{Optimal Value}$\\
$\textbf{Functions}$:\\
$\Theta_{\textup{M}} (m)$ & The optimal value of the original problem \\
$\Theta_{\textup{M}} (m_1)$ & The optimal value of the PCA approximation \\ 
$\Theta_{\textup{L}} (m_1)$ & The optimal value of the first outer approximation \\ 
$\underline{\Theta} (m_1,\mathbf{B})$ & The optimal value of the subproblem of the first outer approximation  \\ 
$\Theta_{\textup{U}} (m_1)$ & The optimal value of the inner approximation \\ 
$\overline{\Theta} (m_1,\mathbf{B})$ & The optimal value of the subproblem of the inner approximation  \\ 
$\Theta_{\textup{L2}} (m_1)$ & The optimal value of the second outer approximation \\ 
\bottomrule
\end{tabular}}
\medskip 
\end{center}
\end{table}

\section{Supplement to Section \ref{Sec:ODRapproach}}

\subsection{Proof of Lemma \ref{Lem:3cons}}
First, we have
\begin{align*}
&\begin{bmatrix}
 \mathbf{I}_m & \mathbf{B} \\ 
 \mathbf{B}^{\top} & \mathbf{I}_{m_1}
\end{bmatrix} \succeq 0 
\iff \mathbf{I}_m - \mathbf{B} \mathbf{I}_{m_1}^{-1} \mathbf{B}^{\top} \succeq 0 
\iff \mathbf{B} \mathbf{B}^{\top} \preceq \mathbf{I}_m,
\end{align*}
where the first equivalence is by Schur complement and the second is because $\mathbf{I}_{m_1}^{-1} = \mathbf{I}_{m_1}$.

Second, we have 
\begin{align*}
 &\begin{bmatrix}
 \mathbf{I}_m & \mathbf{B} \\ 
 \mathbf{B}^{\top} & \mathbf{I}_{m_1}
\end{bmatrix} \succeq 0 
\iff \mathbf{I}_{m_1} -  \mathbf{B}^{\top} \mathbf{I}_m^{-1} \mathbf{B}   \succeq 0
\iff \mathbf{B}^{\top} \mathbf{B} \preceq \mathbf{I}_{m_1},
\end{align*}
where the first equivalence is by Schur complement and the second is because $\mathbf{I}_{m}^{-1} = \mathbf{I}_{m}$.
Thus, the lemma is proved.
\Halmos

\subsection{Proof of Lemma \ref{Lem:VXV}} 
(i) Suppose $\mathbf{X} \succeq \mathbf{Y}$. 
For any $\mathbf{a} \in \mathbb{R}^n$, we have $ \mathbf{V} \mathbf{a} \in \mathbb{R}^m $.
It follows that 
\begin{align*}
    \mathbf{X} \succeq \mathbf{Y}
     \Longrightarrow & (\mathbf{V} \mathbf{a})^{\top} (\mathbf{X}-\mathbf{Y}) (\mathbf{V} \mathbf{a}) \geq 0, \ \forall  \mathbf{a} \in \mathbb{R}^n \\
     \iff &\mathbf{a}^{\top} \left( \mathbf{V}^{\top} (\mathbf{X}-\mathbf{Y}) \mathbf{V} \right) \mathbf{a} \geq 0, \ \forall  \mathbf{a} \in \mathbb{R}^n \\
     \iff &\mathbf{V}^{\top} (\mathbf{X}-\mathbf{Y}) \mathbf{V} \succeq 0 \iff \mathbf{V}^{\top} \mathbf{X} \mathbf{V} \succeq \mathbf{V}^{\top} \mathbf{Y} \mathbf{V}.
\end{align*}

(ii) First, for any $\mathbf{V} \in \mathbb{R}^{m\times m}$, we have
\[
\mathbf{X} \succeq \mathbf{Y} \Longrightarrow \mathbf{V}^{\top} \mathbf{X} \mathbf{V} \succeq \mathbf{V}^{\top} \mathbf{Y} \mathbf{V}
\]
by (i).
Second, suppose $\mathbf{V}^{\top} \mathbf{X} \mathbf{V} \succeq \mathbf{V}^{\top} \mathbf{Y} \mathbf{V}$. Note that $\mathbf{V}^{-1} \in \mathbb{R}^{m\times m}$. 
According to (i), the matrix $\mathbf{V}^{\top} \mathbf{X} \mathbf{V} - \mathbf{V}^{\top} \mathbf{Y} \mathbf{V}$ remains as PSD if it multiplies $(\mathbf{V}^{-1})^{\top}$ before it and $\mathbf{V}^{-1}$ after it, i.e.,
\begin{align*}
    \left( \mathbf{V}^{-1} \right)^{\top} \mathbf{V}^{\top} \mathbf{X} \mathbf{V} \mathbf{V}^{-1} \succeq \left( \mathbf{V}^{-1} \right)^{\top} \mathbf{V}^{\top} \mathbf{Y} \mathbf{V} \mathbf{V}^{-1}.
\end{align*}
It follows that $\mathbf{X} \succeq \mathbf{Y}$ 
because $(\mathbf{V}^{-1})^{\top} \mathbf{V}^{\top} = \mathbf{I}_m$ and $\mathbf{V} \mathbf{V}^{-1} = \mathbf{I}_m$.
Thus, $\mathbf{X} \succeq \mathbf{Y}$ is equivalent to $\mathbf{V}^{\top} \mathbf{X} \mathbf{V} \succeq \mathbf{V}^{\top} \mathbf{Y} \mathbf{V}$ if $\mathbf{V} \in \mathbb{R}^{m\times m}$ is invertible.
\Halmos

\subsection{Proof of Theorem \ref{Theo:lowerbound} }
(i) Given any $\mathbf{x} \in \mathcal{X}$ and $\mathbf{B} \in \mathcal{B}_{m_1}$, i.e., $\mathbf{B}^{\top} \mathbf{B}  = \mathbf{I}_{m_1}$, we define $\boldsymbol{\zeta} = \mathbf{U} \boldsymbol{\Lambda}^{{\frac{1}{2}}} \mathbf{B} \boldsymbol{\xi}_{ \textup{r} } + \boldsymbol{\mu}$ and use $\mathcal{S}_{\zeta}$ and $\mathcal{D}_{\zeta}$ to denote its support and ambiguity set, respectively. 
As $ \mathcal{S}_{ \textup{r} } = \{\boldsymbol{\xi}_{ \textup{r} }\in \mathbb{R}^{m_1} \ | \ \mathbf{U} \boldsymbol{\Lambda}^{{\frac{1}{2}}}\mathbf{B}\boldsymbol{\xi}_{ \textup{r} } + \boldsymbol{\mu} \in \mathcal{S} \}$ and 
$\mathcal{S}_{\zeta} = \{\boldsymbol{\zeta} \in \mathbb{R}^{m} \ | \ \boldsymbol{\zeta} = \mathbf{U} \boldsymbol{\Lambda}^{{\frac{1}{2}}}\mathbf{B}\boldsymbol{\xi}_{ \textup{r} }+\boldsymbol{\mu}, \ \boldsymbol{\xi}_{ \textup{r} }\in \mathcal{S}_{ \textup{r} } \}$, we can deduce $\mathcal{S}_{\zeta} \subseteq  \mathcal{S}$.
We also have
{\small \begin{align}
& \left( \mathbb{E}_{\mathbb{P}_{\zeta}} \left[ \boldsymbol{\zeta} \right] - \boldsymbol{\mu} \right)^{\top} \boldsymbol{\Sigma}^{-1} \left( \mathbb{E}_{\mathbb{P}_{\zeta}} \left[ \boldsymbol{\zeta} \right] - \boldsymbol{\mu} \right) = \left( \mathbb{E}_{\mathbb{P}_{ \textup{r} } } \left[ \mathbf{U} \boldsymbol{\Lambda}^{{\frac{1}{2}}} \mathbf{B} \boldsymbol{\xi}_{ \textup{r} } \right]  \right)^{\top} \boldsymbol{\Sigma}^{-1} \mathbb{E}_{\mathbb{P}_{ \textup{r} } } \left[ \mathbf{U} \boldsymbol{\Lambda}^{{\frac{1}{2}}}\mathbf{B} \boldsymbol{\xi}_{ \textup{r} } \right] \nonumber \\
= \  & \mathbb{E}_{\mathbb{P}_{ \textup{r} } } \left[ \boldsymbol{\xi}_{ \textup{r} }^{\top} \right]  \mathbf{B}^{\top}  \left( \mathbf{U} \boldsymbol{\Lambda}^{{\frac{1}{2}}} \right)^{\top} \boldsymbol{\Sigma}^{-1} \left( \mathbf{U} \boldsymbol{\Lambda}^{{\frac{1}{2}}} \right) \mathbf{B} \mathbb{E}_{\mathbb{P}_{ \textup{r} } } \left[  \boldsymbol{\xi}_{ \textup{r} } \right] = \mathbb{E}_{\mathbb{P}_{ \textup{r} } } \left[ \boldsymbol{\xi}_{ \textup{r} }^{\top} \right]  \mathbf{B}^{\top}  \left( \mathbf{U} \boldsymbol{\Lambda}^{{\frac{1}{2}}} \right)^{\top} \left( \mathbf{U} \boldsymbol{\Lambda} \mathbf{U}^{\top} \right)^{-1} \left( \mathbf{U} \boldsymbol{\Lambda}^{{\frac{1}{2}}} \right) \mathbf{B} \mathbb{E}_{\mathbb{P}_{ \textup{r} } } \left[  \boldsymbol{\xi}_{ \textup{r} } \right] \nonumber \\
= \ & \mathbb{E}_{\mathbb{P}_{ \textup{r} } } \left[ \boldsymbol{\xi}_{ \textup{r} }^{\top} \right]  \mathbf{B}^{\top}  \left( \mathbf{U} \boldsymbol{\Lambda}^{{\frac{1}{2}}}  \right)^{\top} \left( \left( \mathbf{U} \boldsymbol{\Lambda}^{{\frac{1}{2}}} \right) \left( \mathbf{U} \boldsymbol{\Lambda}^{{\frac{1}{2}}} \right)^{\top} \right)^{-1} \left( \mathbf{U} \boldsymbol{\Lambda}^{{\frac{1}{2}}} \right) \mathbf{B} \mathbb{E}_{\mathbb{P}_{ \textup{r} } } \left[ \boldsymbol{\xi}_{ \textup{r} } \right] \nonumber \\
= \ & \mathbb{E}_{\mathbb{P}_{ \textup{r} } } \left[ \boldsymbol{\xi}_{ \textup{r} }^{\top} \right]  \mathbf{B}^{\top}  \left( \mathbf{U} \boldsymbol{\Lambda}^{{\frac{1}{2}}} \right)^{\top} \left( \left(\mathbf{U} \boldsymbol{\Lambda}^{{\frac{1}{2}}} \right)^{\top} \right)^{-1}   \left( \mathbf{U} \boldsymbol{\Lambda}^{{\frac{1}{2}}} \right)^{-1} \left( \mathbf{U} \boldsymbol{\Lambda}^{{\frac{1}{2}}} \right)\mathbf{B} \mathbb{E}_{\mathbb{P}_{ \textup{r} } } \left[  \boldsymbol{\xi}_{ \textup{r} } \right] \nonumber \\
= \ & \mathbb{E}_{\mathbb{P}_{ \textup{r} } } \left[ \boldsymbol{\xi}_{ \textup{r} }^{\top} \right]  \mathbf{B}^{\top}  \mathbf{B} \mathbb{E}_{\mathbb{P}_{ \textup{r} } } \Big[  \boldsymbol{\xi}_{ \textup{r} } \Big] 
=   \mathbb{E}_{\mathbb{P}_{ \textup{r} } } \left[ \boldsymbol{\xi}_{ \textup{r} }^{\top} \right] \mathbb{E}_{\mathbb{P}_{ \textup{r} } } \Big[ \boldsymbol{\xi}_{ \textup{r} } \Big] 
\leq  \gamma_1, \label{eqn:thm1-expect}
\end{align}}%
where the inequality holds because of \eqref{DRSP-RA:LB-amb}.
Meanwhile, we have 
\begin{align}
	&\mathbb{E}_{\mathbb{P}_{\zeta}} \left[ (\boldsymbol{\zeta} -\boldsymbol{\mu})(\boldsymbol{\zeta} - \boldsymbol{\mu})^{\top} \right] 
	=  \mathbb{E}_{\mathbb{P}_{ \textup{r} }} \left[ \mathbf{U} \boldsymbol{\Lambda}^{{\frac{1}{2}}} \mathbf{B} \boldsymbol{\xi}_{\textup{r}} \boldsymbol{\xi}_{\textup{r}}^{\top} \mathbf{B}^{\top}  \boldsymbol{\Lambda}^{{\frac{1}{2}}} \mathbf{U}^{\top} \right] \nonumber\\
	\preceq  & \mathbf{U} \boldsymbol{\Lambda}^{{\frac{1}{2}}}\mathbf{B} \gamma_2 \mathbf{I}_{m_1} \mathbf{B}^{\top}  \boldsymbol{\Lambda}^{{\frac{1}{2}}} \mathbf{U}^{\top} 
	= \gamma_2\mathbf{U} \boldsymbol{\Lambda}^{{\frac{1}{2}}}\mathbf{B} \mathbf{B}^{\top}  \boldsymbol{\Lambda}^{{\frac{1}{2}}} \mathbf{U}^{\top} 
	\preceq \gamma_2\mathbf{U}\boldsymbol{\Lambda} \mathbf{U}^{\top}
	= \gamma_2\boldsymbol{\Sigma}, \label{eqn:thm1-covariance}
\end{align}
where the first inequality holds because of \eqref{DRSP-RA:LB-amb} and the second inequality holds because $\mathbf{B}^{\top} \mathbf{B}  = \mathbf{I}_{m_1} $, leading to $ \mathbf{B}^{\top} \mathbf{B} \preceq \mathbf{I}_{m_1}$, which is further equivalent to $ \mathbf{B} \mathbf{B}^{\top}  \preceq \mathbf{I}_m$ by Lemma \ref{Lem:3cons}.
By $\mathcal{S}_{\zeta} \subseteq  \mathcal{S}$, \eqref{eqn:thm1-expect}, and \eqref{eqn:thm1-covariance}, it follows that $\mathcal{D}_{\zeta}$ lies in $\mathcal{D}_{\text{M0}}$, i.e., $\mathcal{D}_{\zeta} \subseteq \mathcal{D}_{\text{M0}}$. 

Therefore, given any $\mathbf{x} \in \mathcal{X}$ and $\mathbf{B} \in \mathcal{B}_{m_1}$, we have
\begin{align*}
    \max\limits_{\mathbb{P}_{ \textup{r} }  \in \mathcal{D}_{\text{L}} } \ \mathbb{E}_{\mathbb{P}_{\textup{r}} } \left[ f\left(\mathbf{x},\mathbf{U}\boldsymbol{\Lambda}^{{\frac{1}{2}}} \mathbf{B}  \boldsymbol{\xi}_{ \textup{r} } +\boldsymbol{\mu}\right) \right] 
    = \max\limits_{\mathbb{P}_{ \zeta }  \in \mathcal{D}_{\zeta}} \ \mathbb{E}_{\mathbb{P}_{\zeta} } \left[ f\left(\mathbf{x},\boldsymbol{\zeta}\right) \right] 
    \leq \max_{\mathbb{P} \in \mathcal{D}_{\text{M0}}} \ \mathbb{E}_{\mathbb{P}} \left[ f\left(\mathbf{x},\boldsymbol{\xi}\right) \right],
\end{align*}
where the equality holds by change of variables and the inequality holds because $\mathcal{D}_{\zeta} \subseteq \mathcal{D}_{\text{M0}}$. 
It follows that 
\begin{align*}
    \max\limits_{\mathbf{B} \in \mathcal{B}_{m_1}} \ \min\limits_{\mathbf{x} \in \mathcal{X}} \ \max\limits_{\mathbb{P}_{ \textup{r} }  \in \mathcal{D}_{\text{L}} } \ \mathbb{E}_{\mathbb{P}_{\textup{r}} } \left[ f\left(\mathbf{x},\mathbf{U}\boldsymbol{\Lambda}^{{\frac{1}{2}}} \mathbf{B}  \boldsymbol{\xi}_{ \textup{r} } +\boldsymbol{\mu}\right) \right] 
    \leq \min_{\mathbf{x} \in \mathcal{X}} \ \max_{\mathbb{P} \in \mathcal{D}_{\text{M0}}} \ \mathbb{E}_{\mathbb{P}} \left[ f\left(\mathbf{x},\boldsymbol{\xi}\right) \right],
\end{align*}
which demonstrates that the optimal value of Problem \eqref{Equ:DROM3} is a lower bound for that of Problem \eqref{Equ:DROM2} (i.e.,  Problem \eqref{Equ:DROM1}).


(ii) For any $m_1 < m_2 \leq m$, $\mathbf{B}_1 \in \mathbb{R}^{m\times m_1}$, and $\mathbf{C} \in \mathbb{R}^{m \times (m_2-m_1)}$ such that $\mathbf{B}_1^{\top} \mathbf{B}_1 = \mathbf{I}_{m_1}$ and $[\mathbf{B}_1, \mathbf{C}]^{\top} [\mathbf{B}_1, \mathbf{C}]  = \mathbf{I}_{m_2}$, 
we have $\mathbf{B}_2 = [\mathbf{B}_1, \mathbf{C}] \in \mathbb{R}^{m \times m_2}$. 
Meanwhile, we have $ \mathcal{B}_{m_2} = \{\mathbf{B} \in \mathbb{R}^{m\times m_2} \ |  \ \mathbf{B}^\top \mathbf{B} = \mathbf{I}_{m_2} \}$ 
and define 
$\boldsymbol{\zeta}_i = \mathbf{U} \boldsymbol{\Lambda}^{{\frac{1}{2}}} \mathbf{B}_i \boldsymbol{\xi}_{\textup{r}_i} + \boldsymbol{\mu} \in \mathbb{R}^{ m }$ for any $i \in [2]$, where $\boldsymbol{\xi}_{\textup{r}_i} \in \mathbb{R}^{m_i}$.
Clearly, $\mathbf{B}_2 \in \mathcal{B}_{m_2}$ because $\mathbf{B}_2^{\top} \mathbf{B}_2 = \mathbf{I}_{m_2}$.
We further define the ambiguity set of $\boldsymbol{\zeta}_i$ as
\begin{align}
    \mathcal{D}_{\zeta_i}=\left\{\mathbb{P}_{\zeta_i} \ \middle| \ \boldsymbol{\zeta}_i\sim \mathbb{P}_{\zeta_i}, \ \boldsymbol{\zeta}_i = \mathbf{U} \boldsymbol{\Lambda}^{{\frac{1}{2}}} \mathbf{B}_i \boldsymbol{\xi}_{\textup{r}_i} + \boldsymbol{\mu}, \ \boldsymbol{\xi}_{\textup{r}_i} \sim \mathbb{P}_{\textup{r}_i} \in\mathcal{D}_{\textup{r}_i} \right\}, \  \forall i \in [2], \label{eqn:ambiguity-D-zeta-i}
\end{align}
where $\mathcal{D}_{\textup{r}_i}$ represents the ambiguity set of $\boldsymbol{\xi}_{\textup{r}_i}$ for any $i \in [2]$.
Given $\boldsymbol{\zeta}_1 \sim \mathbb{P}_{\zeta_1} \in \mathcal{D}_{\zeta_1}$, there exists a $\boldsymbol{\xi}_{\textup{r}_1}\sim\mathbb{P}_{\textup{r}_1} \in \mathcal{D}_{\textup{r}_1}$ such that 
$\boldsymbol{\zeta}_1 = \mathbf{U} \boldsymbol{\Lambda}^{{\frac{1}{2}}} \mathbf{B}_1 \boldsymbol{\xi}_{\textup{r}_1} + \boldsymbol{\mu}
= \mathbf{U} \boldsymbol{\Lambda}^{{\frac{1}{2}}} \mathbf{B}_2 \bar{\boldsymbol{\xi}}_{\textup{r}_2} + \boldsymbol{\mu},$
where $\bar{\boldsymbol{\xi}}_{\textup{r}_2}=(\boldsymbol{\xi}^{\top}_{\textup{r}_1}, \boldsymbol{0}^{\top}_{m_2-m_1})^{\top} \in \mathbb{R}^{m_2}$.

By using $\mathcal{S}_{\textup{r}_i}$ (see definition in \eqref{eqn:support-sr-b}) to denote the support of $\boldsymbol{\xi}_{\textup{r}_i}$ for any $i \in [2]$, we have  
\[ 
\mathbb{P} \left\{\boldsymbol{\xi}_{\textup{r}_1}\in \mathcal{S}_{\textup{r}_1} \right\} 
= \mathbb{P} \left\{\mathbf{U}\boldsymbol{\Lambda}^{{\frac{1}{2}}} \mathbf{B}_1 \boldsymbol{\xi}_{\textup{r}_1}+\boldsymbol{\mu} \in \mathcal{S} \right\} 
= \mathbb{P} \left\{\mathbf{U} \boldsymbol{\Lambda}^{{\frac{1}{2}}} \mathbf{B}_2 \bar{\boldsymbol{\xi}}_{\textup{r}_2}+\boldsymbol{\mu}
\in \mathcal{S}  \right\} 
= 1, 
\]
where the second equality holds because $\mathbf{U} \boldsymbol{\Lambda}^{{\frac{1}{2}}} \mathbf{B}_1 \boldsymbol{\xi}_{\textup{r}_1} 
= \mathbf{U} \boldsymbol{\Lambda}^{{\frac{1}{2}}} \mathbf{B}_2 \bar{\boldsymbol{\xi}}_{\textup{r}_2}$.
It follows that $\mathbb{P} \{ \bar {\boldsymbol{\xi}}_{\textup{r}_2} \in \mathcal{S}_{\textup{r}_2} \} = 1$ by the definition of $\mathcal{S}_{\textup{r}_2}$.
In addition, we have $\mathbb{E} [ \bar {\boldsymbol{\xi}}_{\textup{r}_2}^{\top} ] \mathbb{E} [ \bar {\boldsymbol{\xi}}_{\textup{r}_2}] = 
\mathbb{E} [ \boldsymbol{\xi}_{\textup{r}_1}^{\top} ] \mathbb{E} [ \boldsymbol{\xi}_{\textup{r}_1}] \leq \gamma_1$ and
\[
\mathbb{E} \left[\bar {\boldsymbol{\xi}}_{\textup{r}_2} \bar{ \boldsymbol{\xi}}_{\textup{r}_2}^{\top}\right]= \left[\begin{array}{cc}\mathbb{E} \left[\boldsymbol{\xi}_{\textup{r}_1}\boldsymbol{\xi}_{\textup{r}_1}^{\top}\right] & \boldsymbol{0}_{m_1 \times (m_2-m_1)}\\ \boldsymbol{0}_{(m_2-m_1) \times m_1}&\boldsymbol{0}_{(m_2-m_1) \times (m_2-m_1)}\end{array}\right] \preceq \gamma_2\mathbf{I}_{m_2}.
\]
Thus, the probability distribution of $\bar {\boldsymbol{\xi}}_{\textup{r}_2}$ belongs to $\mathcal{D}_{\textup{r}_2}$. 
Meanwhile, by the definition of $ \mathcal{D}_{\zeta_i}$ for any $i \in [2]$ in \eqref{eqn:ambiguity-D-zeta-i} and $\boldsymbol{\zeta}_1 
= \mathbf{U} \boldsymbol{\Lambda}^{{\frac{1}{2}}} \mathbf{B}_2 \bar{\boldsymbol{\xi}}_{\textup{r}_2} + \boldsymbol{\mu} $, we have $ \mathbb{P}_{\zeta_1} \in \mathcal{D}_{\zeta_2} $ and further $ \mathcal{D}_{\zeta_1} \subseteq \mathcal{D}_{\zeta_2} $.
Therefore, for any $\mathbf{x} \in \mathcal{X}$, $\mathbf{B}_1 \in \mathbb{R}^{m\times m_1}$, and $\mathbf{C} \in \mathbb{R}^{m \times (m_2-m_1)}$ such that $\mathbf{B}_1^{\top} \mathbf{B}_1 = \mathbf{I}_{m_1}$ and $[\mathbf{B}_1, \mathbf{C}]^{\top} [\mathbf{B}_1, \mathbf{C}]  = \mathbf{I}_{m_2}$, we have
\begin{align} 
 \max\limits_{\mathbb{P}_{\zeta_1}\in \mathcal{D}_{\zeta_1}}\mathbb{E}_{\mathbb{P}_{\zeta_1}} \left[ f\left( \mathbf{x}, \boldsymbol{\zeta}_1 \right) \right] \leq \max\limits_{\mathbb{P}_{\zeta_2}\in \mathcal{D}_{\zeta_2}}\mathbb{E}_{\mathbb{P}_{\zeta_2}} \left[ f \left( \mathbf{x}, \boldsymbol{\zeta}_2 \right) \right]. \label{eqn:thm1-2-monoton-1}
\end{align}
Together with the definitions of $\boldsymbol{\zeta}_i$ ($\forall i \in [2]$) and $\mathbf{B}_2$, inequality \eqref{eqn:thm1-2-monoton-1} leads to
\begin{align*} 
&   \max\limits_{\mathbb{P}_{\textup{r}_1} \in \mathcal{D}_{\textup{r}_1}}\mathbb{E}_{\mathbb{P}_{\zeta_1}} \left[ f\left(\mathbf{x},\mathbf{U}\boldsymbol{\Lambda}^{{\frac{1}{2}}} \mathbf{B}_1  \boldsymbol{\xi}_{ \textup{r}_1 } +\boldsymbol{\mu}\right) \right]    \leq  \max\limits_{\mathbb{P}_{\textup{r}_2} \in \mathcal{D}_{\textup{r}_2}}\mathbb{E}_{\mathbb{P}_{\zeta_2}} \left[f \left(\mathbf{x},\mathbf{U}\boldsymbol{\Lambda}^{{\frac{1}{2}}} [\mathbf{B}_1, \mathbf{C}]  \boldsymbol{\xi}_{ \textup{r}_2 } +\boldsymbol{\mu}\right) \right].
\end{align*}
Considering an optimal solution $\mathbf{B}_1^* \in \mathbb{R}^{m\times m_1}$ of Problem \eqref{Equ:DROM3}, for any $\mathbf{x} \in \mathcal{X}$ and $\mathbf{C} \in \mathbb{R}^{m \times (m_2-m_1)}$ such that $[ \mathbf{B}_1^*, \mathbf{C}]^{\top} [\mathbf{B}_1^*, \mathbf{C}]  = \mathbf{I}_{m_2}$, we have
\begin{align*} 
&   \max\limits_{\mathbb{P}_{\textup{r}_1} \in \mathcal{D}_{\textup{r}_1}}\mathbb{E}_{\mathbb{P}_{\zeta_1}} \left[ f\left(\mathbf{x},\mathbf{U}\boldsymbol{\Lambda}^{{\frac{1}{2}}} \mathbf{B}_1^*  \boldsymbol{\xi}_{ \textup{r}_1 } +\boldsymbol{\mu}\right) \right]    \leq  \max\limits_{\mathbb{P}_{\textup{r}_2} \in \mathcal{D}_{\textup{r}_2}}\mathbb{E}_{\mathbb{P}_{\zeta_2}} \left[f \left(\mathbf{x},\mathbf{U}\boldsymbol{\Lambda}^{{\frac{1}{2}}} [\mathbf{B}_1^*, \mathbf{C}]  \boldsymbol{\xi}_{ \textup{r}_2 } +\boldsymbol{\mu}\right) \right].
\end{align*} 
For any $\mathbf{C} \in \mathbb{R}^{m \times (m_2-m_1)}$ such that $[ \mathbf{B}_1^*, \mathbf{C}]^{\top} [\mathbf{B}_1^*, \mathbf{C}]  = \mathbf{I}_{m_2}$, we have
\begin{align} 
& \min\limits_{\mathbf{x} \in \mathcal{X}} \  \max\limits_{\mathbb{P}_{\textup{r}_1} \in \mathcal{D}_{\textup{r}_1}}\mathbb{E}_{\mathbb{P}_{\zeta_1}} \left[ f\left(\mathbf{x},\mathbf{U}\boldsymbol{\Lambda}^{{\frac{1}{2}}} \mathbf{B}_1^*  \boldsymbol{\xi}_{ \textup{r}_1 } +\boldsymbol{\mu}\right) \right]    \leq  \min\limits_{\mathbf{x} \in \mathcal{X}} \ \max\limits_{\mathbb{P}_{\textup{r}_2} \in \mathcal{D}_{\textup{r}_2}}\mathbb{E}_{\mathbb{P}_{\zeta_2}} \left[f \left(\mathbf{x},\mathbf{U}\boldsymbol{\Lambda}^{{\frac{1}{2}}} [\mathbf{B}_1^*, \mathbf{C}]  \boldsymbol{\xi}_{ \textup{r}_2 } +\boldsymbol{\mu}\right) \right]. \label{Cons:The1-modified-1}
\end{align} 
It follows that
\begin{align*} 
& \max\limits_{\mathbf{B}_1^\top \mathbf{B}_1  = \mathbf{I}_{m_1}} \ \min\limits_{\mathbf{x} \in \mathcal{X}} \ \max\limits_{\mathbb{P}_{\textup{r}_1} \in \mathcal{D}_{\textup{r}_1}}\mathbb{E}_{\mathbb{P}_{\zeta_1}} \left[ f\left(\mathbf{x},\mathbf{U}\boldsymbol{\Lambda}^{{\frac{1}{2}}} \mathbf{B}_1  \boldsymbol{\xi}_{ \textup{r}_1 } +\boldsymbol{\mu}\right) \right]  = \min\limits_{\mathbf{x} \in \mathcal{X}} \  \max\limits_{\mathbb{P}_{\textup{r}_1} \in \mathcal{D}_{\textup{r}_1}}\mathbb{E}_{\mathbb{P}_{\zeta_1}} \left[ f\left(\mathbf{x},\mathbf{U}\boldsymbol{\Lambda}^{{\frac{1}{2}}} \mathbf{B}_1^*  \boldsymbol{\xi}_{ \textup{r}_1 } +\boldsymbol{\mu}\right) \right]   \\ 
& \leq \min\limits_{\mathbf{x} \in \mathcal{X}} \ \max\limits_{\mathbb{P}_{\textup{r}_2} \in \mathcal{D}_{\textup{r}_2}}\mathbb{E}_{\mathbb{P}_{\zeta_2}} \left[f \left(\mathbf{x},\mathbf{U}\boldsymbol{\Lambda}^{{\frac{1}{2}}} [\mathbf{B}_1^*, \mathbf{C}]  \boldsymbol{\xi}_{ \textup{r}_2 } +\boldsymbol{\mu}\right) \right] \leq \max\limits_{\mathbf{B}_2 \in \mathcal{B}_{m_2} } \ \min\limits_{\mathbf{x} \in \mathcal{X}} \ \max\limits_{\mathbb{P}_{\textup{r}_2} \in \mathcal{D}_{\textup{r}_2}}\mathbb{E}_{\mathbb{P}_{\zeta_2}} \left[ f\left(\mathbf{x},\mathbf{U}\boldsymbol{\Lambda}^{{\frac{1}{2}}} \mathbf{B}_2  \boldsymbol{\xi}_{ \textup{r}_2 } +\boldsymbol{\mu}\right) \right],
\end{align*}
where the first inequality holds by \eqref{Cons:The1-modified-1} and the second inequality holds because $[\mathbf{B}_1^*, \mathbf{C}] \in \mathcal{B}_{m_2}$.
That is, the optimal value of Problem \eqref{Equ:DROM3} is nondecreasing in $m_1$.

(iii) When $m_1 = m$, we have $\mathbf{B} \in \mathcal{B}_m \subseteq \mathbb{R}^{m \times m}$, i.e., $\mathbf{B}^{\top} \mathbf{B} = \mathbf{I}_m$. 
First, we have $\Theta_{\textup{L}}(m) \leq \Theta_{\textup{M}}(m)$ by the conclusion (i).
Second, when $\mathbf{B} = \mathbf{I}_m$, Problem \eqref{Equ:DROM3} becomes Problem \eqref{Equ:DROM2}. 
Because $\mathbf{B} = \mathbf{I}_m$ is a feasible solution of Problem \eqref{Equ:DROM3}, it follows that $\Theta_{\textup{L}}(m) \geq \Theta_{\textup{M}}(m)$.
Therefore, we have $\Theta_{\textup{L}}(m) = \Theta_{\textup{M}}(m)$.
\Halmos

\subsection{Proof of Theorem \ref{Theo:sameoptimal}}
We construct a solution $(\mathbf{x}^{\dag}, s^{\dag}, \hat{\boldsymbol{\lambda}}^{\dag}, \mathbf{q}_{\textup{r}}^{\dag}, \mathbf{Q}_{\textup{r}}^{\dag}, \mathbf{B}^{\dag})$ of Problems \eqref{Equ:appro} and \eqref{Equ:appro-1} by setting 
$\mathbf{x}^{\dag} = \mathbf{x}^*$, 
$s^{\dag} = s^*$, 
$\hat{\boldsymbol{\lambda}}^{\dag} = \hat{\boldsymbol{\lambda}}^*$,
$\mathbf{q}_{\textup{r}}^{\dag} = (\boldsymbol{\delta}^{\top}, \mathbf{0}_{m_1-K}^{\top})^{\top}$,
$\mathbf{Q}_{\textup{r}}^{\dag} =  {\scriptsize \begin{bmatrix}
     \mathbf{Y}_{11} & \mathbf{0}_{K \times (m_1 - K)} \\ 
     \mathbf{0}_{(m_1 - K) \times K} & \mathbf{0}_{(m_1 - K) \times (m_1 - K)} 
 \end{bmatrix} }$, 
 and $\mathbf{B}^{\dag} = [\mathbf{V}, \mathbf{C}]$, 
 where $\mathbf{C} \in \mathbb{R}^{m \times (m_1-K)}$ and $[\mathbf{V}, \mathbf{C}]^{\top} [\mathbf{V}, \mathbf{C}] = \mathbf{I}_{m_1}$. 
First, we show this constructed solution is feasible to Problems \eqref{Equ:appro} and \eqref{Equ:appro-1}.
Clearly, this solution satisfies constraints \eqref{Cons:appro2}. 
In addition, from Problem \eqref{Equ:MainProblem}, as $ \mathbf{q}^{\prime} = \mathbf{V} \boldsymbol{\delta} $ and $\mathbf{Q}^{\prime} = \mathbf{V} \mathbf{Y}_{11} \mathbf{V}^{\top}$, for any $k \in [K]$,  we have
\begin{align*}
\begin{bmatrix}  S_k & \hspace{0.1 in}  \frac{1}{2}  \left(\mathbf{V}  \boldsymbol{\delta} +\left(\mathbf{U}\boldsymbol{\Lambda}^{{\frac{1}{2}}}\right)^{\top}\left(\mathbf{A}^{\top}\boldsymbol{\lambda}^*_k-y_k(\mathbf{x}^*)\right)\right)^{\top}\\ 
 \frac{1}{2} \left( \mathbf{V}  \boldsymbol{\delta} +\left(\mathbf{U}\boldsymbol{\Lambda}^{{\frac{1}{2}}}\right)^{\top}\left(\mathbf{A}^{\top}\boldsymbol{\lambda}^*_k-y_k(\mathbf{x}^*)\right) \right) & \mathbf{V} \mathbf{Y}_{11} \mathbf{V}^{\top} \end{bmatrix} \succeq 0, 
\end{align*}
which, by Schur complement, is equivalent to
\begin{small}
\begin{equation}
S_k \left( \mathbf{V} \mathbf{Y}_{11} \mathbf{V}^{\top} \right) \succeq  \frac{1}{4}  \left(\mathbf{V}  \boldsymbol{\delta} +\left(\mathbf{U}\boldsymbol{\Lambda}^{{\frac{1}{2}}}\right)^{\top}\left(\mathbf{A}^{\top}\boldsymbol{\lambda}^*_k-y_k(\mathbf{x}^*)\right)\right) 
\left( \mathbf{V}  \boldsymbol{\delta} +\left(\mathbf{U}\boldsymbol{\Lambda}^{{\frac{1}{2}}}\right)^{\top}\left(\mathbf{A}^{\top}\boldsymbol{\lambda}^*_k-y_k(\mathbf{x}^*)\right) \right)^{\top}. \label{eqn:theorem-2-1}
\end{equation}
\end{small}%
From \eqref{eqn:theorem-2-1}, for any $k \in [K]$, we have the following inequality holds by Lemma \ref{Lem:VXV}:
\begin{align*}
&S_k\left( [\mathbf{V}, \mathbf{C}]^{\top} \mathbf{V} \mathbf{Y}_{11} \mathbf{V}^{\top} [\mathbf{V}, \mathbf{C}] \right)  \\ 
&\hspace{1cm} \succeq \frac{1}{4}  [\mathbf{V}, \mathbf{C}]^{\top}\left(\mathbf{V}  \boldsymbol{\delta} +\left(\mathbf{U}\boldsymbol{\Lambda}^{{\frac{1}{2}}}\right)^{\top}\left(\mathbf{A}^{\top}\boldsymbol{\lambda}^*_k-y_k(\mathbf{x}^*)\right)\right) 
\left( \mathbf{V}  \boldsymbol{\delta} +\left(\mathbf{U}\boldsymbol{\Lambda}^{{\frac{1}{2}}}\right)^{\top}\left(\mathbf{A}^{\top}\boldsymbol{\lambda}^*_k-y_k(\mathbf{x}^*)\right) \right)^{\top} [\mathbf{V}, \mathbf{C}], \nonumber
\end{align*}
which is equivalent to 
\begin{align}
&S_k \mathbf{Q}_{\textup{r}}^{\dag} \succeq  \frac{1}{4}  \left(  \mathbf{q}_{\textup{r}}^{\dag} +\left(\mathbf{U}\boldsymbol{\Lambda}^{{\frac{1}{2}}}\mathbf{B}^{\dag} \right)^{\top}\left(\mathbf{A}^{\top}\boldsymbol{\lambda}^*_k-y_k(\mathbf{x}^*)\right)\right) 
\left( \mathbf{q}_{\textup{r}}^{\dag} +\left(\mathbf{U}\boldsymbol{\Lambda}^{{\frac{1}{2}}}\mathbf{B}^{\dag}\right)^{\top}\left(\mathbf{A}^{\top}\boldsymbol{\lambda}^*_k-y_k(\mathbf{x}^*)\right) \right)^{\top}   \label{eqn:theorem-2-2}
\end{align}
by the construction of the solution $\mathbf{q}_{\textup{r}}^{\dag}, \mathbf{Q}_{\textup{r}}^{\dag}, \mathbf{B}^{\dag}$ and $[\mathbf{V}, \mathbf{C}]^{\top} \mathbf{V} = [\mathbf{I}_{K}, \mathbf{0}_{K \times (m_1-K)}]^{\top}$. By Schur complement, \eqref{eqn:theorem-2-2}  indicates that the constructed solution 
$(\mathbf{x}^{\dag}, s^{\dag}, \hat{\boldsymbol{\lambda}}^{\dag}, \mathbf{q}_{\textup{r}}^{\dag}, \mathbf{Q}_{\textup{r}}^{\dag}, \mathbf{B}^{\dag})$ also satisfies constraints \eqref{Cons:appro1} and thus it is a feasible solution of Problems \eqref{Equ:appro} and \eqref{Equ:appro-1}.

Second, we show the objective value of this feasible solution $(\mathbf{x}^{\dag}, s^{\dag}, \hat{\boldsymbol{\lambda}}^{\dag}, \mathbf{q}_{\textup{r}}^{\dag}, \mathbf{Q}_{\textup{r}}^{\dag}, \mathbf{B}^{\dag})$ is equal to the optimal value of Problem \eqref{Equ:MainProblem}.
The objective value corresponding to this solution is
\begin{align}
s^{\dag} + \gamma_2 \mathbf{I}_{m_1} \bullet \mathbf{Q}_{\textup{r}}^{\dag}  + \sqrt{\gamma_1}\left \| \mathbf{q}_{\textup{r}}^{\dag} \right \|_2 = \ & s^* +\gamma_2\mathbf{I}_{K} \bullet \mathbf{Y}_{11} + \sqrt{\gamma_1}\left \| \boldsymbol{\delta} \right \|_2 \nonumber\\
=\ &s^* +\gamma_2\mathbf{I}_{K} \bullet (\mathbf{Y}_{11}\mathbf{V}^{\top} \mathbf{V}) + \sqrt{\gamma_1}\left \| \boldsymbol{\delta} \right \|_2 \nonumber\\
=\ &s^* +\gamma_2\mathbf{I}_{m} \bullet (\mathbf{V}\mathbf{Y}_{11}\mathbf{V}^{\top}) + \sqrt{\gamma_1}\left \| \boldsymbol{\delta} \right \|_2 \nonumber\\ 
=\ & s^* +\gamma_2\mathbf{I}_{m} \bullet \mathbf{Q}^{\prime} + \sqrt{\gamma_1}\left \| \boldsymbol{\delta} \right \|_2 \nonumber \\
=\ & s^* +\gamma_2\mathbf{I}_{m} \bullet \mathbf{Q}^{\prime} + \sqrt{\gamma_1} \sqrt{\boldsymbol{\delta}^{\top} \boldsymbol{\delta}} \nonumber \\
=\ & s^* +\gamma_2\mathbf{I}_{m} \bullet \mathbf{Q}^{\prime} + \sqrt{\gamma_1} \sqrt{\boldsymbol{\delta}^{\top} \mathbf{V}^{\top} \mathbf{V} \boldsymbol{\delta}} \nonumber \\
=\ & s^* +\gamma_2\mathbf{I}_{m} \bullet \mathbf{Q}^{\prime} + \sqrt{\gamma_1} \sqrt{\mathbf{q}^{\prime \top} \mathbf{q}^{\prime}} \nonumber \\
=\ & s^* +\gamma_2\mathbf{I}_{m} \bullet \mathbf{Q}^{\prime} + \sqrt{\gamma_1} \left \| \mathbf{q}^{\prime} \right \|_2 \nonumber\\ 
=\ &  \Theta_{\textup{M}}(m),  \label{eqn:theorem-2-3}
\end{align}
where the first equality holds by the construction of $(\mathbf{x}^{\dag}, s^{\dag}, \hat{\boldsymbol{\lambda}}^{\dag}, \mathbf{q}_{\textup{r}}^{\dag}, \mathbf{Q}_{\textup{r}}^{\dag}, \mathbf{B}^{\dag})$,
the second equality holds because $\mathbf{V}^{\top} \mathbf{V} = \mathbf{I}_{K}$,
the third equality holds by the cyclic property of a matrix's trace,
the fourth equality holds by the definition of $\mathbf{Q}^{\prime}$ in Theorem \ref{Prop:lowrank}, 
and the seventh equality holds because $ \mathbf{q}^{\prime} = \mathbf{V} \boldsymbol{\delta} $. 
\Halmos

\subsection{Proof of Proposition \ref{prop:bilinear-sdp} } 
We consider the Lagrangian dual of the inner minimization part (i.e., Problem \eqref{Equ:appro-1}) of Problem \eqref{Equ:appro} as follows:
\begin{align}
    \max\limits_{\substack{ \begin{bmatrix}  t_k & \mathbf{p}_k^{\top} \\ \mathbf{p}_k & \mathbf{P}_k \end{bmatrix}  \succeq 0, \\ \forall k \in [K],\\ \mathbf{Z} \succeq 0}} \min\limits_{ \substack{\mathbf{x}, s, {\hat{\boldsymbol{\lambda}}} \geq 0, \\ \mathbf{q}_{ \textup{r} },\mathbf{Q}_{ \textup{r} } \succeq 0 } }   \ \mathcal{L} \left( \mathbf{x}, s, {\hat{\boldsymbol{\lambda}}},  \mathbf{q}_{ \textup{r} }, \mathbf{Q}_{ \textup{r} };  \mathbf{Z},  t_k, \mathbf{p}_k,  \mathbf{P}_k, \forall k \in [K] \right), \label{eqn:lagrangian-dual}
\end{align}
where the Lagrangian function 
{\small \begin{align*}
& \mathcal{L} \left( \mathbf{x}, s, {\hat{\boldsymbol{\lambda}}},  \mathbf{q}_{ \textup{r} }, \mathbf{Q}_{ \textup{r} };  \mathbf{Z},  t_k, \mathbf{p}_k,  \mathbf{P}_k, \forall k \in [K] \right) \\
=  & s + \gamma_2 \mathbf{I}_{m_1} \bullet \mathbf{Q}_{ \textup{r} }  + \sqrt{\gamma_1}\left \| \mathbf{q}_{ \textup{r} } \right \|_2 - \mathbf{Z} \bullet \left( \sum_{i=1}^n (\boldsymbol{\Delta}_i x_i) + \boldsymbol{\Delta}_0 \right) - \sum_{k=1}^K  \begin{bmatrix}  t_k & \mathbf{p}_k^{\top} \\ \mathbf{p}_k & \mathbf{P}_k \end{bmatrix}  \bullet \nonumber\\  
    & \begin{bmatrix}  s-y_k^0(\mathbf{x})-\boldsymbol{\lambda}_k^{\top}b-y_k(\mathbf{x})^{\top}\boldsymbol{\mu} +\boldsymbol{\lambda}_k^{\top}\mathbf{A}\boldsymbol{\mu} & \hspace{0.1 in} \frac{1}{2}  \left(\mathbf{q}_{ \textup{r} } +\left(\mathbf{U}\boldsymbol{\Lambda}^{{\frac{1}{2}}} \mathbf{B} \right)^{\top}\left(\mathbf{A}^{\top}\boldsymbol{\lambda}_k-y_k(\mathbf{x})\right)\right)^{\top} \\ 
\frac{1}{2} \left( \mathbf{q}_{ \textup{r} } +\left(\mathbf{U} \boldsymbol{\Lambda}^{{\frac{1}{2}}}\mathbf{B} \right)^{\top}\left(\mathbf{A}^{\top}\boldsymbol{\lambda}_k-y_k(\mathbf{x})\right) \right) & \mathbf{Q}_{ \textup{r} } \end{bmatrix} \\
= &   \left(1-\sum_{k=1}^K t_k\right) s -  \sum_{k=1}^K \left( t_k \left(\mathbf{A}\boldsymbol{\mu}-\mathbf{b}\right)^{\top} + \mathbf{p}_k^{\top}\left(\mathbf{U} \boldsymbol{\Lambda}^{{\frac{1}{2}}}\mathbf{B} \right)^{\top} \mathbf{A}^{\top}  \right)\boldsymbol{\lambda}_k  + \sqrt{\gamma_1}\left\| \mathbf{q}_{ \textup{r} } \right\|_2 - \sum_{k=1}^K \mathbf{p}_k^{\top} \mathbf{q}_{\textup{r}}
 \\ &   +\left(\gamma_2 \mathbf{I}_{m_1} - \sum_{k=1}^K \mathbf{P}_k\right)  \bullet\mathbf{Q}_{\textup{r}} - \sum_{i=1}^\tau \sum_{j=1}^\tau z_{ij}\left(\mathbf{a}_{ij}\mathbf{x}+a_{ij}^0\right) +\sum_{k=1}^K \left(  t_k y_k^0(x) + \left(t_k \boldsymbol{\mu}^{\top}  +  \mathbf{p}_k^{\top} \left(\mathbf{U} \boldsymbol{\Lambda}^{{\frac{1}{2}}}\mathbf{B} \right)^{\top}\right) y_k(x)\right)  \\
=&    \left(1-\sum_{k=1}^K t_k\right) s -  \sum_{k=1}^K \left( t_k \left(\mathbf{A}\boldsymbol{\mu}-\mathbf{b}\right)^{\top} + \mathbf{p}_k^{\top}\left(\mathbf{U} \boldsymbol{\Lambda}^{{\frac{1}{2}}}\mathbf{B} \right)^{\top} \mathbf{A}^{\top}  \right)\boldsymbol{\lambda}_k  + \sqrt{\gamma_1}\left\| \mathbf{q}_{ \textup{r} } \right\|_2 - \sum_{k=1}^K \mathbf{p}_k^{\top} \mathbf{q}_{\textup{r}}  \\
&  +\left(\gamma_2 \mathbf{I}_{m_1} - \sum_{k=1}^K \mathbf{P}_k\right)  \bullet\mathbf{Q}_{\textup{r}}   - \sum_{i=1}^\tau \sum_{j=1}^\tau z_{ij}\mathbf{a}_{ij}\mathbf{x} - \sum_{i=1}^\tau \sum_{j=1}^\tau z_{ij} a_{ij}^0 + \sum_{k=1}^K \left(  t_k \mathbf{w}_k^0 + \left(t_k \boldsymbol{\mu}^{\top}  +  \mathbf{p}_k^{\top} \left(\mathbf{U} \boldsymbol{\Lambda}^{{\frac{1}{2}}}\mathbf{B} \right)^{\top}\right)  \mathbf{W}_k \right) \mathbf{x}\\
&  + \sum_{k=1}^K \left(t_k d_k^0 + \left(t_k \boldsymbol{\mu}^{\top}  +  \mathbf{p}_k^{\top} \left(\mathbf{U} \boldsymbol{\Lambda}^{{\frac{1}{2}}}\mathbf{B} \right)^{\top}\right) \mathbf{d}_k \right) \\
=&    \left(1-\sum_{k=1}^K t_k\right) s -  \sum_{k=1}^K \left( t_k \left(\mathbf{A}\boldsymbol{\mu}-\mathbf{b}\right)^{\top} + \mathbf{p}_k^{\top}\left(\mathbf{U} \boldsymbol{\Lambda}^{{\frac{1}{2}}}\mathbf{B} \right)^{\top} \mathbf{A}^{\top}  \right)\boldsymbol{\lambda}_k  + \sqrt{\gamma_1}\left\| \mathbf{q}_{ \textup{r} } \right\|_2 - \sum_{k=1}^K \mathbf{p}_k^{\top} \mathbf{q}_{\textup{r}}  \\
&  +\left(\gamma_2 \mathbf{I}_{m_1} - \sum_{k=1}^K \mathbf{P}_k\right)  \bullet\mathbf{Q}_{\textup{r}} + \left( \sum_{k=1}^K \left(  t_k \mathbf{w}_k^0 + \left(t_k \boldsymbol{\mu}^{\top}  +  \mathbf{p}_k^{\top} \left(\mathbf{U} \boldsymbol{\Lambda}^{{\frac{1}{2}}}\mathbf{B} \right)^{\top}\right)  \mathbf{W}_k \right)  - \sum_{i=1}^\tau \sum_{j=1}^\tau z_{ij}\mathbf{a}_{ij} \right) \mathbf{x}   \\
&+ \sum_{k=1}^K \left(t_k d_k^0 + \left(t_k \boldsymbol{\mu}^{\top}  +  \mathbf{p}_k^{\top} \left(\mathbf{U} \boldsymbol{\Lambda}^{{\frac{1}{2}}}\mathbf{B} \right)^{\top}\right) \mathbf{d}_k \right) - \sum_{i=1}^\tau \sum_{j=1}^\tau z_{ij} a_{ij}^0.
\end{align*}}%
To present the objective value of the inner minimization problem of \eqref{eqn:lagrangian-dual} from going to negative infinity, we require
\begin{subeqnarray} 
&& 1-\sum_{k=1}^K t_k = 0,\ \sqrt{\gamma_1} - \left\|\sum_{k=1}^K \mathbf{p}_k \right\|_2 \geq 0,\ \gamma_2 \mathbf{I}_{m_1} - \sum_{k=1}^K \mathbf{P}_k \succeq 0, \slabel{eqn:lagrangian-cons-1} \\
&& t_k \left(\mathbf{A}\boldsymbol{\mu}-\mathbf{b}\right)^{\top} + \mathbf{p}_k^{\top}\left(\mathbf{U} \boldsymbol{\Lambda}^{{\frac{1}{2}}}\mathbf{B} \right)^{\top} \mathbf{A}^{\top}   \leq 0, \ \forall k \in [K], \slabel{eqn:lagrangian-cons-2} \\
&& \sum_{k=1}^K \left(  t_k \mathbf{w}_k^0 + \left(t_k \boldsymbol{\mu}^{\top}  +  \mathbf{p}_k^{\top} \left(\mathbf{U} \boldsymbol{\Lambda}^{{\frac{1}{2}}}\mathbf{B} \right)^{\top} \right)  \mathbf{W}_k \right) - \sum_{i=1}^\tau \sum_{j=1}^\tau z_{ij}\mathbf{a}_{ij} = 0, \slabel{eqn:lagrangian-cons-3} \\
&& \begin{bmatrix}  t_k & \mathbf{p}_k^{\top} \\ \mathbf{p}_k & \mathbf{P}_k \end{bmatrix} \succeq 0, \ \forall k\in [K], \ \mathbf{Z} \succeq 0. \slabel{eqn:lagrangian-cons-4}
\end{subeqnarray}
Then, the dual problem of Problem \eqref{Equ:appro-1} can be described as follows:
\begin{align}
    \max\limits_{\substack{t_k, \mathbf{p}_k, \mathbf{P}_k, \forall k \in [K], \mathbf{Z}  }}  \  & \sum_{k=1}^K \left( t_k d_k^0 + \left( t_k \boldsymbol{\mu}^{\top}  +  \mathbf{p}_k^{\top} \left( \mathbf{U} \boldsymbol{\Lambda}^{{\frac{1}{2}}}\mathbf{B} \right)^{\top} \right) \mathbf{d}_k \right) - \sum_{i=1}^\tau \sum_{j=1}^\tau z_{ij} a_{ij}^0 \label{eqn:lagrangian-obj} \\
\rm{s.t.} \hspace{0.5cm} & \eqref{eqn:lagrangian-cons-1} - \eqref{eqn:lagrangian-cons-4}. \nonumber 
\end{align}

By integrating the outer maximization part of Problem \eqref{Equ:appro} and Problem \eqref{eqn:lagrangian-obj}, 
we obtain the bilinear SDP problem \eqref{Equ:bilinear}. 
Now we would like to prove the strong duality between Problem \eqref{Equ:appro-1} and Problem \eqref{eqn:lagrangian-obj}; 
that is, these two problems share the same optimal value, which further shows that Problem \eqref{Equ:appro} has the same optimal value as Problem \eqref{Equ:bilinear}.
To that end, we find an interior point of Problem \eqref{Equ:appro-1}.

Let $\mathbf{x}^{\prime}$ be an interior point in $\mathcal{X}$, 
we can construct an interior point by setting $\boldsymbol{\hat{\lambda}}^{\prime} = \{\boldsymbol{1}_l,\ldots, \boldsymbol{1}_l\}$, 
$s^{\prime}= \sum_{k=1}^K \left| y_k^0(\mathbf{x}^{\prime}) + \boldsymbol{1}_l^{\top} \mathbf{b} + y_k(\mathbf{x}^{\prime})^{\top} \boldsymbol{\mu} - \boldsymbol{1}_l^{\top}\mathbf{A} \boldsymbol{\mu}\right| + 1$, 
$\mathbf{q}^{\prime}_{\textup{r}} = 0$, 
and $\mathbf{Q}^{\prime}_{\textup{r}} = \sum_{k=1}^K 1/(4(s'-y_k^0(\mathbf{x}^{\prime})-\boldsymbol{1}_l^{\top}b-y_k(\mathbf{x}^{\prime})^{\top}\boldsymbol{\mu} +\boldsymbol{1}_l^{\top}\mathbf{A}\boldsymbol{\mu}))  (\mathbf{U} \boldsymbol{\Lambda}^{{\frac{1}{2}}}\mathbf{B})^{\top} (\mathbf{A}^{\top}\boldsymbol{1}_l-y_k(\mathbf{x}^{\prime})) (\mathbf{A}^{\top}\boldsymbol{1}_l-y_k(\mathbf{x}^{\prime}))^{\top} (\mathbf{U} \boldsymbol{\Lambda}^{{\frac{1}{2}}}\mathbf{B}) + \mathbf{I}_{m_1}$. 
Clearly, $(\mathbf{U} \boldsymbol{\Lambda}^{{\frac{1}{2}}}\mathbf{B})^{\top} (\mathbf{A}^{\top}\boldsymbol{1}_l-y_k(\mathbf{x}^{\prime})) (\mathbf{A}^{\top}\boldsymbol{1}_l-y_k(\mathbf{x}^{\prime}))^{\top} (\mathbf{U} \boldsymbol{\Lambda}^{{\frac{1}{2}}}\mathbf{B}) \succeq 0$. 
Thus, $\mathbf{Q}^{\prime}_{\textup{r}} \succ 0$. 
Now we only need to show that constraints \eqref{Cons:appro1} hold in the positive-definite sense with respect to this constructed solution.

By the construction of $\mathbf{Q}^{\prime}_{ \textup{r} }$, for any $k\in [K]$, we have 
\begin{align}
&\mathbf{Q}^{\prime}_{ \textup{r} } - \frac{\left( \left(\mathbf{U} \boldsymbol{\Lambda}^{{\frac{1}{2}}}\mathbf{B} \right)^{\top}\left(\mathbf{A}^{\top}\boldsymbol{1}_l-y_k(\mathbf{x}^{\prime})\right) \right) \left(\left(\mathbf{U}\boldsymbol{\Lambda}^{{\frac{1}{2}}} \mathbf{B} \right)^{\top}\left(\mathbf{A}^{\top}\boldsymbol{1}_l-y_k(\mathbf{x}^{\prime})\right)\right)^{\top}}{4\left(s'-y_k^0(\mathbf{x}^{\prime})-\boldsymbol{1}_l^{\top}b-y_k(\mathbf{x}^{\prime})^{\top}\boldsymbol{\mu} +\boldsymbol{1}_l^{\top}\mathbf{A}\boldsymbol{\mu}\right)}\nonumber \\
& = \sum_{\forall k^{\prime} \in [K]: k^{\prime} \neq k} \frac{\left( \left(\mathbf{U} \boldsymbol{\Lambda}^{{\frac{1}{2}}}\mathbf{B} \right)^{\top}\left(\mathbf{A}^{\top}\boldsymbol{1}_l-y_{k^{\prime}}(\mathbf{x}^{\prime})\right) \right) \left(\left(\mathbf{U}\boldsymbol{\Lambda}^{{\frac{1}{2}}} \mathbf{B} \right)^{\top}\left(\mathbf{A}^{\top}\boldsymbol{1}_l-y_{k^{\prime}}(\mathbf{x}^{\prime})\right)\right)^{\top}}{4\left(s'-y_{k^{\prime}}^0(\mathbf{x}^{\prime})-\boldsymbol{1}_l^{\top}b-y_{k^{\prime}}(\mathbf{x}^{\prime})^{\top}\boldsymbol{\mu} +\boldsymbol{1}_l^{\top}\mathbf{A}\boldsymbol{\mu}\right)} + \mathbf{I}_{m_1}  \succ 0, \label{propCons:dual}
\end{align}
where $s^{\prime}-y_{k^{\prime}}^0(\mathbf{x}^{\prime})-\boldsymbol{1}_l^{\top}b-y_{k^{\prime}}(\mathbf{x}^{\prime})^{\top}\boldsymbol{\mu} +\boldsymbol{1}_l^{\top}\mathbf{A}\boldsymbol{\mu} > 0$ by the construction of $s^{\prime}$.
By Schur complement, \eqref{propCons:dual} is equivalent to 
\begin{align*}
\begin{bmatrix}  s'-y_k^0(\mathbf{x}^{\prime})-\boldsymbol{1}_l^{\top}b-y_k(\mathbf{x}^{\prime})^{\top}\boldsymbol{\mu} +\boldsymbol{1}_l^{\top}\mathbf{A}\boldsymbol{\mu} & \hspace{0.1 in} \frac{1}{2}  \left(\left(\mathbf{U}\boldsymbol{\Lambda}^{{\frac{1}{2}}} \mathbf{B} \right)^{\top}\left(\mathbf{A}^{\top}\boldsymbol{1}_l-y_k(\mathbf{x}^{\prime})\right)\right)^{\top} \\ 
	 \frac{1}{2} \left( \left(\mathbf{U} \boldsymbol{\Lambda}^{{\frac{1}{2}}}\mathbf{B} \right)^{\top}\left(\mathbf{A}^{\top}\boldsymbol{1}_l-y_k(\mathbf{x}^{\prime})\right) \right) & \mathbf{Q}^{\prime}_{ \textup{r} } \end{bmatrix} \succ 0, \ \forall k \in [K].
\end{align*}
Thus, $(\mathbf{x}^{\prime}, s', \boldsymbol{\hat{\lambda}}^{\prime}, \mathbf{q}^{\prime}_{\textup{r}}, \mathbf{Q}^{\prime}_{\textup{r}})$ is an interior point of Problem \eqref{Equ:appro-1} and the strong duality between Problem \eqref{Equ:appro-1} and Problem \eqref{eqn:lagrangian-obj} holds.
\Halmos

\section{Supplement to Section \ref{sec:upper-bound}}

\subsection{Proof of Theorem \ref{Theo:upperbound}}
(i) For any $\boldsymbol{\xi}_{\textup{I}} \sim \mathbb{P}_{\text{I}} \in \mathcal{D}_{\textup{M}}$, we have 
$\mathbb{E}_{\mathbb{P}_{\text{I}}} [\boldsymbol{\xi}_{\textup{I}} \boldsymbol{\xi}_{\textup{I}}^{\top}] \preceq \gamma_2\mathbf{I}_{m_1}$.
Then, by Lemma \ref{Lem:VXV}, for any given $\mathbf{x} \in \mathcal{X}$ and $\mathbf{B} \in \mathcal{B}_{m_1}$, i.e., $\mathbf{B}^{\top} \mathbf{B} = \mathbf{I}_{m_1}$, we further have 
$ \mathbf{B}^{\top} ( \mathbb{E}_{\mathbb{P}_{\text{I}}} [\boldsymbol{\xi}_{\textup{I}} \boldsymbol{\xi}_{\textup{I}}^{\top}] ) \mathbf{B} \preceq \mathbf{B}^{\top}  ( \gamma_2\mathbf{I}_{m_1} ) \mathbf{B}$, i.e., $\mathbb{E}_{\mathbb{P}_{\text{I}}} [ \mathbf{B}^{\top}\boldsymbol{\xi}_{\textup{I}} \boldsymbol{\xi}_{\textup{I}}^{\top}\mathbf{B}] \preceq \gamma_2\mathbf{B}^{\top} \mathbf{I}_{m_1}\mathbf{B} = \gamma_2\mathbf{I}_{m_1}$.
It follows that $\mathcal{D}_{\textup{M}} \subseteq \mathcal{D}_{\textup{U}}$. 
Thus, given any $\mathbf{x} \in \mathcal{X}$ and $\mathbf{B} \in \mathcal{B}_{m_1}$, we have
\begin{align*}
    \max\limits_{\mathbb{P}_{\text{I}}  \in \mathcal{D}_{\textup{U}} } \ \mathbb{E}_{\mathbb{P}_{\text{I}} } \left[ f\left(\mathbf{x},\mathbf{U}\boldsymbol{\Lambda}^{{\frac{1}{2}}} \boldsymbol{\xi}_{ \textup{I} } +\boldsymbol{\mu}\right) \right]  
    \geq \max_{\mathbb{P}_{\text{I}} \in \mathcal{D}_{\textup{M}}} \ \mathbb{E}_{\mathbb{P}_{\text{I}}} \left[ f\left(\mathbf{x},\mathbf{U}\boldsymbol{\Lambda}^{{\frac{1}{2}}} \boldsymbol{\xi}_{ \textup{I} } +\boldsymbol{\mu}\right) \right].
\end{align*} 
It follows that 
\begin{align*}
    \min\limits_{\mathbf{B} \in \mathcal{B}_{m_1}} \ \min\limits_{\mathbf{x} \in \mathcal{X}} \ \max\limits_{\mathbb{P}_{\text{I}}  \in \mathcal{D}_{\textup{U}} } \ \mathbb{E}_{\mathbb{P}_{\text{I}} } \left[ f\left(\mathbf{x},\mathbf{U}\boldsymbol{\Lambda}^{{\frac{1}{2}}} \boldsymbol{\xi}_{ \textup{I} } +\boldsymbol{\mu}\right) \right] 
    \geq \min_{\mathbf{x} \in \mathcal{X}} \ \max_{\mathbb{P}_{\text{I}} \in \mathcal{D}_{\textup{M}}} \ \mathbb{E}_{\mathbb{P}_{\text{I}}} \left[ f\left(\mathbf{x},\mathbf{U}\boldsymbol{\Lambda}^{{\frac{1}{2}}} \boldsymbol{\xi}_{ \textup{I} } +\boldsymbol{\mu}\right) \right],
\end{align*}
which demonstrates that the optimal value of Problem \eqref{Equ:upperbound-setform} is an upper bound for that of Problem \eqref{Equ:DROM2} (i.e.,  Problem \eqref{Equ:DROM1}).

(ii) Consider any $m_1 < m_2 \leq m$. 
We have $\mathcal{B}_{m_2} := \{\mathbf{B} \in \mathbb{R}^{m \times m_2} \ | \ \mathbf{B}^{\top} \mathbf{B} = \mathbf{I}_{m_2} \}$ and 
consider an optimal solution $(\mathbf{B}^*, \mathbf{x}^*)$ of Problem \eqref{Equ:upperbound-setform}, 
i.e., $\min\limits_{\mathbf{B}  \in \mathcal{B}_{m_1}} \ \min\limits_{\mathbf{x} \in \mathcal{X}} \ \max\limits_{\mathbb{P}_{\text{I}}  \in \mathcal{D}_{\textup{U}} } \ \mathbb{E}_{\mathbb{P}_{\text{I}} } [ f (\mathbf{x}, \mathbf{U}\boldsymbol{\Lambda}^{{\frac{1}{2}}} \boldsymbol{\xi}_{ \textup{I} } +\boldsymbol{\mu} )  ]$.

Note that $(\mathbf{B}^*)^{\top} \mathbf{B}^* = \mathbf{I}_{m_1}$. 
We can then construct $\mathbf{B}^{\prime} = [\mathbf{B}^*, \mathbf{C}] \in \mathbb{R}^{m \times m_2}$ such that $\mathbf{C} \in \mathbb{R}^{m\times (m_2-m_1)}$ and $\mathbf{B}^{\prime} \in \mathcal{B}_{m_2}$, i.e., $ (\mathbf{B}^{\prime})^{\top} \mathbf{B}^{\prime} = \mathbf{I}_{m_2} $.
With $\mathbf{B}^{\prime}$, we use $\mathcal{D}_{\textup{U}}^{\prime}$ to denote the corresponding ambiguity set defined in \eqref{eqn:upperbound-ambiguity-set}.
By the second-moment constraint in $\mathcal{D}_{\textup{U}}^{\prime}$, we have
\begin{align*}
&\mathbb{E}_{\mathbb{P}_{\text{I}}} \left[ (\mathbf{B}^{\prime})^{\top} \boldsymbol{\xi}_{\textup{I}} \boldsymbol{\xi}_{\textup{I}}^{\top} \mathbf{B}^{\prime} \right] \\
= & \mathbb{E}_{\mathbb{P}_{\text{I}}} \left[ [\mathbf{B}^*, \mathbf{C}]^{\top} \boldsymbol{\xi}_{\textup{I}} \boldsymbol{\xi}_{\textup{I}}^{\top} [\mathbf{B}^*, \mathbf{C}] \right]  \\
= & \mathbb{E}_{\mathbb{P}_{\text{I}}} 
\begin{bmatrix}
    (\mathbf{B}^*)^{\top} \boldsymbol{\xi}_{\textup{I}} \boldsymbol{\xi}_{\textup{I}}^{\top} \mathbf{B}^* & (\mathbf{B}^*)^{\top} \boldsymbol{\xi}_{\textup{I}} \boldsymbol{\xi}_{\textup{I}}^{\top} \mathbf{C} \\
    \mathbf{C}^{\top} \boldsymbol{\xi}_{\textup{I}} \boldsymbol{\xi}_{\textup{I}}^{\top} \mathbf{B}^* & \mathbf{C}^{\top} \boldsymbol{\xi}_{\textup{I}} \boldsymbol{\xi}_{\textup{I}}^{\top} \mathbf{C}
\end{bmatrix}   \\
\preceq & \gamma_2 \mathbf{I}_{m_2},
\end{align*}
which implies that $\mathbb{E}_{\mathbb{P}_{\text{I}}}  [ (\mathbf{B}^*)^{\top} \boldsymbol{\xi}_{\textup{I}} \boldsymbol{\xi}_{\textup{I}}^{\top} \mathbf{B}^* ] \preceq \gamma_2 \mathbf{I}_{m_1}$.
It follows that $\mathcal{D}_{\textup{U}}^{\prime} \subseteq \mathcal{D}_{\textup{U}}$.
Therefore, we have
\begin{align}
\max\limits_{\mathbb{P}_{\text{I}} \in \mathcal{D}_{\textup{U}} } \ \mathbb{E}_{\mathbb{P}_{\text{I}}} \left[ f\left(\mathbf{x}^*,\mathbf{U}\boldsymbol{\Lambda}^{{\frac{1}{2}}} \boldsymbol{\xi}_{ \textup{I} } +\boldsymbol{\mu}\right) \right]  
    \geq \max_{\mathbb{P}_{\text{I}} \in \mathcal{D}_{\textup{U}}^{\prime}} \ \mathbb{E}_{\mathbb{P}_{\text{I}}} \left[ f\left(\mathbf{x}^*,\mathbf{U}\boldsymbol{\Lambda}^{{\frac{1}{2}}} \boldsymbol{\xi}_{ \textup{I} } +\boldsymbol{\mu}\right) \right]. \label{Cons:The3-1-1}
\end{align} 
Because $\mathbf{B}^{\prime} \in \mathcal{B}_{m_2}$ and $\mathbf{x}^* \in \mathcal{X}$, the constructed solution $(\mathbf{B}^{\prime}, \mathbf{x}^*)$ is feasible to the problem 
$
\min\limits_{\mathbf{B}  \in \mathcal{B}_{m_2}} \ \min\limits_{\mathbf{x} \in \mathcal{X}} \ \max\limits_{\mathbb{P}_{\text{I}}  \in \mathcal{D}_{\textup{U}}^{\prime} } \ \mathbb{E}_{ \mathbb{P}_{\text{I}} } [ f (\mathbf{x},\mathbf{U}\boldsymbol{\Lambda}^{{\frac{1}{2}}} \boldsymbol{\xi}_{ \textup{I} } +\boldsymbol{\mu} ) ].
$
Then, we have
\begin{align*}
\Theta_{\textup{U}}(m_2) 
& = \min\limits_{\mathbf{B}  \in \mathcal{B}_{m_2}} \ \min\limits_{\mathbf{x} \in \mathcal{X}} \ \max\limits_{\mathbb{P}_{\text{I}}  \in \mathcal{D}_{\textup{U}}^{\prime} } \ 
\mathbb{E}_{ \mathbb{P}_{\text{I}} } \left[ f\left(\mathbf{x}, \mathbf{U}\boldsymbol{\Lambda}^{{\frac{1}{2}}} \boldsymbol{\xi}_{ \textup{I} }  +\boldsymbol{\mu}\right) \right] \\
& \leq \max\limits_{\mathbb{P}_{\text{I}}  \in \mathcal{D}_{\textup{U}}^{\prime} } \ 
\mathbb{E}_{ \mathbb{P}_{\text{I}} } \left[ f\left(\mathbf{x}^*,\mathbf{U}\boldsymbol{\Lambda}^{{\frac{1}{2}}} \boldsymbol{\xi}_{ \textup{I} } +\boldsymbol{\mu}\right) \right]  \\ 
& \leq \max\limits_{\mathbb{P}_{\text{I}} \in \mathcal{D}_{\textup{U}} } \ \mathbb{E}_{\mathbb{P}_{\text{I}}}  \left[ f\left(\mathbf{x}^*,\mathbf{U}\boldsymbol{\Lambda}^{{\frac{1}{2}}} \boldsymbol{\xi}_{ \textup{I} } +\boldsymbol{\mu}\right) \right]  \\
& = \min\limits_{\mathbf{B}  \in \mathcal{B}_{m_1}} \ \min\limits_{\mathbf{x} \in \mathcal{X}} \ \max\limits_{\mathbb{P}_{\text{I}} \in \mathcal{D}_{\textup{U}} } \ \mathbb{E}_{\mathbb{P}_{\text{I}}}  \left[ f\left(\mathbf{x},\mathbf{U}\boldsymbol{\Lambda}^{{\frac{1}{2}}} \boldsymbol{\xi}_{ \textup{I} } +\boldsymbol{\mu}\right) \right]  \\ 
&= \Theta_{\textup{U}}(m_1), 
\end{align*}
where the first inequality holds because $(\mathbf{B}^{\prime}, \mathbf{x}^*)$ is a feasible solution of the problem $
\min\limits_{\mathbf{B}  \in \mathcal{B}_{m_2}} \ \min\limits_{\mathbf{x} \in \mathcal{X}} \ \max\limits_{\mathbb{P}_{\text{I}}  \in \mathcal{D}_{\textup{U}}^{\prime} } \ \mathbb{E}_{ \mathbb{P}_{\text{I}} } [ f (\mathbf{x},\mathbf{U}\boldsymbol{\Lambda}^{{\frac{1}{2}}} \boldsymbol{\xi}_{ \textup{I} } +\boldsymbol{\mu} ) ]
$, 
the second inequality holds by \eqref{Cons:The3-1-1}, 
and the second equality holds because $(\mathbf{B}^*, \mathbf{x}^*)$ is an optimal solution of Problem \eqref{Equ:upperbound-setform}.  
That is, the optimal value of Problem \eqref{Equ:upperbound-setform} is nonincreasing in $m_1$.

(iii) When $m_1 = m$, we have $\mathbf{B} \in \mathcal{B}_m \subseteq \mathbb{R}^{m \times m}$, i.e., $\mathbf{B}^{\top} \mathbf{B} = \mathbf{I}_m$. 
First, we have $\Theta_{\textup{U}}(m) \geq \Theta_{\textup{M}}(m)$ by the conclusion (i). 
Second, when $\mathbf{B} = \mathbf{I}_m$, Problem \eqref{Equ:upperbound-setform} becomes Problem \eqref{Equ:DROM2}. Because $\mathbf{B} = \mathbf{I}_m$ is a feasible solution of Problem \eqref{Equ:upperbound-setform}, it follows that $\Theta_{\textup{U}}(m) \leq \Theta_{\textup{M}}(m)$. 
Therefore, we have $\Theta_{\textup{U}}(m) = \Theta_{\textup{M}}(m)$.
\Halmos

\subsection{Proof of Proposition \ref{prop:ub-sdp-equivalent}}
First, by Theorem 3 in \cite{cheramin2022computationally}, Problem \eqref{Equ:upperbound-setform-2} has the same optimal value as the following problem:
\begin{subeqnarray} \label{Equ:upperboundRef1}
\min\limits_{\mathbf{x}, s, \mathbf{q}, \mathbf{Q}_{\textup{r}} } && s + \gamma_2\mathbf{I}_{m_1} \bullet \mathbf{Q}_{\textup{r}} + \sqrt{\gamma_1}\left \| \mathbf{q}  \right \|_2 \\	
\textnormal{s.t.} && s \geq f\left(\mathbf{x},\mathbf{U}\boldsymbol{\Lambda}^{ \frac{1}{2} }\boldsymbol{\xi}_{\textup{I}}+\boldsymbol{\mu}\right) - \boldsymbol{\xi}_{\textup{I}}^{\top}\mathbf{B}\mathbf{Q}_{\textup{r}} \mathbf{B}^{\top} \boldsymbol{\xi}_{\textup{I}}  - \mathbf{q}^{\top}  \boldsymbol{\xi}_{\textup{I}}, \  \forall \boldsymbol{\xi}_{\textup{I}}\in \mathcal{S}_{\textup{I}}, \slabel{Cons:upperboundRef1-1} \\ 
&&\mathbf{Q}_{\textup{r}}  \succeq 0, \ \mathbf{x} \in \mathcal{X}, \ \mathbf{Q}_{\textup{r}} \in \mathbb{R}^{m_1 \times m_1}, \ \mathbf{q} \in \mathbb{R}^{m}.	
\end{subeqnarray}

Next, we apply the strong duality theorem to constraints \eqref{Cons:upperboundRef1-1}. We define
\begin{align*}
g_k(\boldsymbol{\xi}_{\textup{I}}) = s + \boldsymbol{\xi}_{\textup{I}}^{\top}\mathbf{B}\mathbf{Q}_{\textup{r}} \mathbf{B}^{\top} \boldsymbol{\xi}_{\textup{I}}  + \mathbf{q}^{\top}  \boldsymbol{\xi}_{\textup{I}} - y^0_k(\mathbf{x}) - y_k(\mathbf{x})^{\top} \left(\mathbf{U}\boldsymbol{\Lambda}^{ \frac{1}{2} }\boldsymbol{\xi}_{\textup{I}}+\boldsymbol{\mu} \right), \ \forall k \in [K].
\end{align*}
As function $f(\mathbf{x},\boldsymbol{\xi})$ is piecewise linear convex, we can reformulate \eqref{Cons:upperboundRef1-1} as
\begin{align}
g_k(\boldsymbol{\xi}_{\textup{I}}) \geq 0, \  \forall \boldsymbol{\xi}_{\textup{I}}\in \mathcal{S}_{\textup{I}}, \ \forall k \in [K], \nonumber
\end{align}
which is equivalent to 
\begin{align}
    \min\limits_{\mathbf{A}\left(\mathbf{U} \boldsymbol{\Lambda}^{ \frac{1}{2} } \boldsymbol{\xi}_{\textup{I}} + \boldsymbol{\mu} \right) \leq  \mathbf{b}, \ \boldsymbol{\xi}_{\textup{I}} \in \mathbb{R}^m } \ g_k(\boldsymbol{\xi}_{\textup{I}}) \geq 0, \ \forall k \in [K]. \label{Cons:Prop-SDP-1}
\end{align}
For any $k \in [K]$, the Lagrangian dual problem of $\min_{\mathbf{A} (\mathbf{U}\boldsymbol{\Lambda}^{ \frac{1}{2} }\boldsymbol{\xi}_{\textup{I}}+\boldsymbol{\mu} ) \leq  \mathbf{b}, \ \boldsymbol{\xi}_{\textup{I}} \in \mathbb{R}^m } \ g_k(\boldsymbol{\xi}_{\textup{I}})$ is 
\begin{align}
\max\limits_{\boldsymbol{\lambda}_k \geq 0} \ \min\limits_{\boldsymbol{\xi}_{\textup{I}} \in \mathbb{R}^{m}} \ g_k(\boldsymbol{\xi}_{\textup{I}}) + \boldsymbol{\lambda}_k^{\top} \left( \mathbf{A}\left(\mathbf{U}\boldsymbol{\Lambda}^{ \frac{1}{2} }\boldsymbol{\xi}_{\textup{I}}+\boldsymbol{\mu} \right) - \mathbf{b} \right), \nonumber
\end{align}
where $\boldsymbol{\lambda}_k \in \mathbb{R}^{l}$.
Because there exists an interior point for the primal problem, the strong duality holds. 
Thus, constraints \eqref{Cons:Prop-SDP-1} are equivalent to 
\begin{align}
\max\limits_{\boldsymbol{\lambda}_k \geq 0} \ \min\limits_{\boldsymbol{\xi}_{\textup{I}}} \ g_k(\boldsymbol{\xi}_{\textup{I}}) + \boldsymbol{\lambda}_k^{\top} \left( \mathbf{A}\left(\mathbf{U}\boldsymbol{\Lambda}^{ \frac{1}{2} }\boldsymbol{\xi}_{\textup{I}}+\boldsymbol{\mu} \right) - \mathbf{b} \right) \geq 0, \ \forall k \in [K], \nonumber
\end{align}
which are further equivalent to
\begin{align}
\exists \boldsymbol{\lambda}_k \geq 0: \ & s + \boldsymbol{\xi}_{\textup{I}}^{\top}\mathbf{B}\mathbf{Q}_{\textup{r}} \mathbf{B}^{\top} \boldsymbol{\xi}_{\textup{I}}  + \mathbf{q}^{\top} \boldsymbol{\xi}_{\textup{I}} - y^0_k(\mathbf{x}) - y_k(\mathbf{x})^{\top} \left(\mathbf{U}\boldsymbol{\Lambda}^{ \frac{1}{2} }\boldsymbol{\xi}_{\textup{I}}+\boldsymbol{\mu} \right) \nonumber \\ 
&\hspace{1cm} + \boldsymbol{\lambda}_k^{\top} \left( \mathbf{A}\left(\mathbf{U}\boldsymbol{\Lambda}^{ \frac{1}{2} }\boldsymbol{\xi}_{\textup{I}}+\boldsymbol{\mu} \right) - \mathbf{b} \right) \geq 0,\  \forall \boldsymbol{\xi}_{\textup{I}}\in \mathbb{R}^m, \ \forall k \in [K]. \label{Cons:Prop-SDP-2}
\end{align}

Note that $\mathbf{B}^{\top} \mathbf{B} = \mathbf{I}_{m_1}$; that is, all the column vectors of $\mathbf{B}$ are orthogonal. 
We can then extend $\mathbf{B}$ to $[\mathbf{B},\bar{\mathbf{B}} ] \in \mathbb{R}^{m\times m}$ with $\bar{\mathbf{B}} \in \mathbb{R}^{m\times (m-m_1)}$ such that all the column vectors of $[\mathbf{B},\bar{\mathbf{B}} ]$ span the space of $\mathbb{R}^m$. 
Thus, we can always find $\boldsymbol{\xi}_1 \in \mathbb{R}^{m_1}$ and $\boldsymbol{\xi}_2 \in \mathbb{R}^{m-m_1}$ such that 
\begin{align*}
\boldsymbol{\xi}_{\textup{I}} = \mathbf{B} \boldsymbol{\xi}_1 + \bar{\mathbf{B}} \boldsymbol{\xi}_2.
\end{align*}

\noindent
It follows that constraints \eqref{Cons:Prop-SDP-2} become
\begin{align}
\exists \boldsymbol{\lambda}_k \geq 0: \ & s + \boldsymbol{\xi}_1^{\top} \mathbf{Q}_{\textup{r}} \boldsymbol{\xi}_1  + \mathbf{q}^{\top}\left( \mathbf{B} \boldsymbol{\xi}_1 + \bar{\mathbf{B}} \boldsymbol{\xi}_2 \right) - y^0_k(\mathbf{x}) - y_k(\mathbf{x})^{\top} \left(\mathbf{U}\boldsymbol{\Lambda}^{ \frac{1}{2} }\left(\mathbf{B} \boldsymbol{\xi}_1 + \bar{\mathbf{B}} \boldsymbol{\xi}_2 \right)+\boldsymbol{\mu} \right) \nonumber \\ 
&\hspace{0.5cm} + \boldsymbol{\lambda}_k^{\top} \left( \mathbf{A}\left(\mathbf{U}\boldsymbol{\Lambda}^{ \frac{1}{2} }\left(\mathbf{B} \boldsymbol{\xi}_1 + \bar{\mathbf{B}} \boldsymbol{\xi}_2 \right)+\boldsymbol{\mu} \right) - \mathbf{b} \right) \geq 0,\  \forall \boldsymbol{\xi}_1 \in \mathbb{R}^{m_1}, \boldsymbol{\xi}_2 \in \mathbb{R}^{m-m_1}, \ \forall k \in [K]. \label{Cons:Prop-SDP-2.5}
\end{align}

\noindent
We further define 
\begin{align*}
\mathbf{Z}_k = 
\begin{bmatrix}  
s-y_k^0(\mathbf{x})-\boldsymbol{\lambda}_k^{\top}\mathbf{b}-y_k(\mathbf{x})^{\top}\boldsymbol{\mu} +\boldsymbol{\lambda}_k^{\top}\mathbf{A}\boldsymbol{\mu} & \hspace{0.1 in} \frac{1}{2}  \left(\mathbf{B}^{\top}\mathbf{q} +\left(\mathbf{U}\boldsymbol{\Lambda}^{{\frac{1}{2}}} \mathbf{B} \right)^{\top}\left(\mathbf{A}^{\top}\boldsymbol{\lambda}_k-y_k(\mathbf{x})\right) \right)^{\top} \\ 
	 \frac{1}{2} \left( \mathbf{B}^{\top}\mathbf{q} +\left(\mathbf{U}\boldsymbol{\Lambda}^{{\frac{1}{2}}} \mathbf{B} \right)^{\top}\left(\mathbf{A}^{\top}\boldsymbol{\lambda}_k-y_k(\mathbf{x})\right)  \right) & \mathbf{Q}_{ \textup{r} } 
\end{bmatrix}, \ \forall k \in [K].
\end{align*}

\noindent
Thus, we have
\begin{align}
\eqref{Cons:Prop-SDP-2.5} & \iff \exists \boldsymbol{\lambda}_k \geq 0: \ \left( 1,\boldsymbol{\xi}_1^{\top} \right) \mathbf{Z}_k \left( 1,\boldsymbol{\xi}_1^{\top} \right)^{\top} +  \boldsymbol{\xi}_2^{\top} \left( \bar{\mathbf{B}}^{\top} \mathbf{q} + \left(\mathbf{U}\boldsymbol{\Lambda}^{{\frac{1}{2}}} \bar{\mathbf{B}} \right)^{\top} \left(\mathbf{A}^{\top}\boldsymbol{\lambda}_k-y_k(\mathbf{x})\right) \right) \geq 0, \nonumber \\
& \hspace{8cm} \forall \boldsymbol{\xi}_1 \in \mathbb{R}^{m_1}, \boldsymbol{\xi}_2 \in \mathbb{R}^{m-m_1}, \ \forall k \in [K]. \nonumber \\ 
& \iff \exists \boldsymbol{\lambda}_k \geq 0: \ \left( 1,\boldsymbol{\xi}_1^{\top} \right) \mathbf{Z}_k \left( 1,\boldsymbol{\xi}_1^{\top} \right)^{\top} \geq 0, \ \forall \boldsymbol{\xi}_1 \in \mathbb{R}^{m_1}, \ \forall k \in [K]; \label{Cons:Prop-SDP-3} \\
&  \hspace{2.5cm} \bar{\mathbf{B}}^{\top} \mathbf{q} +\left(\mathbf{U}\boldsymbol{\Lambda}^{{\frac{1}{2}}} \bar{\mathbf{B}} \right)^{\top}\left(\mathbf{A}^{\top}\boldsymbol{\lambda}_k-y_k(\mathbf{x})\right) = 0,\ \forall k \in [K]. \nonumber \\
&\iff  \exists \boldsymbol{\lambda}_k \geq 0: \  \mathbf{Z}_k \succeq 0,\  \bar{\mathbf{B}}^{\top} \mathbf{q} +\left(\mathbf{U}\boldsymbol{\Lambda}^{{\frac{1}{2}}} \bar{\mathbf{B}} \right)^{\top}\left(\mathbf{A}^{\top}\boldsymbol{\lambda}_k-y_k(\mathbf{x})\right) = 0,\ \forall k \in [K]. \nonumber \\ 
&\iff  \exists \boldsymbol{\lambda}_k \geq 0: \  \mathbf{Z}_k \succeq 0,\  \bar{\mathbf{B}}^{\top} \left(\mathbf{q}+\left(\mathbf{U}\boldsymbol{\Lambda}^{{\frac{1}{2}}}  \right)^{\top}\left(\mathbf{A}^{\top}\boldsymbol{\lambda}_k-y_k(\mathbf{x})\right) \right) = 0,\ \forall k \in [K]. \label{Cons:Prop-SDP-4} \\ 
&\iff  \exists \boldsymbol{\lambda}_k \geq 0, \mathbf{u}_k \in \mathbb{R}^{m_1}: \ \mathbf{Z}_k \succeq 0, \ \mathbf{q}+\left(\mathbf{U}\boldsymbol{\Lambda}^{{\frac{1}{2}}}  \right)^{\top}\left(\mathbf{A}^{\top}\boldsymbol{\lambda}_k-y_k(\mathbf{x})\right) = \mathbf{B} \mathbf{u}_k,\ \forall k \in [K]. \label{Cons:Prop-sdp-end}
\end{align}
The first equivalence holds due to the definition of $\mathbf{Z}_k$.
For the third equivalence, clearly $\Longleftarrow$ follows from the definition of a PSD matrix.
To prove $\Longrightarrow$, we consider two possible cases for any  $(\eta_0  \in \mathbb{R}, \boldsymbol{\eta}^{\top}\in \mathbb{R}^{m_1})^{\top} \in \mathbb{R}^{m_1+1}$:
(i) if $\eta_0=0$, then  $(\eta_0, \boldsymbol{\eta}^{\top}) \mathbf{Z}_k (\eta_0, \boldsymbol{\eta}^{\top})^{\top} = \boldsymbol{\eta}^{\top} \mathbf{Q}_{\textup{r}} \boldsymbol{\eta} \geq 0$ because $\mathbf{Q}_{\textup{r}}$ is PSD;
(ii) if $\eta_0 \neq 0$, then we have $(\eta_0, \boldsymbol{\eta}^{\top}) \mathbf{Z}_k (\eta_0,  \boldsymbol{\eta}^{\top})^{\top} = \eta_0^2 ( 1, \frac{\boldsymbol{\eta}^{\top}}{\eta_0} ) \mathbf{Z}_k ( 1, \frac{\boldsymbol{\eta}^{\top}}{\eta_0} )^{\top} \geq 0$ according to \eqref{Cons:Prop-SDP-3}. Therefore, $\Longrightarrow$ holds.
For the fifth equivalence, \eqref{Cons:Prop-SDP-4} shows that $\mathbf{q} + (\mathbf{U} \boldsymbol{\Lambda}^{{\frac{1}{2}}}  )^{\top} (\mathbf{A}^{\top} \boldsymbol{\lambda}_k - y_k(\mathbf{x}))$ is in the null space of $\bar{\mathbf{B}}$ and thus cannot be represented by basis vectors in the space of $\bar{\mathbf{B}}$.
Because $[\mathbf{B},\bar{\mathbf{B}} ]$ span the space of $\mathbb{R}^m$, we have $\mathbf{q} + (\mathbf{U} \boldsymbol{\Lambda}^{{\frac{1}{2}}}  )^{\top} (\mathbf{A}^{\top} \boldsymbol{\lambda}_k - y_k(\mathbf{x}))$ should be in the space of $ \mathbf{B}$.
That is, there exists $\mathbf{u}_k \in \mathbb{R}^{m_1} $ such that $\mathbf{q}+(\mathbf{U}\boldsymbol{\Lambda}^{{\frac{1}{2}}}  )^{\top}(\mathbf{A}^{\top}\boldsymbol{\lambda}_k-y_k(\mathbf{x})) = \mathbf{B} \mathbf{u}_k$ for any $k \in [K]$.
Meanwhile, because $\mathbf{B}^{\top} \mathbf{B} = \mathbf{I}_{m_1}$, we have 
\begin{align}
\mathbf{B}^{\top}\mathbf{q} +\left(\mathbf{U}\boldsymbol{\Lambda}^{{\frac{1}{2}}} \mathbf{B} \right)^{\top}\left(\mathbf{A}^{\top}\boldsymbol{\lambda}_k-y_k(\mathbf{x})\right) = \mathbf{B}^{\top} \mathbf{B}\mathbf{u}_k=\mathbf{u}_k, \ \forall k \in [K]. \nonumber
\end{align}

Finally, we obtain Problem \eqref{Equ:upperbound-sdpform} by replacing constraints \eqref{Cons:upperboundRef1-1} with \eqref{Cons:Prop-sdp-end} and replacing $\mathbf{B}^{\top}\mathbf{q} + (\mathbf{U}\boldsymbol{\Lambda}^{{\frac{1}{2}}} \mathbf{B} )^{\top} (\mathbf{A}^{\top} \boldsymbol{\lambda}_k - y_k(\mathbf{x}) )$ with $\mathbf{u}_k$.
\Halmos

\subsection{Proof of Theorem \ref{Theo:upperboundSameoptimal}}
Consider $m_1 = K$.
We construct a solution $(\mathbf{x}^{\dag}, s^{\dag}, \hat{\boldsymbol{\lambda}}^{\dag}, \mathbf{q}^{\dag}, \mathbf{Q}_{\textup{r}}^{\dag}, \hat{\mathbf{u}}^{\dag}, \mathbf{B}^{\dag})$ of Problem \eqref{Equ:upperbound-setform}
by setting 
$\mathbf{x}^{\dag} = \mathbf{x}^*$, 
$s^{\dag} = s^*$, 
$\hat{\boldsymbol{\lambda}}^{\dag} = \hat{\boldsymbol{\lambda}}^*$,
$\mathbf{q}^{\dag} =\mathbf{q}^{\prime}= \mathbf{V}\boldsymbol{\delta}$,
$\mathbf{Q}_{\textup{r}}^{\dag} = \mathbf{Y}_{11}$, $\mathbf{B}^{\dag} = \mathbf{V}$, and 
$\hat{\mathbf{u}}_k^{\dag} = \boldsymbol{\delta} + \boldsymbol{\nu}_k \ (k \in [K])$.

First, we show this constructed solution is feasible to Problem \eqref{Equ:upperbound-setform}.
Clearly, this solution satisfies constraints \eqref{Cons:upperbound-sdpform-2.5}--\eqref{Cons:upperbound-sdpform-3}. 
By the construction of the solution, for any  $ k \in [K] $, we further have 
\begin{align}
\mathbf{q}^{\dag}+\left(\mathbf{U}\boldsymbol{\Lambda}^{{\frac{1}{2}}}\right)^{\top}\left(\mathbf{A}^{\top}\boldsymbol{\lambda}^{\dag}_k-y_k(\mathbf{x}^{\dag})\right) & =  \mathbf{V}  \boldsymbol{\delta} +\left(\mathbf{U}\boldsymbol{\Lambda}^{{\frac{1}{2}}}\right)^{\top}\left(\mathbf{A}^{\top}\boldsymbol{\lambda}^*_k-y_k(\mathbf{x}^*)\right) \nonumber \\ 
& = \mathbf{V}  \boldsymbol{\delta} + \mathbf{V} \boldsymbol{\nu}_k  = \mathbf{B}^{\dag} \hat{\mathbf{u}}_k^{\dag}, \nonumber
\end{align}
where the first equality holds by the construction of $\mathbf{q}^{\dag}$, the second equality holds by \eqref{Cons:Prop2-2}, and the third equality holds by the construction of $\hat{\mathbf{u}}_k^{\dag} $.
Thus, this solution satisfies constraints \eqref{Cons:upperbound-sdpform-2}. 
Meanwhile, $ \mathbf{V}^{\top} \mathbf{V} = \mathbf{I}_{K} =  \mathbf{I}_{m_1} $.
It follows that $ (\mathbf{B}^{\dag})^{\top} {\mathbf{B}^{\dag}} =  \mathbf{I}_{m_1} $. 

In addition, from Problem \eqref{Equ:MainProblem}, as $ \mathbf{q}^{\prime} = \mathbf{V} \boldsymbol{\delta} $ and $\mathbf{Q}^{\prime} = \mathbf{V} \mathbf{Y}_{11} \mathbf{V}^{\top}$, for any $k \in [K]$,  we have
\begin{align*}
\begin{bmatrix}  S_k & \hspace{0.1 in}  \frac{1}{2}  \left(\mathbf{V}  \boldsymbol{\delta} + \left( \mathbf{U}\boldsymbol{\Lambda}^{{\frac{1}{2}}} \right)^{\top} \left( \mathbf{A}^{\top}\boldsymbol{\lambda}^*_k - y_k(\mathbf{x}^*) \right) \right)^{\top} \\ 
 \frac{1}{2} \left( \mathbf{V}  \boldsymbol{\delta} + \left( \mathbf{U} \boldsymbol{\Lambda}^{{\frac{1}{2}}} \right)^{\top} \left( \mathbf{A}^{\top} \boldsymbol{\lambda}^*_k - y_k(\mathbf{x}^*) \right) \right) & \mathbf{V} \mathbf{Y}_{11} \mathbf{V}^{\top} \end{bmatrix} \succeq 0, 
\end{align*}
which, by Schur complement, is equivalent to
\begin{align}
4 S_k \left( \mathbf{V} \mathbf{Y}_{11} \mathbf{V}^{\top} \right) & \succeq    \left(\mathbf{V}  \boldsymbol{\delta} +\left(\mathbf{U}\boldsymbol{\Lambda}^{{\frac{1}{2}}}\right)^{\top}\left(\mathbf{A}^{\top}\boldsymbol{\lambda}^*_k-y_k(\mathbf{x}^*)\right)\right) 
\left( \mathbf{V}  \boldsymbol{\delta} +\left(\mathbf{U}\boldsymbol{\Lambda}^{{\frac{1}{2}}}\right)^{\top}\left(\mathbf{A}^{\top}\boldsymbol{\lambda}^*_k-y_k(\mathbf{x}^*)\right) \right)^{\top} \nonumber \\ 
&  = \left(\mathbf{V} \mathbf{u}_k^{\dag}\right) 
\left( \mathbf{V} \mathbf{u}_k^{\dag} \right)^{\top}, \label{eqn:theorem-3-1}
\end{align}
where the equality holds by \eqref{Cons:Prop2-2} and the construction of $\hat{\mathbf{u}}_k^{\dag} $.
From \eqref{eqn:theorem-3-1}, for any $k \in [K]$, we have the following inequality holds by Lemma \ref{Lem:VXV}:
\begin{small}
\begin{equation}
4S_k\left( \mathbf{V}^{\top} \mathbf{V} \mathbf{Y}_{11} \mathbf{V}^{\top} \mathbf{V} \right) \succeq    \mathbf{V}^{\top}\left(\mathbf{V} \mathbf{u}_k^{\dag}\right) 
\left( \mathbf{V} \mathbf{u}_k^{\dag} \right)^{\top} \mathbf{V}, \nonumber
\end{equation}
\end{small}%
which is equivalent to 
\begin{align}
4S_k \mathbf{Y}_{11} \succeq  \mathbf{u}_k^{\dag} \mathbf{u}_k^{\dag \top}   \label{eqn:theorem-3-2}
\end{align}
because $\mathbf{V}^{\top} \mathbf{V} =  \mathbf{I}_{K}$. 
By Schur complement, for any $k \in [K]$, \eqref{eqn:theorem-3-2} further becomes
\begin{align*}
&\begin{bmatrix}  S_k &   \frac{1}{2}\mathbf{u}_k^{\dag \top}\\ 
\frac{1}{2}  \mathbf{u}_k^{\dag} & \mathbf{Y}_{11}  \end{bmatrix} \succeq 0, 
\end{align*}
which indicates that the constructed solution 
$(\mathbf{x}^{\dag}, s^{\dag}, \hat{\boldsymbol{\lambda}}^{\dag}, \mathbf{q}^{\dag}, \mathbf{Q}_{\textup{r}}^{\dag}, \hat{\mathbf{u}}^{\dag}, \mathbf{B}^{\dag})$ also satisfies constraints \eqref{Cons:upperbound-sdpform-1} 
and thus it is a feasible solution of Problem \eqref{Equ:upperbound-setform}.

Second, we show this feasible solution $(\mathbf{x}^{\dag}, s^{\dag}, \hat{\boldsymbol{\lambda}}^{\dag}, \mathbf{q}^{\dag}, \mathbf{Q}_{\textup{r}}^{\dag}, \hat{\mathbf{u}}^{\dag}, \mathbf{B}^{\dag})$ is an optimal solution of Problem \eqref{Equ:upperbound-setform}.
The objective value corresponding to this solution is
\begin{align}
 & s^{\dag} + \gamma_2 \mathbf{I}_{m_1} \bullet \mathbf{Q}_{\textup{r}}^{\dag}  + \sqrt{\gamma_1}\left \| \mathbf{q}^{\dag} \right \|_2 = s^* +\gamma_2\mathbf{I}_{m_1} \bullet \mathbf{Y}_{11} + \sqrt{\gamma_1}\left \| \mathbf{q}^{\prime} \right \|_2 \nonumber\\
 = \ & s^* +\gamma_2\mathbf{I}_{m_1} \bullet (\mathbf{Y}_{11}\mathbf{V}^{\top} \mathbf{V}) + \sqrt{\gamma_1}\left \| \mathbf{q}^{\prime} \right \|_2 = s^* +\gamma_2\mathbf{I}_{m} \bullet (\mathbf{V}\mathbf{Y}_{11}\mathbf{V}^{\top}) + \sqrt{\gamma_1}\left \| \mathbf{q}^{\prime} \right \|_2 \nonumber\\ 
 = \ & s^* +\gamma_2\mathbf{I}_{m} \bullet \mathbf{Q}^{\prime} + \sqrt{\gamma_1}\left \| \mathbf{q}^{\prime} \right \|_2 = \Theta_{\textup{M}}(m), \nonumber 
\end{align}
where the first equality holds by the construction of $(\mathbf{x}^{\dag}, s^{\dag}, \hat{\boldsymbol{\lambda}}^{\dag}, \mathbf{q}^{\dag}, \mathbf{Q}_{\textup{r}}^{\dag}, \hat{\mathbf{u}}^{\dag}, \mathbf{B}^{\dag})$,
the second equality holds because $\mathbf{V}^{\top} \mathbf{V} = \mathbf{I}_{K}$,
the third equality holds by the cyclic property of a matrix's trace, 
and the fourth equality holds by the definition of 
$\mathbf{Q}^{\prime}$ in Theorem \ref{Prop:lowrank}. 
Therefore, the solution $(\mathbf{x}^{\dag}, s^{\dag}, \hat{\boldsymbol{\lambda}}^{\dag}, \mathbf{q}^{\dag}, \mathbf{Q}_{\textup{r}}^{\dag}, \hat{\mathbf{u}}^{\dag}, \mathbf{B}^{\dag})$ is an optimal solution of Problem \eqref{Equ:upperbound-setform}.

Finally, when $ m_1 >K $ and $m_1 \leq m$, we have $\Theta_{\textup{U}}(m_1) \geq \Theta_{\textup{M}}(m)$ by the conclusion (i) in Theorem \ref{Theo:upperbound} and
$\Theta_{\textup{U}}(m_1) \leq \Theta_{\textup{U}}(K) = \Theta_{\textup{M}}(m)$ by the conclusion (ii) in Theorem \ref{Theo:upperbound}.
It follows that $\Theta_{\textup{U}}(m_1) = \Theta_{\textup{M}}(m)$.
\Halmos

\subsection{Proof of Proposition \ref{Lem:VXVwithRank}}
First, by Lemma \ref{Lem:VXV}, for any $\mathbf{B} \in \mathcal{B}_{m_1}$, we have 
\begin{align}
\mathbf{X} \preceq \mathbf{I}_{m} \Longrightarrow \mathbf{B}^{\top} \mathbf{X} \mathbf{B} \preceq \mathbf{B}^{\top} \mathbf{I}_{m} \mathbf{B} = \mathbf{I}_{m_1}. \nonumber
\end{align}

Second, we perform eigenvalue decomposition on $\mathbf{X}$, i.e., $\mathbf{X} = \mathbf{Q} \boldsymbol{\Lambda} \mathbf{Q}^{\top}$,
where $\mathbf{Q} \in \mathbb{R}^{m\times m}$ is a matrix with orthonormal column vectors and $\boldsymbol{\Lambda} \in \mathbb{R}^{m\times m}$ is a diagonal matrix. Without loss of generality, we assume that the diagonal elements of $\boldsymbol{\Lambda}$ are arranged in a nonincreasing order and let $\boldsymbol{\Lambda}_{m_1 \times m_1}$ represent the upper-left submatrix of $\boldsymbol{\Lambda}$. 


Now we let $\mathbf{B} = \mathbf{Q}_{m\times m_1}$, where $\mathbf{Q}_{m\times m_1}$ is the left submatrix of $\mathbf{Q}$. Then we have $\mathbf{B} \in \mathcal{B}_{m_1}$ and 
\begin{align*}
\mathbf{B}^{\top} \mathbf{X} \mathbf{B} \preceq \mathbf{I}_{m_1}  & \Longrightarrow \mathbf{B}^{\top} \mathbf{Q} \boldsymbol{\Lambda} \mathbf{Q}^{\top} \mathbf{B} \preceq \mathbf{I}_{m_1}  \Longrightarrow \mathbf{Q}_{m\times m_1}^{\top} \mathbf{Q} \boldsymbol{\Lambda} \mathbf{Q}^{\top} \mathbf{Q}_{m\times m_1} \preceq \mathbf{I}_{m_1} \\
& \Longrightarrow [\mathbf{I}_{m_1}, \mathbf{0}_{m_1 \times (m-m_1)}] \boldsymbol{\Lambda}  [\mathbf{I}_{m_1}, \mathbf{0}_{m_1 \times (m-m_1)}]^{\top} \preceq \mathbf{I}_{m_1} \\ 
& \Longrightarrow  \boldsymbol{\Lambda}_{m_1 \times m_1}    \preceq \mathbf{I}_{m_1}   \Longrightarrow  \boldsymbol{\Lambda}  \preceq \mathbf{I}_{m}  \Longrightarrow \mathbf{Q} \boldsymbol{\Lambda} \mathbf{Q}^{\top}\preceq \mathbf{Q}\mathbf{I}_{m} \mathbf{Q}^{\top} \Longrightarrow \mathbf{X} \preceq \mathbf{I}_{m} 
\end{align*} 
where the first deduction holds by the eigenvalue decomposition of $\mathbf{X}$, 
the second deduction holds by the construction of $\mathbf{B}$, 
the third deduction holds because all the column vectors in $\mathbf{Q}$ are orthonormal, 
the fourth deduction holds by the definition of $\boldsymbol{\Lambda}_{m_1 \times m_1}$, 
the fifth deduction holds because ${\rm{rank}}(\mathbf{X}) \leq m_1$, 
the sixth deduction holds by Lemma \ref{Lem:VXV}.
Thus, if $\mathbf{B}^{\top} \mathbf{X} \mathbf{B} \preceq \mathbf{I}_{m_1} $ for any $ \mathbf{B} \in \mathcal{B}_{m_1}$, then we have $ \mathbf{X} \preceq \mathbf{I}_{m} $.
The proof is complete.
\Halmos

\section{Supplement to Section \ref{sec:lower-bound-2}}

\subsection{Proof of Theorem \ref{thm-new-lower-bound}}
By dualizing the inner maximization problem of Problem \eqref{Equ:lowerbound-setform-2} and integrating it with the outer minimization operators, we first obtain the following formulation:
\begin{subeqnarray} \label{Equ:lowerboundRef1}
\min\limits_{ \substack{ \mathbf{x}, s, \mathbf{q}, \mathbf{Q}_{\textup{r}}^{\prime}, \mathbf{Q}_{\textup{r}}^{\prime \prime} \\ \mathbf{B}_1, \mathbf{B}_2 } } && s + \gamma_2\mathbf{I}_{m_1} \bullet \mathbf{Q}_{\textup{r}}^{\prime} + \sqrt{\gamma_1} \left \| \mathbf{q}  \right \|_2 \\
\textnormal{s.t.} && s \geq f\left(\mathbf{x},\mathbf{U}\boldsymbol{\Lambda}^{ \frac{1}{2} }\boldsymbol{\xi}_{\textup{I}}+\boldsymbol{\mu}\right) - \boldsymbol{\xi}_{\textup{I}}^{\top}\mathbf{B}_1 \mathbf{Q}_{\textup{r}}^{\prime} \mathbf{B}_1^{\top} \boldsymbol{\xi}_{\textup{I}}  
- \boldsymbol{\xi}_{\textup{I}}^{\top}\mathbf{B}_2 \mathbf{Q}_{\textup{r}}^{\prime \prime} \mathbf{B}_2^{\top} \boldsymbol{\xi}_{\textup{I}} 
- \mathbf{q}^{\top} 
\boldsymbol{\xi}_{\textup{I}}, \  \forall \boldsymbol{\xi}_{\textup{I}}\in \mathcal{S}_{\textup{I}}, \slabel{Cons:lowerboundRef1-1} \\ 
&&\mathbf{Q}_{\textup{r}}^{\prime}  \succeq 0, \ \mathbf{Q}_{\textup{r}}^{\prime \prime}  \succeq 0, \ \mathbf{x} \in \mathcal{X}, \ \mathbf{Q}_{\textup{r}}^{\prime} \in \mathbb{R}^{m_1 \times m_1}, \ \mathbf{Q}_{\textup{r}}^{\prime \prime}  \in \mathbb{R}^{(K-m_1) \times (K-m_1)},\ \mathbf{q} \in \mathbb{R}^{m}, \\
&& \mathbf{B}_1 \in \mathbb{R}^{m \times m_1},  \mathbf{B}_2  \in \mathbb{R}^{m \times (K-m_1)}, \ [  \mathbf{B}_1,  \mathbf{B}_2 ]^{\top} [\mathbf{B}_1,  \mathbf{B}_2] = \mathbf{I}_K.
\end{subeqnarray}

Next, we apply the strong duality theorem to constraints \eqref{Cons:lowerboundRef1-1}. We define
\begin{align*}
g_k(\boldsymbol{\xi}_{\textup{I}}) = s + \boldsymbol{\xi}_{\textup{I}}^{\top}\mathbf{B}_1 \mathbf{Q}_{\textup{r}}^{\prime} \mathbf{B}_1^{\top} \boldsymbol{\xi}_{\textup{I}}  
+ \boldsymbol{\xi}_{\textup{I}}^{\top}\mathbf{B}_2 \mathbf{Q}_{\textup{r}}^{\prime \prime} \mathbf{B}_2^{\top} \boldsymbol{\xi}_{\textup{I}}  + \mathbf{q}^{\top}  \boldsymbol{\xi}_{\textup{I}} - y^0_k(\mathbf{x}) - y_k(\mathbf{x})^{\top} \left(\mathbf{U}\boldsymbol{\Lambda}^{ \frac{1}{2} }\boldsymbol{\xi}_{\textup{I}}+\boldsymbol{\mu} \right), \ \forall k \in [K].
\end{align*}
As function $f(\mathbf{x},\boldsymbol{\xi})$ is piecewise linear convex, we can reformulate \eqref{Cons:lowerboundRef1-1} as
\begin{align}
g_k(\boldsymbol{\xi}_{\textup{I}}) \geq 0, \  \forall \boldsymbol{\xi}_{\textup{I}}\in \mathcal{S}_{\textup{I}}, \ \forall k \in [K], \nonumber
\end{align}
which is equivalent to 
\begin{align}
    \min\limits_{\mathbf{A}\left(\mathbf{U} \boldsymbol{\Lambda}^{ \frac{1}{2} } \boldsymbol{\xi}_{\textup{I}} + \boldsymbol{\mu} \right) \leq  \mathbf{b}, \ \boldsymbol{\xi}_{\textup{I}} \in \mathbb{R}^m } \ g_k(\boldsymbol{\xi}_{\textup{I}}) \geq 0, \ \forall k \in [K]. \label{Cons:lowerboundProp-SDP-1}
\end{align}
For any $k \in [K]$, the Lagrangian dual problem of $\min_{\mathbf{A} (\mathbf{U}\boldsymbol{\Lambda}^{ \frac{1}{2} }\boldsymbol{\xi}_{\textup{I}}+\boldsymbol{\mu} ) \leq  \mathbf{b}, \ \boldsymbol{\xi}_{\textup{I}} \in \mathbb{R}^m } \ g_k(\boldsymbol{\xi}_{\textup{I}})$ is 
\begin{align}
\max\limits_{\boldsymbol{\lambda}_k \geq 0} \ \min\limits_{\boldsymbol{\xi}_{\textup{I}} \in \mathbb{R}^{m}} \ g_k(\boldsymbol{\xi}_{\textup{I}}) + \boldsymbol{\lambda}_k^{\top} \left( \mathbf{A}\left(\mathbf{U}\boldsymbol{\Lambda}^{ \frac{1}{2} }\boldsymbol{\xi}_{\textup{I}}+\boldsymbol{\mu} \right) - \mathbf{b} \right), \nonumber
\end{align}
where $\boldsymbol{\lambda}_k \in \mathbb{R}^{l}$.
Because there exists an interior point for the primal problem, the strong duality holds. 
Thus, constraints \eqref{Cons:lowerboundProp-SDP-1} are equivalent to 
\begin{align}
\max\limits_{\boldsymbol{\lambda}_k \geq 0} \ \min\limits_{\boldsymbol{\xi}_{\textup{I}}} \ g_k(\boldsymbol{\xi}_{\textup{I}}) + \boldsymbol{\lambda}_k^{\top} \left( \mathbf{A}\left(\mathbf{U}\boldsymbol{\Lambda}^{ \frac{1}{2} }\boldsymbol{\xi}_{\textup{I}}+\boldsymbol{\mu} \right) - \mathbf{b} \right) \geq 0, \ \forall k \in [K], \nonumber
\end{align}
which are further equivalent to
\begin{align}
\exists \boldsymbol{\lambda}_k \geq 0: \ & s + \boldsymbol{\xi}_{\textup{I}}^{\top}\mathbf{B}_1 \mathbf{Q}_{\textup{r}}^{\prime} \mathbf{B}_1^{\top} \boldsymbol{\xi}_{\textup{I}}  
+ \boldsymbol{\xi}_{\textup{I}}^{\top}\mathbf{B}_2 \mathbf{Q}_{\textup{r}}^{\prime \prime} \mathbf{B}_2^{\top} \boldsymbol{\xi}_{\textup{I}}  + \mathbf{q}^{\top}  \boldsymbol{\xi}_{\textup{I}} - y^0_k(\mathbf{x}) - y_k(\mathbf{x})^{\top} \left(\mathbf{U}\boldsymbol{\Lambda}^{ \frac{1}{2} }\boldsymbol{\xi}_{\textup{I}}+\boldsymbol{\mu} \right) \nonumber \\ 
&\hspace{1cm} + \boldsymbol{\lambda}_k^{\top} \left( \mathbf{A}\left(\mathbf{U}\boldsymbol{\Lambda}^{ \frac{1}{2} }\boldsymbol{\xi}_{\textup{I}}+\boldsymbol{\mu} \right) - \mathbf{b} \right) \geq 0,\  \forall \boldsymbol{\xi}_{\textup{I}}\in \mathbb{R}^m, \ \forall k \in [K]. \label{Cons:lowerboundProp-SDP-2}
\end{align}

Note that $\mathbf{B}^{\top} \mathbf{B} = \mathbf{I}_{m_1}$; that is, all the column vectors of $\mathbf{B}$ are orthogonal. 
We can then extend $\mathbf{B}$ to $[\mathbf{B},\bar{\mathbf{B}} ] \in \mathbb{R}^{m\times m}$ with $\bar{\mathbf{B}} \in \mathbb{R}^{m\times (m-K)}$ such that all the column vectors of $[\mathbf{B},\bar{\mathbf{B}} ]$ span the space of $\mathbb{R}^m$. 
Thus, we can always find $\boldsymbol{\xi}_1 \in \mathbb{R}^{m_1}$, $\boldsymbol{\xi}_2 \in \mathbb{R}^{K-m_1}$, and $\boldsymbol{\xi}_3 \in \mathbb{R}^{m-K}$ such that 
\begin{align*}
\boldsymbol{\xi}_{\textup{I}} = \mathbf{B}_1 \boldsymbol{\xi}_1 + \mathbf{B}_2 \boldsymbol{\xi}_2 + \bar{\mathbf{B}} \boldsymbol{\xi}_3.
\end{align*}

\noindent
It follows that constraints \eqref{Cons:lowerboundProp-SDP-2} become 
\begin{align}
\exists \boldsymbol{\lambda}_k \geq 0: \ & s + \boldsymbol{\xi}_1^{\top} \mathbf{Q}_{\textup{r}}^{\prime} \boldsymbol{\xi}_1 + \boldsymbol{\xi}_2^{\top} \mathbf{Q}_{\textup{r}}^{\prime \prime} \boldsymbol{\xi}_2  + \mathbf{q}^{\top} \left( \mathbf{B}_1 \boldsymbol{\xi}_1 + \mathbf{B}_2 \boldsymbol{\xi}_2 + \bar{\mathbf{B}} \boldsymbol{\xi}_3 \right) - y^0_k (\mathbf{x}) \nonumber \\
& - y_k(\mathbf{x})^{\top} \left(\mathbf{U}\boldsymbol{\Lambda}^{ \frac{1}{2} }\left(\mathbf{B}_1 \boldsymbol{\xi}_1 + \mathbf{B}_2 \boldsymbol{\xi}_2 + \bar{\mathbf{B}} \boldsymbol{\xi}_3 \right) +\boldsymbol{\mu} \right) \nonumber  \\
& + \boldsymbol{\lambda}_k^{\top} \left( \mathbf{A}\left(\mathbf{U}\boldsymbol{\Lambda}^{ \frac{1}{2} }\left(\mathbf{B}_1 \boldsymbol{\xi}_1 + \mathbf{B}_2 \boldsymbol{\xi}_2 + \bar{\mathbf{B}} \boldsymbol{\xi}_3 \right)+\boldsymbol{\mu} \right) - \mathbf{b} \right) \geq 0, \nonumber \\
&\hspace{2cm} \forall \boldsymbol{\xi}_1 \in \mathbb{R}^{m_1}, \boldsymbol{\xi}_2 \in \mathbb{R}^{K-m_1}, \boldsymbol{\xi}_3 \in \mathbb{R}^{m-K}, \ \forall k \in [K]. \label{Cons:lowerboundProp-SDP-2.5}
\end{align}

\noindent
We further define 
\begin{align*}
\mathbf{Z}_k = 
\begin{bmatrix}  
s - y_k^0(\mathbf{x}) - \boldsymbol{\lambda}_k^{\top}\mathbf{b} - y_k(\mathbf{x})^{\top} \boldsymbol{\mu} + \boldsymbol{\lambda}_k^{\top}\mathbf{A}\boldsymbol{\mu} &   \frac{1}{2}    (\mathbf{h}_k^{\prime})^{\top} & \frac{1}{2}(\mathbf{h}_k^{\prime \prime})^{\top}    \\ 
\frac{1}{2} \mathbf{h}_k^{\prime} & \mathbf{Q}_{ \textup{r} }^{\prime} & \mathbf{0}_{m_1 \times (K-m_1)} \\ 
\frac{1}{2} \mathbf{h}_k^{\prime \prime} & \mathbf{0}_{(K-m_1) \times m_1} & \mathbf{Q}_{ \textup{r} }^{\prime\prime} 
\end{bmatrix}, \ \forall k \in [K],
\end{align*}
where $\mathbf{h}_k^{\prime} = \mathbf{B}_1^{\top} \mathbf{q} + (\mathbf{U} \boldsymbol{\Lambda}^{{\frac{1}{2}}} \mathbf{B}_1 )^{\top} (\mathbf{A}^{\top} \boldsymbol{\lambda}_k - y_k(\mathbf{x}) )$ and $\mathbf{h}_k^{\prime \prime} = \mathbf{B}_2^{\top}\mathbf{q} +(\mathbf{U}\boldsymbol{\Lambda}^{{\frac{1}{2}}} \mathbf{B}_2 )^{\top} (\mathbf{A}^{\top}\boldsymbol{\lambda}_k-y_k(\mathbf{x}) )$ for any $k \in [K]$.
It follows that
\begin{align}
\eqref{Cons:lowerboundProp-SDP-2.5} & \iff \exists \boldsymbol{\lambda}_k \geq 0: \ \left( 1,\boldsymbol{\xi}_1^{\top}, \boldsymbol{\xi}_2^{\top} \right) \mathbf{Z}_k \left( 1,\boldsymbol{\xi}_1^{\top}, \boldsymbol{\xi}_2^{\top} \right)^{\top} +  \boldsymbol{\xi}_3^{\top} \left( \bar{\mathbf{B}}^{\top} \mathbf{q} + \left(\mathbf{U}\boldsymbol{\Lambda}^{{\frac{1}{2}}} \bar{\mathbf{B}} \right)^{\top} \left(\mathbf{A}^{\top}\boldsymbol{\lambda}_k-y_k(\mathbf{x})\right) \right) \geq 0, \nonumber \\
& \hspace{8cm} \forall \boldsymbol{\xi}_1 \in \mathbb{R}^{m_1}, \boldsymbol{\xi}_2 \in \mathbb{R}^{K-m_1}, \boldsymbol{\xi}_3 \in \mathbb{R}^{m-K}, \ \forall k \in [K]. \nonumber \\ 
& \iff \exists \boldsymbol{\lambda}_k \geq 0: \ \left( 1,\boldsymbol{\xi}_1^{\top},\boldsymbol{\xi}_2^{\top} \right) \mathbf{Z}_k \left( 1,\boldsymbol{\xi}_1^{\top},\boldsymbol{\xi}_2^{\top} \right)^{\top} \geq 0, \ \forall \boldsymbol{\xi}_1 \in \mathbb{R}^{m_1}, \boldsymbol{\xi}_2 \in \mathbb{R}^{K-m_1}, \ \forall k \in [K]; \label{Cons:lowerboundProp-SDP-3} \\
&  \hspace{2.5cm} \bar{\mathbf{B}}^{\top} \mathbf{q} +\left(\mathbf{U}\boldsymbol{\Lambda}^{{\frac{1}{2}}} \bar{\mathbf{B}} \right)^{\top}\left(\mathbf{A}^{\top}\boldsymbol{\lambda}_k-y_k(\mathbf{x})\right) = 0,\ \forall k \in [K]. \nonumber \\
&\iff  \exists \boldsymbol{\lambda}_k \geq 0: \  \mathbf{Z}_k \succeq 0,\  \bar{\mathbf{B}}^{\top} \mathbf{q} +\left(\mathbf{U}\boldsymbol{\Lambda}^{{\frac{1}{2}}} \bar{\mathbf{B}} \right)^{\top}\left(\mathbf{A}^{\top}\boldsymbol{\lambda}_k-y_k(\mathbf{x})\right) = 0,\ \forall k \in [K]. \nonumber \\ 
& \iff  \exists \boldsymbol{\lambda}_k \geq 0: \  \mathbf{Z}_k \succeq 0,\  \bar{\mathbf{B}}^{\top} \left(\mathbf{q}+\left(\mathbf{U}\boldsymbol{\Lambda}^{{\frac{1}{2}}}  \right)^{\top}\left(\mathbf{A}^{\top}\boldsymbol{\lambda}_k-y_k(\mathbf{x})\right) \right) = 0,\ \forall k \in [K]. \label{Cons:lowerboundProp-SDP-4} \\ 
& \iff  \exists \boldsymbol{\lambda}_k \geq 0, \mathbf{u}_k^{\prime} \in \mathbb{R}^{m_1}, \mathbf{u}_k^{\prime \prime} \in \mathbb{R}^{K-m_1}: \nonumber\\
& \hspace{2cm} \mathbf{Z}_k \succeq 0, \ \mathbf{q}+\left(\mathbf{U}\boldsymbol{\Lambda}^{{\frac{1}{2}}}  \right)^{\top}\left(\mathbf{A}^{\top}\boldsymbol{\lambda}_k-y_k(\mathbf{x})\right) = \mathbf{B}_1 \mathbf{u}_k^{\prime} + \mathbf{B}_2 \mathbf{u}_k^{\prime \prime},\ \forall k \in [K]. \label{Cons:lowerboundProp-sdp-end}
\end{align}
The first equivalence holds due to the definition of $\mathbf{Z}_k$.
For the third equivalence, clearly $\Longleftarrow$ follows from the definition of a PSD matrix.
To prove $\Longrightarrow$, we consider two possible cases for any 
$(\eta_0  \in \mathbb{R}, \boldsymbol{\eta}_1^{\top} \in \mathbb{R}^{m_1}, \boldsymbol{\eta}_2^{\top} \in \mathbb{R}^{K-m_1} )^{\top} \in \mathbb{R}^{K+1}$:
(i) if $\eta_0=0$, then $( \eta_0, \boldsymbol{\eta}_1^{\top}, \boldsymbol{\eta}_2^{\top} ) \mathbf{Z}_k ( \eta_0, \boldsymbol{\eta}_1^{\top}, \boldsymbol{\eta}_2^{\top} )^{\top} 
= \boldsymbol{\eta}_1^{\top} \mathbf{Q}_{\textup{r}}^{\prime} \boldsymbol{\eta}_1 + \boldsymbol{\eta}_2^{\top} \mathbf{Q}_{\textup{r}}^{\prime \prime} \boldsymbol{\eta}_2 \geq 0$ 
because $\mathbf{Q}_{\textup{r}}^{\prime}$ and $\mathbf{Q}_{\textup{r}}^{\prime \prime}$ are PSD;
(ii) if $\eta_0 \neq 0$, then we have $( \eta_0, \boldsymbol{\eta}_1^{\top}, \boldsymbol{\eta}_2^{\top} ) \mathbf{Z}_k ( \eta_0, \boldsymbol{\eta}_1^{\top}, \boldsymbol{\eta}_2^{\top} )^{\top} 
= \eta_0^2 ( 1, \frac{\boldsymbol{\eta}_1^{\top}}{\eta_0}, \frac{\boldsymbol{\eta}_2^{\top}}{\eta_0} ) \mathbf{Z}_k ( 1, \frac{\boldsymbol{\eta}_1^{\top}}{\eta_0}, \frac{\boldsymbol{\eta}_2^{\top}}{\eta_0} )^{\top} \geq 0$ according to \eqref{Cons:lowerboundProp-SDP-3}. Therefore, $\Longrightarrow$ holds.
For the fifth equivalence, \eqref{Cons:lowerboundProp-SDP-4} shows that $\mathbf{q} + (\mathbf{U} \boldsymbol{\Lambda}^{{\frac{1}{2}}}  )^{\top} (\mathbf{A}^{\top} \boldsymbol{\lambda}_k - y_k(\mathbf{x}))$ is in the null space of $\bar{\mathbf{B}}$ and thus cannot be represented by basis vectors in the space of $\bar{\mathbf{B}}$.
Because $[\mathbf{B},\bar{\mathbf{B}} ]$ span the space of $\mathbb{R}^m$, we have $\mathbf{q} + (\mathbf{U} \boldsymbol{\Lambda}^{{\frac{1}{2}}}  )^{\top} (\mathbf{A}^{\top} \boldsymbol{\lambda}_k - y_k(\mathbf{x}))$ should be in the space of $ \mathbf{B}$.
That is, there exists $\mathbf{u}_k^{\prime} \in \mathbb{R}^{m_1} $ and $\mathbf{u}_k^{\prime \prime} \in \mathbb{R}^{K-m_1} $ such that $\mathbf{q}+(\mathbf{U}\boldsymbol{\Lambda}^{{\frac{1}{2}}}  )^{\top}(\mathbf{A}^{\top}\boldsymbol{\lambda}_k-y_k(\mathbf{x})) = \mathbf{B}_1 \mathbf{u}_k^{\prime} + \mathbf{B}_2 \mathbf{u}_k^{\prime \prime}$ for any $k \in [K]$. Meanwhile, because $\mathbf{B}^{\top} \mathbf{B} = \mathbf{I}_{K}$, we have 
\begin{align}
\mathbf{h}_k^{\prime} = \mathbf{B}_1^{\top}\mathbf{q} +\left(\mathbf{U}\boldsymbol{\Lambda}^{{\frac{1}{2}}} \mathbf{B}_1 \right)^{\top}\left(\mathbf{A}^{\top}\boldsymbol{\lambda}_k-y_k(\mathbf{x})\right) = \mathbf{B}_1^{\top} \mathbf{B}_1\mathbf{u}_k^{\prime}=\mathbf{u}_k^{\prime}, \ \forall k \in [K], \nonumber \\ 
\mathbf{h}_k^{\prime \prime} = \mathbf{B}_2^{\top}\mathbf{q} + \left( \mathbf{U}\boldsymbol{\Lambda}^{{\frac{1}{2}}} \mathbf{B}_2  \right)^{\top} \left( \mathbf{A}^{\top} \boldsymbol{\lambda}_k - y_k (\mathbf{x}) \right) = \mathbf{B}_2^{\top} \mathbf{B}_2\mathbf{u}_k^{\prime \prime}=\mathbf{u}_k^{\prime \prime}, \ \forall k \in [K]. \nonumber
\end{align} 
By replacing constraints \eqref{Cons:lowerboundRef1-1} with \eqref{Cons:lowerboundProp-sdp-end}, we obtain the following problem:
\begin{subeqnarray}
& \min\limits_{ \substack{\mathbf{x}, s, {\hat{\boldsymbol{\lambda}}}, \mathbf{q},  \\ \mathbf{Q}_{ \textup{r} }^{\prime}, \mathbf{Q}_{ \textup{r} }^{\prime\prime}, \hat{\mathbf{u}}^{\prime}, \hat{\mathbf{u}}^{\prime \prime}, \\ \mathbf{B}_1, \mathbf{B}_2 } } & s + \gamma_2 \mathbf{I}_{m_1} \bullet \mathbf{Q}_{ \textup{r} }^{\prime}  + \sqrt{\gamma_1}\left \| \mathbf{q} \right \|_2 \slabel{eqn:obj:newlb-proof} \\
&\textnormal{s.t.} & \small{  \begin{bmatrix}  
s-y_k^0(\mathbf{x})-\boldsymbol{\lambda}_k^{\top}\mathbf{b}-y_k(\mathbf{x})^{\top}\boldsymbol{\mu} +\boldsymbol{\lambda}_k^{\top}\mathbf{A}\boldsymbol{\mu} &   \frac{1}{2}    (\mathbf{u}_k^{\prime})^{\top} & \frac{1}{2}(\mathbf{u}_k^{\prime \prime})^{\top}    \\ 
\frac{1}{2}   \mathbf{u}_k^{\prime} & \mathbf{Q}_{ \textup{r} }^{\prime} & \mathbf{0}_{m_1 \times (K-m_1)} \\ 
\frac{1}{2}\mathbf{u}_k^{\prime \prime} & \mathbf{0}_{(K-m_1)\times m_1}  & \mathbf{Q}_{ \textup{r} }^{\prime\prime} 
\end{bmatrix}} \succeq 0, \ \forall k\in[K], \slabel{Cons:MotivateUpperbound-1}  \\
&& \mathbf{q} +  \left(\mathbf{U}\boldsymbol{\Lambda}^{{\frac{1}{2}}}  \right)^{\top}\left(\mathbf{A}^{\top}\boldsymbol{\lambda}_k-y_k(\mathbf{x})\right) = \mathbf{B}_1 \mathbf{u}_k^{\prime} +  \mathbf{B}_2 \mathbf{u}_k^{\prime \prime} ,\ \forall k \in [K] ,   \\
&& \mathbf{x} \in \mathcal{X},\ \mathbf{q}  \in \mathbb{R}^{m}, \ \mathbf{Q}_{ \textup{r} }^{\prime} \in \mathbb{R}^{ m_1 \times m_1 }, \ \mathbf{Q}_{ \textup{r} }^{\prime \prime} \in \mathbb{R}^{ (K-m_1) \times (K-m_1) },  \\ 
&& \mathbf{B}_1 \in \mathbb{R}^{m \times m_1},  \mathbf{B}_2  \in \mathbb{R}^{m \times (K-m_1)}, \ [  \mathbf{B}_1,  \mathbf{B}_2 ]^{\top} [\mathbf{B}_1,  \mathbf{B}_2] = \mathbf{I}_K, \\
&& \hat{\boldsymbol{\lambda}} = \left\{\boldsymbol{\lambda}_1, \dots, \boldsymbol{\lambda}_K \right\}, \ \boldsymbol{\lambda}_k \in \mathbb{R}_+^{l}, \ \forall k \in [K], \\
&&\hat{\mathbf{u}}^{\prime} = \left\{\mathbf{u}_1^{\prime}, \dots, \mathbf{u}_K^{\prime} \right\}, \ \mathbf{u}_k^{\prime} \in \mathbb{R}^{m_1}, \ \forall k\in[K],  \\ 
&&\hat{\mathbf{u}}^{\prime \prime} = \left\{\mathbf{u}_1^{\prime \prime}, \dots, \mathbf{u}_K^{\prime \prime} \right\}, \ \mathbf{u}_k^{\prime \prime} \in \mathbb{R}^{K-m_1}, \ \forall k \in [K].
\end{subeqnarray}
Note that the value of $\mathbf{Q}_{ \textup{r} }^{\prime \prime}$ does not contribute to the objective function \eqref{eqn:obj:newlb-proof}.
We can then let $M$ be an arbitrarily large positive number and $\mathbf{Q}_{ \textup{r} }^{\prime\prime} = M \mathbf{I}_{(K-m_1)\times (K-m_1)}$ be an optimal solution, by which constraints \eqref{Cons:MotivateUpperbound-1} become
\begin{align}
\begin{bmatrix}  
s-y_k^0(\mathbf{x})-\boldsymbol{\lambda}_k^{\top}\mathbf{b}-y_k(\mathbf{x})^{\top}\boldsymbol{\mu} +\boldsymbol{\lambda}_k^{\top}\mathbf{A}\boldsymbol{\mu} &   \frac{1}{2}    (\mathbf{u}_k^{\prime})^{\top}     \\ 
\frac{1}{2}   \mathbf{u}_k^{\prime} & \mathbf{Q}_{ \textup{r} }^{\prime}  
\end{bmatrix}\succeq 0, \ \forall k \in [K]. \label{Cons:MotivateUpperbound-1-new}  
\end{align}
By replacing \eqref{Cons:MotivateUpperbound-1} with \eqref{Cons:MotivateUpperbound-1-new}, we obtain the formulation of Problem \eqref{Equ:.lowerbound-sdpform-rewrite}.

Based on the formulation of Problem \eqref{Equ:.lowerbound-sdpform-rewrite}, now we show that the three conclusions hold. 
Note that for any $\mathbf{B} \in \mathcal{B}_K = \{\mathbf{B}\in \mathbb{R}^{m \times K} \ | \ \mathbf{B}^{\top} \mathbf{B} = \mathbf{I}_{K} \}$, the optimal value of Problem \eqref{Equ:upperbound-setform}, i.e.,  $\Theta_{\textup{U}}(K)$, reaches the optimal value of the original Problem \eqref{Equ:MainProblem}, i.e., $ \Theta_{\textup{M}}(m)$.
We would like to show that by relaxing the constraints in Problem \eqref{Equ:upperbound-setform}, we can obtain the exact formulation of Problem \eqref{Equ:.lowerbound-sdpform-rewrite}, thereby the three conclusions hold.

First, we rewrite constraints \eqref{Cons:upperbound-sdpform-1}--\eqref{Cons:upperbound-sdpform-2} in Problem \eqref{Equ:upperbound-setform} with $m_1 = K$  by dividing $\mathbf{B}$ into $[\mathbf{B}_1, \mathbf{B}_2]$ and $\mathbf{u}_k$ into $((\mathbf{u}_k^{\prime})^{\top},(\mathbf{u}_k^{\prime \prime})^{\top} )^{\top}$. 
Thus, we obtain the following formulation:
\begin{subeqnarray}  \label{Equ:upperbound-sdpform-rewrite}
& \min\limits_{ \substack{\mathbf{x}, s, {\hat{\boldsymbol{\lambda}}}, \\ \mathbf{q}, \mathbf{Q}_{ \textup{r} }, \hat{\mathbf{u}}^{\prime}, \hat{\mathbf{u}}^{\prime \prime},\\ \mathbf{B}_1, \mathbf{B}_2 } } & s + \gamma_2 \mathbf{I}_{K} \bullet \mathbf{Q}_{ \textup{r} }  + \sqrt{\gamma_1}\left \| \mathbf{q} \right \|_2 \slabel{Obj:upperbound-sdpform-rewrite}\\
&\textnormal{s.t.} & \small{  \begin{bmatrix} s-y_k^0(\mathbf{x})-\boldsymbol{\lambda}_k^{\top}\mathbf{b}-y_k(\mathbf{x})^{\top}\boldsymbol{\mu} +\boldsymbol{\lambda}_k^{\top}\mathbf{A}\boldsymbol{\mu} & \hspace{0.1 in} \frac{1}{2} \left((\mathbf{u}_k^{\prime})^{\top},(\mathbf{u}_k^{\prime \prime})^{\top} \right)  \\ 
\frac{1}{2} \left((\mathbf{u}_k^{\prime})^{\top},(\mathbf{u}_k^{\prime \prime})^{\top} \right)^{\top} & \mathbf{Q}_{ \textup{r} } \end{bmatrix}} \succeq 0, \ \forall k\in[K], \slabel{Cons:upperbound-sdpform-rewrite-1} \\
&& \mathbf{q} +  \left(\mathbf{U}\boldsymbol{\Lambda}^{{\frac{1}{2}}}  \right)^{\top}\left(\mathbf{A}^{\top}\boldsymbol{\lambda}_k-y_k(\mathbf{x})\right) = \mathbf{B}_1 \mathbf{u}_k^{\prime} +  \mathbf{B}_2 \mathbf{u}_k^{\prime \prime} ,\ \forall k \in [K] , \slabel{Cons:upperbound-sdpform-rewrite-2} \\
&& \mathbf{x} \in \mathcal{X}, \ [
  \mathbf{B}_1,  \mathbf{B}_2 
]^{\top} [\mathbf{B}_1,  \mathbf{B}_2] = \mathbf{I}_K, \slabel{Cons:upperbound-sdpform-rewrite-3}\\
&& \mathbf{q}  \in \mathbb{R}^{m}, \ \mathbf{Q}_{ \textup{r} } \in \mathbb{R}^{ K \times K }, \ \mathbf{B}_1 \in \mathbb{R}^{m\times m_1},  \mathbf{B}_2  \in \mathbb{R}^{m\times (K-m_1)},  \\ 
&& \hat{\boldsymbol{\lambda}} = \left\{\boldsymbol{\lambda}_1, \dots, \boldsymbol{\lambda}_K \right\}, \ \boldsymbol{\lambda}_k \in \mathbb{R}_+^{l}, \ \forall k \in [K], \\
&&\hat{\mathbf{u}}^{\prime} = \left\{\mathbf{u}_1^{\prime}, \dots, \mathbf{u}_K^{\prime} \right\}, \ \mathbf{u}_k^{\prime} \in \mathbb{R}^{m_1}, \ \forall k\in[K],  \\ 
&&\hat{\mathbf{u}}^{\prime \prime} = \left\{\mathbf{u}_1^{\prime \prime}, \dots, \mathbf{u}_K^{\prime \prime} \right\}, \ \mathbf{u}_k^{\prime \prime} \in \mathbb{R}^{K-m_1}, \ \forall k \in [K]. \slabel{Cons:upperbound-sdpform-rewrite-end}
\end{subeqnarray}

Second, we relax constraints \eqref{Cons:upperbound-sdpform-rewrite-1} into 
\begin{align}
\begin{bmatrix} s-y_k^0(\mathbf{x})-\boldsymbol{\lambda}_k^{\top}\mathbf{b}-y_k(\mathbf{x})^{\top}\boldsymbol{\mu} +\boldsymbol{\lambda}_k^{\top}\mathbf{A}\boldsymbol{\mu} & \hspace{0.1 in} \frac{1}{2}   (\mathbf{u}_k^{\prime})^{\top}  \\ 
\frac{1}{2}  \mathbf{u}_k^{\prime}  & \mathbf{Q}_{ \textup{r} }^{\prime} \end{bmatrix} \succeq 0, \ \forall k\in[K], \label{eqn:sdp-relax-upper-left}
\end{align}
where $\mathbf{Q}_{ \textup{r} }^{\prime} \in \mathbb{R}^{m_1 \times m_1}$ is the upper-left submatrix of $\mathbf{Q}_{ \textup{r} }$. 
Note that if we use \eqref{eqn:sdp-relax-upper-left} to replace \eqref{Cons:upperbound-sdpform-rewrite-1}, we obtain a relaxation and accordingly lower bound for Problem \eqref{Equ:upperbound-sdpform-rewrite}.
In addition, we further reduce the optimal value of the relaxation by replacing $\mathbf{Q}_{ \textup{r} }$ in the objective function \eqref{Obj:upperbound-sdpform-rewrite} with $\mathbf{Q}_{ \textup{r} }^{\prime}$.
That is, we obtain a lower bound for the optimal value of Problem \eqref{Equ:upperbound-setform} with $m_1 = K$ (i.e., Problem \eqref{Equ:MainProblem}).
After these two steps of relaxations, we obtain the exact formulation of Problem \eqref{Equ:.lowerbound-sdpform-rewrite}.
Thus, we can conclude that Problem \eqref{Equ:.lowerbound-sdpform-rewrite} is a relaxation of Problem \eqref{Equ:upperbound-setform} with $m_1=K$. Therefore, by the conclusion in Theorem \ref{Theo:upperboundSameoptimal}, we have
\begin{align*}
\Theta_{\textup{L2}}(m_1) \leq \Theta_{\textup{U}} (K) = \Theta_{\textup{M}}(m).
\end{align*}
That is, the conclusion (i) holds.

For the conclusion (ii): For any $0 \leq m_1 < m_2 \leq K$, we can follow the above two steps of relaxations to relax Problem \eqref{Equ:upperbound-sdpform-rewrite} to the problem with the optimal value $\Theta_{\textup{L2}}(m_2)$, and based on this relaxed problem, we can further relax it to the problem with the optimal value $\Theta_{\textup{L2}}(m_1)$. 
Because all these problems are minimization problems, we have $\Theta_{\textup{L2}}(m_1) \leq \Theta_{\textup{L2}}(m_2)$.

For the conclusion (iii):  When $m_1 = K$, Problem \eqref{Equ:.lowerbound-sdpform-rewrite} becomes Problem \eqref{Equ:upperbound-setform} with $m_1=K$. Thus, by the conclusion in Theorem \ref{Theo:upperboundSameoptimal}, we have 
\begin{align*}
\Theta_{\textup{L2}}(K) = \Theta_{\textup{U}} (K) = \Theta_{\textup{M}}(m). \Halmos 
\end{align*}

\section{Supplement to Section \ref{Sec:Algorithm}}

\subsection{Proof of Proposition \ref{Prop:algorithm}}
Given $ (\mathbf{x}, s, \mathbf{q}, \mathbf{Q}_{ \textup{r} }, \boldsymbol{\lambda}_k, 
\tilde{\mathbf{u}}_k, \mathbf{u}_k, \boldsymbol{\beta}_k, \forall k \in [K])$, 
we can omit the constant in the objective function of Problem \eqref{Equ:ALProblem} and rewrite 
this problem as follows:
\begin{align}
\min\limits_{\mathbf{B} \in  \mathbb{R}^{m\times m_1} } \left\{ 
\sum_{k=1}^K   - \boldsymbol{\beta}_k^{\top}  \mathbf{B} \mathbf{u}_k   + \sum_{k=1}^K \left( - \rho  \tilde{\mathbf{u}}_k^{\top} \mathbf{B} \mathbf{u}_k + \frac{\rho}{2} \mathbf{u}_k^{\top} \mathbf{B}^{\top}  \mathbf{B} \mathbf{u}_k \right) 
\ \middle| \
\mathbf{B}^{\top} \mathbf{B} = \mathbf{I}_{m_1}
\right\}. \label{Equ:proofAl}
\end{align}
By $\mathbf{B}^{\top} \mathbf{B} = \mathbf{I}_{m_1}$, the term $ \mathbf{u}_k^{\top} \mathbf{B}^{\top}  \mathbf{B} \mathbf{u}_k$ is also a constant. 
Because $\boldsymbol{\beta}_k^{\top}  \mathbf{B} \mathbf{u}_k = (\boldsymbol{\beta}_k \mathbf{u}_k^{\top}) \bullet \mathbf{B}$ 
and $\tilde{\mathbf{u}}_k^{\top} \mathbf{B} \mathbf{u}_k = (\tilde{\mathbf{u}}_k  \mathbf{u}_k^{\top}) \bullet \mathbf{B}$, 
we can further rewrite Problem \eqref{Equ:proofAl} as follows:
\begin{align*}
\max\limits_{\mathbf{B} \in  \mathbb{R}^{m\times m_1} } \left\{ 
\sum_{k=1}^K \left( \boldsymbol{\beta}_k \mathbf{u}_k^{\top} + \rho \tilde{\mathbf{u}}_k \mathbf{u}_k^{\top} \right)  \bullet \mathbf{B} 
\ \middle| \
\mathbf{B}^{\top} \mathbf{B} = \mathbf{I}_{m_1}
\right\}.    
\end{align*}

\noindent
By the SVD, i.e., $\sum_{k=1}^K ( \boldsymbol{\beta}_k \mathbf{u}_k^{\top} + \rho \tilde{\mathbf{u}}_k 
\mathbf{u}_k^{\top} ) = \tilde{\mathbf{U}}  \tilde{\boldsymbol{\Sigma}}  \tilde{\mathbf{V}}^{\top}$, 
we have 
\begin{align*}
\mathbf{B}^* = & \argmax\limits_{ \mathbf{B}^{\top} \mathbf{B} = \mathbf{I}_{m_1} } \ \left( \tilde{\mathbf{U}}  \tilde{\boldsymbol{\Sigma}}  \tilde{\mathbf{V}}^{\top} \right) \bullet \mathbf{B} 
= \argmax\limits_{ \mathbf{B}^{\top} \mathbf{B} = \mathbf{I}_{m_1} } \ tr \left( \tilde{\mathbf{U}}  \tilde{\boldsymbol{\Sigma}}  \tilde{\mathbf{V}}^{\top} \mathbf{B}^{\top} \right)  \\
= & \argmax\limits_{ \mathbf{B}^{\top} \mathbf{B} = \mathbf{I}_{m_1} } \ tr \left(   \tilde{\boldsymbol{\Sigma}}  \tilde{\mathbf{V}}^{\top} \mathbf{B}^{\top} \tilde{\mathbf{U}} \right)  
= \argmax\limits_{ \mathbf{B}^{\top} \mathbf{B} = \mathbf{I}_{m_1} } \  \tilde{\boldsymbol{\Sigma}} \bullet \left( \tilde{\mathbf{U}}^{\top} \mathbf{B}  \tilde{\mathbf{V}} \right),
\end{align*}
where the second and fourth equalities hold by the definition of a matrix's trace and the third equality holds by the cyclic property of a matrix's trace.
\cite{elden1999procrustes} show that $\mathbf{B}^* = \tilde{\mathbf{U}}  \tilde{\mathbf{V}}^{\top}$ is an optimal solution.
\Halmos

\subsection{Proof of Proposition \ref{Prop:upperbound}}
First, we have 
$\Theta_{\textup{M}}(m) \geq 
\Theta_{\textup{L}} (m_1) \geq 
\underline{\Theta} (m_1, \mathbf{B}^{\prime}) = s^*  + \gamma_2 \mathbf{I}_{m_1} \bullet \mathbf{Q}^*_{\textup{r}} + \sqrt{\gamma_1} \| \mathbf{q}^*_{\textup{r}} \|_2$, 
where the first inequality holds by conclusion (i) of Theorem \ref{Theo:lowerbound} and 
the second inequality holds because $ \mathbf{B}^{\prime} $ is a feasible solution of Problem \eqref{Equ:appro} and this problem is a maximization problem. 

Next, we would like to construct a feasible solution $(\mathbf{x}^{\prime}, s^{\prime}, \hat{\boldsymbol{\lambda}}^{\prime},\mathbf{q}^{\prime},\mathbf{Q}^{\prime})$ of Problem \eqref{Equ:MainProblem}. 
We set $\mathbf{x}^{\prime}=\mathbf{x}^*$,
$\hat{\boldsymbol{\lambda}}^{\prime} = \hat{\boldsymbol{\lambda}}^*$, 
$s^{\prime} = s^*+s_0$, 
$\mathbf{q}^{\prime} = \mathbf{B}^{\prime} \mathbf{q}^*_{\textup{r}}$, 
and $\mathbf{Q}^{\prime} = \mathbf{B}^{\prime} \mathbf{Q}^*_{\textup{r}} (\mathbf{B}^{\prime})^{\top} + \mathbf{Q}_0$, where $s_0 \geq 0$ and $\mathbf{Q}_0 \succeq 0$ and their values will be decided later. 
Clearly, this solution satisfies constraints \eqref{Cons:MainProblemcons2}. 
For this solution to satisfy constraints \eqref{Cons:SDP}, the values $s_0$ and $ \mathbf{Q}_0$ should satisfy 
\begin{align}
& \left(S_k + s_0 \right) \left( \mathbf{B}^{\prime} \mathbf{Q}^*_{\textup{r}} (\mathbf{B}^{\prime})^{\top} + \mathbf{Q}_0 \right) \succeq  \frac{1}{4} \left(\mathbf{B}^{\prime} \mathbf{q}^*_{\textup{r}} +\left(\mathbf{U}\boldsymbol{\Lambda}^{{\frac{1}{2}}}\right)^{\top}\left(\mathbf{A}^{\top}\boldsymbol{\lambda}^*_k-y_k(\mathbf{x}^*)\right)\right)  \nonumber \\
& \hspace{5cm} \times
\left( \mathbf{B}^{\prime} \mathbf{q}^*_{\textup{r}} +\left(\mathbf{U}\boldsymbol{\Lambda}^{{\frac{1}{2}}}\right)^{\top}\left(\mathbf{A}^{\top}\boldsymbol{\lambda}^*_k-y_k(\mathbf{x}^*)\right) \right)^{\top} = \frac{1}{4}  \mathbf{M}_k, \ \forall k \in [K]. \label{Equ:UB1-1-inproof}
\end{align} 

Note that, if $ ( S + s_0 ) \mathbf{Q}_0 \succeq  (1/4)  \mathbf{M}_k$ for any $k \in [K]$, then \eqref{Equ:UB1-1-inproof} holds.
By the definition of $\mathbf{M}_k$, we have $\mathbf{M}_k \succeq 0$ for any $ k \in [K]$. Therefore, for any $s_0 \geq 0$, we can construct 
\begin{equation}  
\mathbf{Q}_0 = \sum_{k=1}^K \frac{1}{4 (S+s_0)}   \mathbf{M}_k   \nonumber
\end{equation}
such that \eqref{Equ:UB1-1-inproof} holds and hence $(\mathbf{x}^{\prime}, s^{\prime}, \hat{\boldsymbol{\lambda}}^{\prime},\mathbf{q}^{\prime},\mathbf{Q}^{\prime})$ is a feasible solution of Problem \eqref{Equ:MainProblem}. 
The objective value (denoted by $\Theta_{\textup{M}}^{\prime}$) with respect to this constructed solution is
\begin{align*}
s^{\prime} + \gamma_2\mathbf{I}_{m} \bullet \mathbf{Q}^{\prime} + \sqrt{\gamma_1} \left\| \mathbf{q}^{\prime} \right\|_2 
&= s^*+s_0 +\gamma_2\mathbf{I}_{m} \bullet \mathbf{B}^{\prime} \mathbf{Q}^*_{\textup{r}} (\mathbf{B}^{\prime})^{\top} + \gamma_2\mathbf{I}_{m} \bullet \mathbf{Q}_0 + \sqrt{\gamma_1} \left\| \mathbf{B}^{\prime} \mathbf{q}^*_{\textup{r}} \right\|_2 \\ 
&= s^*+s_0 +\gamma_2\mathbf{I}_{m_1} \bullet  \mathbf{Q}^*_{\textup{r}} + \gamma_2\mathbf{I}_{m} \bullet \mathbf{Q}_0 + \sqrt{\gamma_1} \left\|  \mathbf{q}^*_{\textup{r}} \right\|_2\\
&= \underline{\Theta} (m_1,\mathbf{B}^{\prime}) +s_0+ \sum_{k=1}^K \frac{\gamma_2}{4  (S+s_0) } \mathbf{I}_m \bullet  \mathbf{M}_k,
\end{align*}
where the second equality holds because $\mathbf{I}_{m} \bullet \mathbf{B}^{\prime} \mathbf{Q}^*_{\textup{r}} (\mathbf{B}^{\prime})^{\top} = \mathbf{I}_{m_1} \bullet  \mathbf{Q}^*_{\textup{r}} (\mathbf{B}^{\prime})^{\top} \mathbf{B}^{\prime} = \mathbf{I}_{m_1} \bullet  \mathbf{Q}^*_{\textup{r}}$ and $(\mathbf{q}^{*}_{\textup{r}})^{\top} (\mathbf{B}^{\prime})^{\top} \mathbf{B}^{\prime} \mathbf{q}^*_{\textup{r}} = (\mathbf{q}^{*}_{\textup{r}})^{\top}  \mathbf{q}^*_{\textup{r}}$. 
As this constructed solution is a feasible solution of Problem \eqref{Equ:MainProblem}, which is a minimization problem, we have $\Theta_{\textup{M}}(m) \leq \Theta_{\textup{M}}^{\prime}$
It follows that 
\begin{align}
& \Theta_{\textup{M}}(m) - \Theta_{\textup{L}} (m_1)
\leq \Theta_{\textup{M}}^{\prime} - \underline{\Theta} (m_1,\mathbf{B}^{\prime}) = s_0+ \sum_{k=1}^K \frac{\gamma_2}{4  (S+s_0) } \mathbf{I}_m \bullet  \mathbf{M}_k. \label{eqn:gap-upper-bound-1}
\end{align}

We further choose a value of $s_0$ to minimize the the right-hand side (RHS) of \eqref{eqn:gap-upper-bound-1}.
Note that (i) If $\sqrt{P} - S < 0$, then the RHS of \eqref{eqn:gap-upper-bound-1} is minimized at $P/S$ with $s_0 = 0$;
(ii) If $\sqrt{P} - S \geq 0$, then the RHS of \eqref{eqn:gap-upper-bound-1} is minimized at $2\sqrt{P} - S$ with $s_0 = \sqrt{P} - S$.
Therefore, we conclude that the proposition holds.
\Halmos 

\section{Supplement to Section \ref{Sec:experiments}}

\subsection{Multiproduct Newsvendor Problem} \label{subsec-apx:newsvendor}
By Proposition \ref{Prop:cheramin2022}, Problem \eqref{Equ:Newsvendor} has the same optimal value as the following SDP formulation: 
\begin{subeqnarray} \label{Equ:NewsvendorSDP}
&\min\limits_{ \substack{ \mathbf{x},s,\boldsymbol{\lambda}_1, \\  \boldsymbol{\lambda}_2,\mathbf{q},\mathbf{Q} } } & s +\gamma_2\mathbf{I}_{m} \bullet \mathbf{Q} + \sqrt{\gamma_1} \left \| \mathbf{q}\right \|_2\\ 
&\text{s.t.} &\small{ \begin{bmatrix}
s-(\mathbf{c}-\mathbf{v})^{\top} \mathbf{x}-\boldsymbol{\lambda}_1^{\top}(\mathbf{b}-\mathbf{A}\boldsymbol{\mu}) &\ \frac{1}{2}  \left(\mathbf{q} +\left(\mathbf{U}\boldsymbol{\Lambda}^{{\frac{1}{2}}}\right)^{\top}\mathbf{A}^{\top}\boldsymbol{\lambda}_1\right)^{\top} \\
\frac{1}{2}  \left(\mathbf{q} +\left(\mathbf{U}\boldsymbol{\Lambda}^{{\frac{1}{2}}}\right)^{\top}\mathbf{A}^{\top}\boldsymbol{\lambda}_1\right) & \mathbf{Q}
\end{bmatrix}} \succeq 0, \\
&& \small{\begin{bmatrix}
s-(\mathbf{c}-\mathbf{g})^{\top} \mathbf{x}-\boldsymbol{\lambda}_2^{\top}(\mathbf{b}-\mathbf{A}\boldsymbol{\mu}) + (\mathbf{v}-\mathbf{g})^{\top} \boldsymbol{\mu} &\ \frac{1}{2}  \left(\mathbf{q} +\left(\mathbf{U}\boldsymbol{\Lambda}^{{\frac{1}{2}}}\right)^{\top}\left(\mathbf{A}^{\top}\boldsymbol{\lambda}_2+\mathbf{v}-\mathbf{g}\right)\right)^{\top}\\
\frac{1}{2}  \left(\mathbf{q} +\left(\mathbf{U}\boldsymbol{\Lambda}^{{\frac{1}{2}}}\right)^{\top}\left(\mathbf{A}^{\top}\boldsymbol{\lambda}_2+\mathbf{v}-\mathbf{g}\right)\right)& \mathbf{Q} 
\end{bmatrix}} \succeq 0,  \\ 
&& \mathbf{x} \in \mathbb{R}^{m}_+, \ \boldsymbol{\lambda}_1  \in \mathbb{R}_+^{l}, \ \boldsymbol{\lambda}_2 \in \mathbb{R}_+^{l}, \ \mathbf{q} \in \mathbb{R}^m, \ \mathbf{Q} \in \mathbb{R}^{m \times m}.
\end{subeqnarray}

By the first outer approximation \eqref{Equ:bilinear}, 
the following problem provides a lower bound for the optimal value of Problem \eqref{Equ:NewsvendorSDP}:
\begin{subeqnarray} \label{Equ:NewsvendorDual}
&\max\limits_{\substack{\mathbf{B}, t_1, \mathbf{p}_1,\mathbf{P}_1,\\ t_2,\mathbf{p}_2,\mathbf{P}_2} } &   \left(t_2 \boldsymbol{\mu}^{\top}  +  \mathbf{p}_2^{\top} \left(\mathbf{U} \boldsymbol{\Lambda}^{{\frac{1}{2}}}\mathbf{B} \right)^{\top}\right) \left(\mathbf{g}-\mathbf{v}\right)  \\
&\text{s.t.} & 1- t_1 - t_2 = 0, \ \sqrt{\gamma_1} - \left\| \mathbf{p}_1 + \mathbf{p}_2 \right\|_2 \geq 0, \\
&& t_1 (\mathbf{A}\boldsymbol{\mu}-\mathbf{b})^{\top} + \mathbf{p}_1^{\top} \left( \mathbf{U} \boldsymbol{\Lambda}^{{\frac{1}{2}}}\mathbf{B} \right)^{\top} \mathbf{A}^{\top} \leq 0,  \\
&& t_2 (\mathbf{A}\boldsymbol{\mu}-\mathbf{b})^{\top} + \mathbf{p}_2^{\top}\left(\mathbf{U} \boldsymbol{\Lambda}^{{\frac{1}{2}}}\mathbf{B} \right)^{\top} \mathbf{A}^{\top}   \leq 0,  \\
&& \gamma_2 \mathbf{I}_{m_1} -  \mathbf{P}_1-\mathbf{P}_2 \succeq 0 ,\ t_1 \left(\mathbf{c} - \mathbf{v}\right) +   t_2 \left(\mathbf{c} - \mathbf{g}\right)  \geq 0, \\
&&  \begin{bmatrix}
 t_1 & \mathbf{p}_1^{\top} \\
 \mathbf{p}_1 & \mathbf{P}_1 
\end{bmatrix} \succeq 0, \ 
\begin{bmatrix}
 t_2 & \mathbf{p}_2^{\top} \\
 \mathbf{p}_2 & \mathbf{P}_2 
\end{bmatrix} \succeq 0, \ \mathbf{B}^{\top} \mathbf{B} = \mathbf{I}_{m_1}, \\
&& \mathbf{B} \in \mathbb{R}^{m \times m_1}, \ \mathbf{p}_1 \in \mathbb{R}^{m_1}, \ \mathbf{p}_2 \in \mathbb{R}^{m_1}, \ \mathbf{P}_1 \in \mathbb{R}^{m_1\times m_1}, \ \mathbf{P}_2 \in \mathbb{R}^{m_1 \times m_1}.
\end{subeqnarray}

By the inner approximation  \eqref{Equ:upperbound-setform}, the following problem provides an upper bound for the optimal value of Problem \eqref{Equ:NewsvendorSDP} and achieves the optimal value of Problem \eqref{Equ:NewsvendorSDP} when $m_1 \geq 2$:
\begin{subeqnarray} \label{Equ:NewsvendorUB}
& \min\limits_{\substack{ \mathbf{B}, \mathbf{x},s,\boldsymbol{\lambda}_1, \boldsymbol{\lambda}_2,\\ \mathbf{q},\mathbf{Q}_{\textup{r}}, \mathbf{u}_1, \mathbf{u}_2 }} & s +\gamma_2\mathbf{I}_{m_1} \bullet \mathbf{Q}_{\textup{r}} + \sqrt{\gamma_1} \left \| \mathbf{q}\right \|_2\\ 
&\text{s.t.} & \small{ \begin{bmatrix}
s-(\mathbf{c}-\mathbf{v})^{\top} \mathbf{x}-\boldsymbol{\lambda}_1^{\top}(\mathbf{b}-\mathbf{A}\boldsymbol{\mu}) &\ \frac{1}{2} \mathbf{u}_1^{\top} \\
\frac{1}{2}  \mathbf{u}_1 & \mathbf{Q}_{\textup{r}}
\end{bmatrix}} \succeq 0, \\
&& \small{\begin{bmatrix}
s-(\mathbf{c}-\mathbf{g})^{\top} \mathbf{x}-\boldsymbol{\lambda}_2^{\top}(\mathbf{b}-\mathbf{A}\boldsymbol{\mu}) + (\mathbf{v}-\mathbf{g})^{\top} \boldsymbol{\mu} &\ \frac{1}{2}  \mathbf{u}_2^{\top}\\
\frac{1}{2}  \mathbf{u}_2& \mathbf{Q}_{\textup{r}} 
\end{bmatrix}} \succeq 0,  \\ 
&& \mathbf{q} +\left(\mathbf{U}\boldsymbol{\Lambda}^{{\frac{1}{2}}}\right)^{\top}\mathbf{A}^{\top}\boldsymbol{\lambda}_1 = \mathbf{B} \mathbf{u}_1, \\ 
&& \mathbf{q} +\left(\mathbf{U}\boldsymbol{\Lambda}^{{\frac{1}{2}}}\right)^{\top}\left(\mathbf{A}^{\top}\boldsymbol{\lambda}_2+\mathbf{v}-\mathbf{g}\right) = \mathbf{B} \mathbf{u}_2, \\ 
&& \mathbf{x} \in \mathbb{R}^{m}_+, \ \boldsymbol{\lambda}_1  \in \mathbb{R}_+^{l}, \ \boldsymbol{\lambda}_2 \in \mathbb{R}_+^{l}, \ \mathbf{B}^{\top} \mathbf{B} = \mathbf{I}_{m_1}, \\ 
&& \mathbf{q} \in \mathbb{R}^m, \ \mathbf{Q}_{\textup{r}} \in \mathbb{R}^{m_1 \times m_1}, \ \mathbf{B} \in \mathbb{R}^{m \times m_1}, \ \mathbf{u}_1 \in \mathbb{R}^{m_1}, \  \mathbf{u}_2 \in \mathbb{R}^{m_1}.
\end{subeqnarray}

By the second outer approximation  \eqref{Equ:.lowerbound-sdpform-rewrite}, 
the following problem with $m_1 \leq 2$ provides another lower bound for the optimal value of Problem \eqref{Equ:NewsvendorSDP} and achieves the optimal value of Problem \eqref{Equ:NewsvendorSDP} when $m_1 = 2$:
\begin{subeqnarray} \label{Equ:NewsvendorRLB}
& \min\limits_{\substack{ \mathbf{B}, \bar{\mathbf{B}}, \mathbf{x},s,\boldsymbol{\lambda}_1, \boldsymbol{\lambda}_2,\\ \mathbf{q},\mathbf{Q}_{\textup{r}}, \mathbf{u}_1, \mathbf{u}_2,\mathbf{h}_1, \mathbf{h}_2 }} & s + \gamma_2 \mathbf{I}_{m_1} \bullet \mathbf{Q}_{\textup{r}} + \sqrt{\gamma_1} \left \| \mathbf{q} \right \|_2 \\ 
& \text{s.t.} & \small{ \begin{bmatrix}
s - (\mathbf{c} - \mathbf{v})^{\top} \mathbf{x} - \boldsymbol{\lambda}_1^{\top}(\mathbf{b}-\mathbf{A}\boldsymbol{\mu}) &\ \frac{1}{2} \mathbf{u}_1^{\top} \\
\frac{1}{2}  \mathbf{u}_1 & \mathbf{Q}_{\textup{r}}
\end{bmatrix} }  \succeq 0, \\
&& \small{ \begin{bmatrix}
s-(\mathbf{c}-\mathbf{g})^{\top} \mathbf{x}-\boldsymbol{\lambda}_2^{\top}(\mathbf{b}-\mathbf{A}\boldsymbol{\mu}) + (\mathbf{v}-\mathbf{g})^{\top} \boldsymbol{\mu} &\ \frac{1}{2}  \mathbf{u}_2^{\top}\\
\frac{1}{2}  \mathbf{u}_2& \mathbf{Q}_{\textup{r}} 
\end{bmatrix} } \succeq 0,  \\ 
&& \mathbf{q} +\left(\mathbf{U}\boldsymbol{\Lambda}^{{\frac{1}{2}}}\right)^{\top}\mathbf{A}^{\top}\boldsymbol{\lambda}_1 = \mathbf{B} \mathbf{u}_1 + \bar{\mathbf{B}} \mathbf{h}_1, \\ 
&& \mathbf{q} +\left(\mathbf{U}\boldsymbol{\Lambda}^{{\frac{1}{2}}}\right)^{\top}\left(\mathbf{A}^{\top}\boldsymbol{\lambda}_2+\mathbf{v}-\mathbf{g}\right) = \mathbf{B} \mathbf{u}_2 + \bar{\mathbf{B}} \mathbf{h}_2, \\ 
&& \mathbf{x} \in \mathbb{R}^{m}_+, \ 
\boldsymbol{\lambda}_1  \in \mathbb{R}_+^{l}, \ 
\boldsymbol{\lambda}_2 \in \mathbb{R}_+^{l}, \ 
[\mathbf{B}, \bar{\mathbf{B}}]^{\top} [\mathbf{B}, \bar{\mathbf{B}}] = \mathbf{I}_{K}, \\ 
&& \mathbf{q} \in \mathbb{R}^m, \ 
\mathbf{Q}_{\textup{r}} \in \mathbb{R}^{m_1 \times m_1}, \ 
\mathbf{B} \in \mathbb{R}^{m \times m_1}, \ 
\bar{\mathbf{B}} \in \mathbb{R}^{m \times (K-m_1)}, \\
&& \mathbf{u}_1 \in \mathbb{R}^{m_1}, \ 
\mathbf{u}_2 \in \mathbb{R}^{m_1}, \ 
\mathbf{h}_1 \in \mathbb{R}^{K-m_1}, \ \mathbf{h}_2 \in \mathbb{R}^{K-m_1}.
\end{subeqnarray}

\subsection{Production-Transportation Problem} \label{subsec-apx:production-transportation}
By Theorem 1 in \cite{cheramin2022computationally}, Problem \eqref{Equ:two-stage-dro-transportation-1} has the same optimal value as the following problem:
\begin{subeqnarray} \label{Equ:two-stage-dro-transportation-1-ref1}
\min\limits_{\mathbf{x}, s, \mathbf{q}, \mathbf{Q} } && s + \gamma_2\mathbf{I}_{mn} \bullet \mathbf{Q}  + \sqrt{\gamma_1}\left \| \mathbf{q}  \right \|_2 \\	
\textnormal{s.t.} && s \geq \mathbf{c}^{\top} \mathbf{x} +  \mathcal{U}\left(\mathcal{Q} \left(\mathbf{x}, \mathbf{U} \boldsymbol{\Lambda}^{{\frac{1}{2}}} \boldsymbol{\xi}_{\text{I}} + \boldsymbol{\mu} \right)\right)  - \boldsymbol{\xi}_{\textup{I}}^{\top} \mathbf{Q} \boldsymbol{\xi}_{\textup{I}}  - \mathbf{q}^{\top}  \boldsymbol{\xi}_{\textup{I}},  \  \forall \boldsymbol{\xi}_{\textup{I}}\in \mathcal{S}_{\textup{I}}, \slabel{Equ:two-stage-dro-transportation-1-ref1-1} \\
&&\mathbf{Q}  \succeq 0, \ \mathbf{0} \leq  \mathbf{x} \leq \mathbf{1}, \ \mathbf{Q} \in \mathbb{R}^{mn \times mn}, \ \mathbf{q} \in \mathbb{R}^{mn}.	\slabel{Equ:two-stage-dro-transportation-1-ref1-2}
\end{subeqnarray} 

\noindent 
By the definition of $\mathcal{U} (\mathcal{Q} (\mathbf{x}, \boldsymbol{\xi} ) )$, constraints \eqref{Equ:two-stage-dro-transportation-1-ref1-1} are equivalent to 
\begin{align}
s \geq \mathbf{c}^{\top} \mathbf{x} + \alpha_k \mathcal{Q} \left( \mathbf{x}, \mathbf{U} \boldsymbol{\Lambda}^{{\frac{1}{2}}} \boldsymbol{\xi}_{\text{I}}+\boldsymbol{\mu} \right) 
+ \beta_k - \boldsymbol{\xi}_{\textup{I}}^{\top} \mathbf{Q} \boldsymbol{\xi}_{\textup{I}}  - \mathbf{q}^{\top}  \boldsymbol{\xi}_{\textup{I}}, \ \forall \boldsymbol{\xi}_{\textup{I}}\in \mathcal{S}_{\textup{I}}, \ k \in [K].   \label{Equ:two-stage-dro-transportation-1-ref1-1ref}  
\end{align}

\noindent
By the definition of $\mathcal{Q}(\mathbf{x}, \boldsymbol{\xi})$, constraints \eqref{Equ:two-stage-dro-transportation-1-ref1-1ref} are further equivalent to 
\begin{subeqnarray}
&& s \geq \mathbf{c}^{\top} \mathbf{x} + \alpha_k \mathbf{z}_k^{\top} \left( \mathbf{U}\boldsymbol{\Lambda}^{{\frac{1}{2}}} \boldsymbol{\xi}_{\text{I}}+\boldsymbol{\mu} \right) 
+ \beta_k - \boldsymbol{\xi}_{\textup{I}}^{\top} \mathbf{Q} \boldsymbol{\xi}_{\textup{I}}  - \mathbf{q}^{\top}  \boldsymbol{\xi}_{\textup{I}}, \  \forall \boldsymbol{\xi}_{\textup{I}}\in \mathcal{S}_{\textup{I}}, \ k \in [K],   \slabel{Equ:two-stage-dro-transportation-1-ref1-3}  \\ 
&& \sum_{i=1}^m z_{ijk} = d_j, \ \forall j \in [n], \ k \in [K], \slabel{Equ:two-stage-dro-transportation-1-ref1-4} \\
&& \sum_{j=1}^n z_{ijk} = x_i, \ \forall i \in [m], \ k \in [K], \slabel{Equ:two-stage-dro-transportation-1-ref1-5} \\ 
&& z_{ijk} \geq 0, \ \forall i \in [m], \ j \in [n], \ k \in [K], \slabel{Equ:two-stage-dro-transportation-1-ref1-6}
\end{subeqnarray}
where $\mathbf{z}_k$ is a vector whose $((i-1)m+j)$-th element is $z_{ijk}$. Therefore, Problem \eqref{Equ:two-stage-dro-transportation-1} has the same optimal value as the following problem:
\begin{align}
\min\limits_{\mathbf{x}, \mathbf{z}_k, \forall k \in [K], s, \mathbf{q}, \mathbf{Q} }
\left\{ s + \gamma_2\mathbf{I}_{mn} \bullet \mathbf{Q}  + \sqrt{\gamma_1}\left \| \mathbf{q}  \right \|_2
\ \middle| \
\eqref{Equ:two-stage-dro-transportation-1-ref1-2}, 
\eqref{Equ:two-stage-dro-transportation-1-ref1-3} - \eqref{Equ:two-stage-dro-transportation-1-ref1-6}  
\right\}.  \label{Equ:two-stage-dro-transportation-1-ref2}
\end{align}

Note that, by Theorem 1 in \cite{cheramin2022computationally}, 
the following problem can be reformulated as Problem \eqref{Equ:two-stage-dro-transportation-1-ref2}:
\begin{subeqnarray} \label{Equ:one-stage-dro-transportation}
\min\limits_{ \mathbf{x}, \mathbf{z}_k, \forall k \in [K]}  \max\limits_{\mathbb{P}_{\text{I}} \in \mathcal{D}_{\textup{M}}}  && \mathbb{E}_{\mathbb{P}} \left[ \mathbf{c}^{\top} \mathbf{x} + \alpha_k \mathbf{z}_k^{\top} \left( \mathbf{U}\boldsymbol{\Lambda}^{{\frac{1}{2}}} \boldsymbol{\xi}_{\text{I}}+\boldsymbol{\mu} \right) + \beta_k \right] \\ 
\text{s.t.} && \eqref{Equ:two-stage-dro-transportation-1-ref1-4} - \eqref{Equ:two-stage-dro-transportation-1-ref1-6}, \ \mathbf{0} \leq  \mathbf{x} \leq \mathbf{1}.
\end{subeqnarray}

\noindent 
It implies that the two-stage DRO problem \eqref{Equ:two-stage-dro-transportation} has the same optimal value 
as the single-stage DRO problem \eqref{Equ:one-stage-dro-transportation}. 
Then, by Proposition \ref{Prop:cheramin2022}, we have the following SDP reformulation:
{\small \begin{subeqnarray} \label{Equ:one-stage-dro-transportation-sdp}
&\min\limits_{\substack{ \mathbf{x}, \mathbf{z}_k (\forall k\in [K]), \\ s,  \boldsymbol{\lambda}_k (\forall k\in [K]), \mathbf{q}, \mathbf{Q} }} & s + \gamma_2 \mathbf{I}_{mn} \bullet \mathbf{Q} + \sqrt{\gamma_1}\left \| \mathbf{q} \right \|_2 \\
& {\normalfont \text{s.t.}} & \begin{bmatrix}  s-\mathbf{c}^{\top} \mathbf{x} - \beta_k -\boldsymbol{\lambda}_k^{\top} \mathbf{b} - \alpha_k \mathbf{z}_k^{\top} \boldsymbol{\mu}+\boldsymbol{\lambda}_k^{\top}\mathbf{A}\boldsymbol{\mu} & \hspace{0.1 in}  \frac{1}{2}  \left(\mathbf{q} +\left(\mathbf{U}\boldsymbol{\Lambda}^{{\frac{1}{2}}}\right)^{\top}\left(\mathbf{A}^{\top}\boldsymbol{\lambda}_k-\alpha_k \mathbf{z}_k)\right)\right)^{\top}\\ 
 \frac{1}{2} \left( \mathbf{q} + \left( \mathbf{U}\boldsymbol{\Lambda}^{{\frac{1}{2}}} \right)^{\top} \left( \mathbf{A}^{\top}\boldsymbol{\lambda}_k-\alpha_k \mathbf{z}_k ) \right) \right) & \mathbf{Q}\end{bmatrix} \succeq 0, \nonumber \\
&& \hspace{4.25in} \forall k\in[K],   \\
&& \boldsymbol{\lambda}_k  \in \mathbb{R}_+^{l}, \ \forall k \in [K], \mathbf{q} \in \mathbb{R}^{mn}, \mathbf{Q} \in \mathbb{R}^{mn \times mn}, \\ 
&& \eqref{Equ:two-stage-dro-transportation-1-ref1-4} - \eqref{Equ:two-stage-dro-transportation-1-ref1-6}, \ \mathbf{0} \leq  \mathbf{x} \leq \mathbf{1}. 
\end{subeqnarray}}%

By the first outer approximation \eqref{Equ:bilinear}, the following problem provides a lower bound for the optimal value of Problem \eqref{Equ:one-stage-dro-transportation-sdp}: 
{\small \begin{subeqnarray} \label{Equ:transportation-dual}
&\max\limits_{\substack{ t_k, \mathbf{p}_k,\mathbf{P}_k, \forall k \in [K], \\ 
\mathbf{w}_k, \mathbf{u}_k, \forall k \in [K], \mathbf{v}, \mathbf{B} }} &   \sum_{k=1}^K t_k \beta_k - \sum_{i=i}^m v_i + \sum_{k=1}^K \sum_{j=1}^n w_{jk} d_j  \\
&\text{s.t.} & 1 - \sum_{k=1}^K t_k = 0, \ \sqrt{\gamma_1} - \left\| \sum_{k=1}^K \mathbf{p}_k \right\|_2 \geq 0, \\
&& t_k (\mathbf{A}\boldsymbol{\mu}-\mathbf{b})^{\top} + \mathbf{p}_k^{\top} \left( \mathbf{U} \boldsymbol{\Lambda}^{{\frac{1}{2}}}\mathbf{B} \right)^{\top} \mathbf{A}^{\top} \leq 0, \ \forall k \in [K], \\
&& \gamma_2 \mathbf{I}_{m_1} -  \sum_{k=1}^K \mathbf{P}_k \succeq 0 ,\ \sum_{k=1}^K \left( t_k \mathbf{c} + \mathbf{u}_k \right) + \mathbf{v} \geq 0, \\
&& \alpha_k t_k \boldsymbol{\mu}^{\top} + \alpha_k \mathbf{p}_k^{\top} \left( \mathbf{U} \boldsymbol{\Lambda}^{{\frac{1}{2}}}\mathbf{B} \right)^{\top} - ( \underbrace{\mathbf{w}_k^{\top}, \ldots, \mathbf{w}_k^{\top}}_\text{repeat m times} ) - ( \underbrace{u_{1k}, \ldots, u_{1k}}_\text{repeat n times}, \ldots, u_{mk}, \ldots, u_{mk} ) \geq 0, \nonumber  \\
&& \hspace{11cm} \forall k \in [K], \\ 
&&  \begin{bmatrix}
 t_k & \mathbf{p}_k^{\top} \\
 \mathbf{p}_k & \mathbf{P}_k 
\end{bmatrix} \succeq 0, \ \forall k \in [K], \ 
\mathbf{B}^{\top} \mathbf{B} = \mathbf{I}_{m_1}, \ \mathbf{v} \in \mathbb{R}^{m}_{+}, \ \mathbf{B} \in \mathbb{R}^{mn \times m_1},  \\
&& \mathbf{p}_k \in \mathbb{R}^{m_1}, \ \mathbf{P}_k \in \mathbb{R}^{m_1\times m_1}, \ \mathbf{w}_k \in \mathbb{R}^{n}, \ \mathbf{u}_k \in \mathbf{R}^{m}, \ \forall k \in [K]. 
\end{subeqnarray}}%

By the inner approximation  \eqref{Equ:upperbound-setform}, the following problem provides an upper bound 
for the optimal value of Problem \eqref{Equ:one-stage-dro-transportation-sdp} and achieves the optimal 
value of Problem \eqref{Equ:one-stage-dro-transportation-sdp} when $m_1 \geq K$:
\begin{subeqnarray} \label{Equ:transportation-upperbound}
&\min\limits_{\substack{ \mathbf{B}, \mathbf{x}, \mathbf{z}_k (\forall k\in [K]), \\ \mathbf{u}_k,  \boldsymbol{\lambda}_k, \forall k\in [K], s, \mathbf{q}, \mathbf{Q}_\text{r} }} & s + \gamma_2 \mathbf{I}_{m_1} \bullet \mathbf{Q}_{\textup{r}} + \sqrt{\gamma_1}\left \| \mathbf{q} \right \|_2 \\
& {\normalfont \text{s.t.}} & \begin{bmatrix}  s-\mathbf{c}^{\top} \mathbf{x} - \beta_k -\boldsymbol{\lambda}_k^{\top}\mathbf{b}-\alpha_k \mathbf{z}_k^{\top}\boldsymbol{\mu}+\boldsymbol{\lambda}_k^{\top}\mathbf{A}\boldsymbol{\mu} & \hspace{0.1 in}  \frac{1}{2}  \mathbf{u}_k^{\top}\\ 
 \frac{1}{2} \mathbf{u}_k & \mathbf{Q}_{\textup{r}}  \end{bmatrix} \succeq 0, \ \forall k\in[K],   \\
&& \mathbf{q} +\left(\mathbf{U}\boldsymbol{\Lambda}^{{\frac{1}{2}}}\right)^{\top}\left(\mathbf{A}^{\top}\boldsymbol{\lambda}_k-\alpha_k \mathbf{z}_k)\right) = \mathbf{B} \mathbf{u}_k, \ \forall k \in [K], \\
&& \boldsymbol{\lambda}_k  \in \mathbb{R}_+^{l}, \ \mathbf{u}_k \in \mathbb{R}^{m_1}, \ \forall k \in [K], \mathbf{q} \in \mathbb{R}^{mn}, \mathbf{Q}_{\textup{r}} \in \mathbb{R}^{m_1 \times m_1}, \\ 
&& \eqref{Equ:two-stage-dro-transportation-1-ref1-4} - \eqref{Equ:two-stage-dro-transportation-1-ref1-6}, \ \mathbf{0} \leq  \mathbf{x} \leq \mathbf{1},\ \mathbf{B}^{\top} \mathbf{B} = \mathbf{I}_{m_1}, \ \mathbf{B} \in \mathbb{R}^{mn \times m_1}.  
\end{subeqnarray}

By the second outer approximation  \eqref{Equ:.lowerbound-sdpform-rewrite}, 
the following problem with $m_1 \leq K$ provides another lower bound for the optimal value of Problem 
\eqref{Equ:one-stage-dro-transportation-sdp} and achieves the optimal value of Problem \eqref{Equ:one-stage-dro-transportation-sdp} when $m_1 = K$:
\begin{subeqnarray} \label{Equ:transportation-LBR}
&\min\limits_{\substack{ \mathbf{B}, \bar{\mathbf{B}}, \mathbf{x}, \mathbf{z}_k (\forall k\in [K]), \\ \mathbf{u}_k, \mathbf{h}_k,  \boldsymbol{\lambda}_k, \forall k\in [K], s, \mathbf{q}, \mathbf{Q}_\text{r} }} & s + \gamma_2 \mathbf{I}_{m_1} \bullet \mathbf{Q}_{\textup{r}} + \sqrt{\gamma_1}\left \| \mathbf{q} \right \|_2 \\
& {\normalfont \text{s.t.}} & \begin{bmatrix}  s-\mathbf{c}^{\top} \mathbf{x} - \beta_k -\boldsymbol{\lambda}_k^{\top}\mathbf{b}-\alpha_k \mathbf{z}_k^{\top}\boldsymbol{\mu}+\boldsymbol{\lambda}_k^{\top}\mathbf{A}\boldsymbol{\mu} & \hspace{0.1 in}  \frac{1}{2}  \mathbf{u}_k^{\top}\\ 
 \frac{1}{2} \mathbf{u}_k & \mathbf{Q}_{\textup{r}}  \end{bmatrix} \succeq 0, \ \forall k\in[K],   \\
&& \mathbf{q} +\left(\mathbf{U}\boldsymbol{\Lambda}^{{\frac{1}{2}}}\right)^{\top}\left(\mathbf{A}^{\top}\boldsymbol{\lambda}_k-\alpha_k \mathbf{z}_k)\right) = \mathbf{B} \mathbf{u}_k + \bar{\mathbf{B}} \mathbf{h}_k, \ \forall k \in [K], \\
&& \boldsymbol{\lambda}_k  \in \mathbb{R}_+^{l}, \ \mathbf{u}_k \in \mathbb{R}^{m_1}, \ \mathbf{h}_k \in \mathbb{R}^{K-m_1}, \ \forall k \in [K], \mathbf{q} \in \mathbb{R}^{mn}, \mathbf{Q}_{\textup{r}} \in \mathbb{R}^{m_1 \times m_1}, \\ 
&& \eqref{Equ:two-stage-dro-transportation-1-ref1-4} - \eqref{Equ:two-stage-dro-transportation-1-ref1-6}, \ \mathbf{0} \leq  \mathbf{x} \leq \mathbf{1},\ [\mathbf{B}, \bar{\mathbf{B}}]^{\top} [\mathbf{B}, \bar{\mathbf{B}}] = \mathbf{I}_{m_1}, \ \mathbf{B} \in \mathbb{R}^{mn \times m_1}, \ \bar{\mathbf{B}} \in \mathbb{R}^{mn \times (K-m_1)}.  
\end{subeqnarray}

\subsection{Proof of Proposition \ref{prop:solve-int-dim-opt}}
Because $\mathbf{B}^{\top} \mathbf{B} = \mathbf{I}_{m_1}$, we have $\mathbf{B}^{\top} \mathbf{B} \preceq \mathbf{I}_{m_1}$, which implies $\mathbf{B} \mathbf{B}^{\top} \preceq \mathbf{I}_m$ by Lemma \ref{Lem:3cons}. 
It follows that $\mathbf{r}^{\top} \mathbf{B} \mathbf{B}^{\top} \mathbf{r} \leq \mathbf{r}^{\top} \mathbf{r}$. 
Meanwhile, we have 
\begin{align*}
\mathbf{r}^{\top} \mathbf{B}^* \mathbf{B}^{*\top} \mathbf{r} = 
\begin{bmatrix}
 \frac{\mathbf{r}^{\top} \mathbf{r}}{\| \mathbf{r} \|_2} & \mathbf{0}_{1 \times (m_1-1)}
\end{bmatrix} 
\begin{bmatrix}
 \frac{\mathbf{r}^{\top} \mathbf{r}}{\| \mathbf{r} \|_2} \\ \mathbf{0}_{ (m_1-1)\times 1} 
\end{bmatrix} = \mathbf{r}^{\top} \mathbf{r},
\end{align*}
indicating that $\mathbf{B}^* = \begin{bmatrix}
  \mathbf{r}/ \| \mathbf{r} \|_2, \mathbf{0}_{m \times (m_1-1)}
\end{bmatrix}$ is an optimal solution of Problem \eqref{Equ:insight}.
\Halmos

\subsection{Sensitivity Analyses} \label{subsec:sensi_analysis}

\begin{table}[!htbp]   
\centering
\caption{Sensitivity Analyses on the Production-Transportation Problem with $K=5$}\label{Table:production-transportation-sensiAna-K5}
\centering
\resizebox{0.92\textwidth}{!}{  
\begin{tabular}{|c|cc|c|c|c|c|c|c|c|} 
\hline\rule{0pt}{2.5ex} 
&& Size ($(m,n)$) & (4,25) & (5,20) & (5,40) & (8,25) & (10,40) & (20,30) & (20,40)  \\ [0.2em] 
\hline\rule{0pt}{2.5ex}
\multirow{11}{*}{$m_1 = 3$} & \multirow{4}{*}{ODR-LB} & Gap1 ($\%$) &0.37  & 0.71  & 0.67 & 0.47 & - & - & - \\[0.2em]
&& Time (secs) &2.86  & 2.81  & 4.87  & 4.32  & 15.16  & 27.51  & 57.45  \\[0.2em]
&& Interval Gap ($\%$) &0.69  & 0.79  & 0.73  & 0.51  & 0.49  & 0.51  & 0.57  \\[0.2em] 
&& Theoretical Gap ($\%$) & 1.27 & 0.62 & 2.20 & 1.36 & - & - & - \\ [0.2em]
\cline{2-10}  \rule{0pt}{2.5ex}
&\multirow{4}{*}{ODR-RLB} & Gap1 ($\%$) &0.07  & 0.17  & 0.02 & 0.02 & - & - & -  \\[0.2em]
&& Time (secs) & 13.53  & 14.19  & 43.93  & 57.79  & 205.99  & 442.26  & 1092.43  \\[0.2em]
&& Interval Gap ($\%$) &0.39  & 0.25  & 0.08  & 0.06  & 0.03  & 0.02  & 0.02  \\[0.2em]
&& Theoretical Gap ($\%$) & 1.08 & 3.19 & 1.96 & 0.96 & - & - & -  \\ [0.2em]
\cline{2-10}  \rule{0pt}{2.5ex}
&\multirow{3}{*}{ODR-UB} & Gap2 ($\%$) &0.32  & 0.08  & 0.06 & 0.04 & - & - & -  \\[0.2em]
&& Time (secs) & 6.03  & 5.60  & 22.31  & 20.53  & 122.04  & 332.10  & 660.08  \\[0.2em]
&& Theoretical Gap ($\%$) & 1.28 & 1.03 & 3.20 & 2.05 & - & - & -  \\ [0.2em]
\hline\rule{0pt}{2.5ex}
\multirow{11}{*}{$m_1 = 5$} & \multirow{4}{*}{ODR-LB} & Gap1 ($\%$) &0.14  & 0.34  & 0.22 & 0.29 & - & - & - \\[0.2em] 
&& Time (secs) &4.03  & 3.63  & 6.01  & 4.82  & 12.41  & 25.03  & 57.79  \\[0.2em]
&& Interval Gap ($\%$) &0.15  & 0.34  & 0.43  & 0.29  & 0.52  & 0.55  & 0.51   \\[0.2em]
&& Theoretical Gap ($\%$) & 1.24 & 1.07 & 4.91 & 0.77 & 2.71 & 1.25 & 1.09  \\ [0.2em]
\cline{2-10}  \rule{0pt}{2.5ex}
&\multirow{4}{*}{ODR-RLB} & Gap1 ($\%$) &0.01  & 0.02  & 0.01 & 0.00 & - & - & -  \\[0.2em]
&& Time (secs) & 5.42  & 5.36  & 22.40  & 20.86  & 123.64  & 330.29  & 665.43   \\[0.2em]
&& Interval Gap ($\%$) &0.02  & 0.02  & 0.01  & 0.01  & 0.00  & 0.00  & 0.00   \\[0.2em]
&& Theoretical Gap ($\%$) & 1.13 & 1.17 & 1.80 & 0.81 & 1.22 & 1.92 & 1.32  \\ [0.2em]
\cline{2-10}  \rule{0pt}{2.5ex}
&\multirow{3}{*}{ODR-UB} & Gap2 ($\%$) &0.01  & 0.01  & 0.00 & 0.00 & - & - & -  \\[0.2em]
&& Time (secs) & 5.48  & 5.32  & 22.41  & 20.86  & 123.48  & 329.95  & 665.35   \\[0.2em]
&& Theoretical Gap ($\%$) & 1.13 & 1.17 & 1.80 & 0.81 & 1.22 & 1.92 & 1.32 \\ [0.2em]
\hline\rule{0pt}{2.5ex}
\multirow{7}{*}{$m_1 = 7$} & \multirow{4}{*}{ODR-LB} & Gap1 ($\%$) &0.16  & 0.21  & 0.53 & 0.23 & - & - & -  \\[0.2em]
&& Time (secs) &3.56  & 4.29  & 5.39  & 5.70  & 19.27  & 22.89  & 58.38  \\[0.2em]
&& Interval Gap ($\%$) &0.17  & 0.22  & 0.53  & 0.24  & 0.31  & 0.56  & 0.50  \\[0.2em]
&& Theoretical Gap ($\%$) & 1.22 & 0.88 & 1.67 & 0.79 & - & - & - \\ [0.2em]
\cline{2-10}  \rule{0pt}{2.5ex}
&\multirow{3}{*}{ODR-UB} & Gap2 ($\%$) &0.02  & 0.01  & 0.00 & 0.00 & - & - & -  \\[0.2em]
&& Time (secs) & 5.62  & 5.48  & 22.42  & 21.24  & 121.61  & 330.94  & 660.85  \\[0.2em]
&& Theoretical Gap ($\%$) & 1.15 & 1.04 & 1.28 & 0.59 & - & - & -  \\ [0.2em]  
\hline 
\end{tabular}} 
\end{table}

\begin{table}[!htbp]   
\centering
\caption{Sensitivity Analyses on the Production-Transportation Problem with $K=10$}\label{Table:production-transportation-sensiAna-K10}
\centering
\resizebox{0.92\textwidth}{!}{  
\begin{tabular}{|c|cc|c|c|c|c|c|c|c|} 
\hline\rule{0pt}{2.5ex} 
&& Size ($(m,n)$) & (4,25) & (5,20) & (5,40) & (8,25) & (10,40) & (20,30) & (20,40)  \\ [0.2em] 
\hline\rule{0pt}{2.5ex}
\multirow{11}{*}{$m_1 = 8$} & \multirow{4}{*}{ODR-LB} & Gap1 ($\%$) &0.17  & 0.15  & 0.16 & 0.28 & - & - & - \\[0.2em] 
&& Time (secs) &6.26  & 5.75  & 9.80  & 8.87  & 29.94  & 54.71  & 133.58  \\[0.2em]
&& Interval Gap ($\%$) &0.17  & 0.15  & 0.16  & 0.28  & 0.23  & 0.22  & 0.27  \\[0.2em] 
&& Theoretical Gap ($\%$) & 4.32 & 4.53 & 6.22 & 6.09 & - & - & - \\ [0.2em]
\cline{2-10}  \rule{0pt}{2.5ex}
&\multirow{4}{*}{ODR-RLB} & Gap1 ($\%$) &0.01  & 0.00  & 0.01 & 0.00 & - & - & -  \\[0.2em]
&& Time (secs) & 13.03  & 12.91  & 53.84  & 35.92  & 205.07  & 593.38  & 1334.02  \\[0.2em]
&& Interval Gap ($\%$) &0.01  & 0.00  & 0.01  & 0.01  & 0.00  & 0.00  & 0.00  \\[0.2em]
&& Theoretical Gap ($\%$) & 3.53 & 3.02 & 4.24 & 4.86 & - & - & -  \\ [0.2em]
\cline{2-10}  \rule{0pt}{2.5ex}
&\multirow{3}{*}{ODR-UB} & Gap2 ($\%$) &0.00  & 0.00  & 0.00 & 0.00 & - & - & -  \\[0.2em]
&& Time (secs) & 10.31  & 10.23  & 39.26  & 31.59  & 119.58  & 339.84  & 755.18  \\[0.2em]
&& Theoretical Gap ($\%$) & 2.21 & 2.40 & 3.98 & 3.04 & - & - & -  \\ [0.2em] 
\hline\rule{0pt}{2.5ex}
\multirow{11}{*}{$m_1 = 10$} & \multirow{4}{*}{ODR-LB} & Gap1 ($\%$) & 0.12  & 0.19  & 0.13 & 0.25 & - & - & - \\[0.2em] 
&& Time (secs) & 6.65  & 6.41  & 11.11 & 10.31 & 31.42 & 61.41 & 131.10 \\[0.2em]
&& Interval Gap ($\%$) & 0.12  & 0.19  & 0.13 & 0.25 & 0.22 & 0.19 & 0.25 \\[0.2em]
&& Theoretical Gap ($\%$) & 1.78  & 1.97  & 2.96 & 1.19 & 2.10 & 1.49 & 1.22 \\[0.2em]
\cline{2-10}\rule{0pt}{2.5ex}
&\multirow{4}{*}{ODR-RLB} & Gap1 ($\%$) & 0.00  & 0.00  & 0.00 & 0.00 & - & - & - \\[0.2em]
&& Time (secs) & 11.22  & 10.65  & 40.34 & 32.95 & 122.54 & 342.80 & 748.11 \\[0.2em]
&& Interval Gap ($\%$) & 0.00  & 0.00  & 0.00 & 0.00 & 0.00 & 0.00 & 0.00 \\[0.2em]
&& Theoretical Gap ($\%$) & 1.86  & 2.27  & 1.74 & 1.22 & 2.10 & 1.39 & 2.30 \\[0.2em]
\cline{2-10}\rule{0pt}{2.5ex} 
&\multirow{3}{*}{ODR-UB} & Gap2 ($\%$) & 0.00  & 0.00  & 0.00 & 0.00 & - & - & - \\[0.2em]
&& Time (secs) & 11.14  & 10.69  & 40.28 & 32.92 & 122.48 & 344.10 & 747.90 \\[0.2em]
&& Theoretical Gap ($\%$) & 1.86  & 2.27  & 1.74 & 1.22 & 2.10 & 1.39 & 2.30\\[0.2em]
\hline\rule{0pt}{2.5ex}
\multirow{7}{*}{$m_1 = 12$} & \multirow{4}{*}{ODR-LB} & Gap1 ($\%$) &0.11  & 0.14  & 0.15 & 0.23 & - & - & -  \\[0.2em] 
&& Time (secs) &8.89  & 7.90  & 11.78  & 13.44  & 32.16  & 63.96  & 122.91  \\[0.2em]
&& Interval Gap ($\%$) &0.11  & 0.14  & 0.15  & 0.23  & 0.19  & 0.19  & 0.27  \\[0.2em]
&& Theoretical Gap ($\%$) & 3.64 & 5.01 & 7.41 & 7.85 & - & - & - \\ [0.2em]
\cline{2-10}  \rule{0pt}{2.5ex}
&\multirow{3}{*}{ODR-UB} & Gap2 ($\%$) &0.00  & 0.00  & 0.00 & 0.00 & - & - & -  \\[0.2em]
&& Time (secs) & 11.71  & 11.55  & 38.59  & 29.68  & 123.40  & 346.68  & 761.21  \\[0.2em]
&& Theoretical Gap ($\%$) & 0.81 & 1.32 & 2.20 & 2.75 & - & - & -  \\ [0.2em]  
\hline 
\end{tabular}} 
\end{table}

\begin{table}[!htbp]    
\centering
\caption{Sensitivity Analyses on the Production-Transportation Problem with $K=15$}\label{Table:production-transportation-sensiAna-K15}
\centering
\resizebox{0.92\textwidth}{!}{  
\begin{tabular}{|c|cc|c|c|c|c|c|c|c|} 
\hline\rule{0pt}{2.5ex} 
&& Size ($(m,n)$) & (4,25) & (5,20) & (5,40) & (8,25) & (10,40) & (20,30) & (20,40)  \\ [0.2em] 
\hline\rule{0pt}{2.5ex}
\multirow{11}{*}{$m_1 = 13$} & \multirow{4}{*}{ODR-LB} & Gap1 ($\%$) &0.08  & 0.21  & 0.15 & 0.18 & - & - & - \\[0.2em]
&& Time (secs) &11.12  & 11.10  & 18.00  & 24.05  & 31.32  & 71.04  & 142.88  \\[0.2em]
&& Interval Gap ($\%$) &0.08  & 0.21  & 0.15  & 0.18  & 0.25  & 0.19  & 0.23  \\[0.2em] 
&& Theoretical Gap ($\%$) & 4.56 & 7.01 & 5.79 & 6.95 & - & - & - \\ [0.2em]
\cline{2-10}  \rule{0pt}{2.5ex}
&\multirow{4}{*}{ODR-RLB} & Gap1 ($\%$) &0.00  & 0.03  & 0.00 & 0.00 & - & - & -  \\[0.2em]
&& Time (secs) & 26.19  & 25.91  & 97.35  & 68.99  & 213.44  & 641.74  & 1451.80  \\[0.2em]
&& Interval Gap ($\%$) &0.00  & 0.03  & 0.00  & 0.00  & 0.00  & 0.01  & 0.00  \\[0.2em]
&& Theoretical Gap ($\%$) & 6.09 & 5.83 & 5.88 & 6.12 & - & - & -  \\ [0.2em]
\cline{2-10}  \rule{0pt}{2.5ex}
&\multirow{3}{*}{ODR-UB} & Gap2 ($\%$) &0.00  & 0.00  & 0.00 & 0.00 & - & - & -  \\[0.2em]
&& Time (secs) & 19.46  & 19.15  & 71.71  & 50.62  & 128.92  & 367.33  & 833.81  \\[0.2em]
&& Theoretical Gap ($\%$) & 2.44 & 3.00 & 2.52 & 3.07 & - & - & -  \\ [0.2em]
\hline\rule{0pt}{2.5ex}
\multirow{11}{*}{$m_1 = 15$} & \multirow{4}{*}{ODR-LB} & Gap1 ($\%$) &0.05  & 0.10  & 0.12 & 0.15 & - & - & - \\[0.2em] 
&& Time (secs) & 13.92  & 16.56  & 22.11  & 26.16  & 43.02  & 80.51  & 168.90  \\[0.2em]
&& Interval Gap ($\%$) &0.05  & 0.10  & 0.12  & 0.15  & 0.11  & 0.13  & 0.09   \\[0.2em]
&& Theoretical Gap ($\%$) & 4.69 & 5.22 & 5.57 & 6.38 & 3.29 & 4.41 & 3.95  \\ [0.2em]
\cline{2-10}  \rule{0pt}{2.5ex}
&\multirow{4}{*}{ODR-RLB} & Gap1 ($\%$) &0.00  & 0.00  & 0.00 & 0.00 & - & - & -  \\[0.2em]
&& Time (secs) & 22.60  & 21.45  & 77.18  & 63.99  & 241.03  & 689.24  & 1550.21   \\[0.2em]
&& Interval Gap ($\%$) &0.00  & 0.00  & 0.00  & 0.00  & 0.00  & 0.00  & 0.00   \\[0.2em]
&& Theoretical Gap ($\%$) & 1.38 & 2.61 & 1.79 & 1.73 & 2.14 & 1.94 & 2.25  \\ [0.2em]
\cline{2-10}  \rule{0pt}{2.5ex}
&\multirow{3}{*}{ODR-UB} & Gap2 ($\%$) &0.00  & 0.00  & 0.00 & 0.00 & - & - & -  \\[0.2em]
&& Time (secs) & 22.61  & 21.44  & 76.99  & 64.08  & 149.21  & 401.63  & 878.68   \\[0.2em]
&& Theoretical Gap ($\%$) & 1.38 & 2.61 & 1.79 & 1.73 & 2.14 & 1.94 & 2.25  \\ [0.2em]
\hline\rule{0pt}{2.5ex}
\multirow{7}{*}{$m_1 = 17$} & \multirow{4}{*}{ODR-LB} & Gap1 ($\%$) &0.10  & 0.09  & 0.13 & 0.16 & - & - & -  \\[0.2em]
&& Time (secs) &18.10  & 21.04  & 28.81  & 32.61  & 52.71  & 84.10  & 188.21  \\[0.2em]
&& Interval Gap ($\%$) &0.10  & 0.09  & 0.13  & 0.16  & 0.10  & 0.08  & 0.13  \\[0.2em]
&& Theoretical Gap ($\%$) & 4.75 & 4.13 & 5.50 & 6.16 & - & - & - \\ [0.2em]
\cline{2-10}  \rule{0pt}{2.5ex}
&\multirow{3}{*}{ODR-UB} & Gap2 ($\%$) &0.00  & 0.00  & 0.00 & 0.00 & - & - & -  \\[0.2em]
&& Time (secs) & 26.83  & 27.84  & 84.70  & 63.17  & 161.79  & 447.03  & 926.33  \\[0.2em]
&& Theoretical Gap ($\%$) & 2.20 & 3.61 & 2.58 & 2.20 & - & - & -  \\ [0.2em]  
\hline 
\end{tabular}} 
\end{table} 


\end{APPENDICES}

\end{document}